\documentclass[10pt]{article}
\usepackage{amssymb}
\usepackage{amsthm}
\usepackage{enumerate}
\usepackage{amsmath}
\usepackage{mathtools}
\usepackage{extarrows}
\usepackage{fullpage}
\usepackage{xcolor}
\usepackage{hyperref}
\usepackage{tikz}
\usetikzlibrary{decorations.markings, arrows.meta}
\usepackage{caption}
\usepackage{float}

\newtheorem{thm}{Theorem}[section]
\newtheorem{lm}[thm]{Lemma}
\newtheorem{defn}[thm]{Definition}
\newtheorem{prop}[thm]{Proposition}
\newtheorem{cor}[thm]{Corollary}

\newtheorem{rmk}[thm]{Remark}

\numberwithin{equation}{section}
\newcommand{\nn}{n}
\newcommand{\kk}{k}
\newcommand{\mm}{m}
\newcommand{\rr}{r}

\newcommand{\Airy}{\mathcal{A}}

\newcommand{\bbeta}{\boldsymbol{\beta}}

\newcommand{\bn}{{\boldsymbol{\nn}}}

\newcommand{\bk}{{\boldsymbol{\kk}}}

\newcommand{\bxi}{\boldsymbol{\xi}}
\newcommand{\bzeta}{\boldsymbol{\zeta}}
\newcommand{\cC}{\mathcal{C}}
\newcommand{\complexC}{\mathbb{C}}

\newcommand{\dL}{\mathcal{L}}
\newcommand{\fdd}{\mathrm{f.d.d.}}

\newcommand{\hKPZ}{\mathcal{H}^{\mathrm{KPZ}}}
\newcommand{\hKPZhB}{\hKPZ(\alpha,\tau;\mathfrak{h}_{\mathrm{hB}})}

\newcommand{\identity}{\mathbf{1}}

\newcommand{\inn}{\mathrm{in}}

\newcommand{\out}{\mathrm{out}}
\newcommand{\prob}{\mathbb{P}}

\newcommand{\rB}{\mathrm{B}}
\newcommand{\rd}{\mathrm{d}}
\newcommand{\realR}{\mathbb{R}}

\newcommand{\rmi}{\mathrm{i}}

\newcommand{\rL}{\mathrm{L}}

\newcommand{\rR}{\mathrm{R}}
\newcommand{\rT}{\mathrm{T}}
\newcommand{\rK}{\mathrm{K}}

\newcommand{\Var}{\mathrm{Var}}
\newcommand{\vecteta}{\boldsymbol{\eta}}
\renewcommand{\Re}{\mathrm{Re}}
\renewcommand{\Im}{\mathrm{Im}}
\newcommand{\hatD}{\hat{D}}
\newcommand{\hatK}{\widehat{\rK}}
\newcommand{\Res}{\operatorname*{Res}}

\title{Edge behavior of the KPZ fixed point under upper-tail conditioning}
\author{Yuanyuan Liu\footnote{School of Mathematics and Statistics, Hunan Research Center of the Basic Discipline for
Analytical Mathematics, HNP-LAMA, Central South University. Email: \texttt{liuyy@csu.edu.cn}} \and Zhipeng Liu\footnote{Department of Mathematics and Statistics, University of Maine. Email: \texttt{zhipeng.liu@maine.edu} } \and Haoyang Ying\footnote{School of Mathematics and Statistics, Hunan Research Center of the Basic Discipline for
Analytical Mathematics, HNP-LAMA, Central South University. Email: \texttt{haoyangying@csu.edu.cn}}}
\begin{document}

\hypersetup{hidelinks}
\maketitle

\begin{abstract}
We consider the KPZ fixed point $\hKPZ(x,t)$ with the narrow-wedge initial profile, conditioned on $\{\hKPZ(0,1)\ge L\}$ as $L\to\infty$. Our main result describes its fluctuations near the boundary of the region affected by the conditioning. More precisely, in the spatial window $x=t\sqrt L+O(1)$, the centered field converges in finite-dimensional distributions to the KPZ fixed point with a half-Brownian initial profile. The convergence holds jointly at arbitrary positive times, including times before and after the conditioning time. The proof is based on an asymptotic correspondence between conditional narrow-wedge multipoint formulas and those for the half-Brownian KPZ fixed point.
\end{abstract}
\section{Introduction and main results}\label{sec:introduction}

\subsection{The KPZ fixed point and upper-tail conditioning}\label{sec:background}

The KPZ fixed point is a space-time random field expected to describe the limiting height fluctuations of models in the Kardar--Parisi--Zhang (KPZ) universality class. These include random growth models, interacting particle systems, and directed polymers, whose height fluctuations and spatial correlation lengths have scales $t^{1/3}$ and $t^{2/3}$, respectively \cite{Corwin12,Baik23,Remenik23,Quastel26}. The KPZ fixed point was constructed by Matetski, Quastel, and Remenik \cite{Matetski-Quastel-Remenik21} as the limit of the totally asymmetric simple exclusion process (TASEP) with general initial profiles under $1:2:3$ scaling. It is a Markovian process with transition probabilities given by explicit Fredholm determinant formulas.

The KPZ fixed point can also be represented through the directed landscape, the random metric constructed by Dauvergne, Ortmann, and Vir\'ag \cite{Dauvergne-Ortmann-Virag22} as the limit of Brownian last-passage percolation. Let $\dL$ denote the directed landscape. For an admissible initial profile $\mathfrak{h}$, this representation is \cite{Nica-Quastel-Remenik20,Dauvergne-Virag21}
\begin{equation*}
  \hKPZ(x,t;\mathfrak{h})=\sup_{y\in\realR}\{\mathfrak{h}(y)+\dL(y,0;x,t)\}.
\end{equation*}
Convergence to the KPZ fixed point or the directed landscape is now known for a number of models \cite{Nica-Quastel-Remenik20,Virag20,Dauvergne-Virag21,Wu26,Aggarwal-Corwin-Hegde24,Dauvergne-Zhang24,Aggarwal-Corwin-Hegde24a,Zhang25,Zhang26}. Related limiting objects have also been studied in half-space and periodic settings \cite{Zhang24,Dauvergne-Zhang26,Baik-Liao-Liu26,Aggarwal-Corwin-Schmid26}.

In this paper, we study the KPZ fixed point conditioned on a large height at one point. Let $\mathfrak{h}_{\mathrm{NW}}$ be the narrow-wedge initial profile, defined by $\mathfrak{h}_{\mathrm{NW}}(0)=0$ and $\mathfrak{h}_{\mathrm{NW}}(x)=-\infty$ for $x\ne0$. We write $\hKPZ(x,t)=\hKPZ(x,t;\mathfrak{h}_{\mathrm{NW}})$, so that $\hKPZ(x,t)=\dL(0,0;x,t)$. More precisely, we consider the conditional law given $\{\hKPZ(0,1)\ge L\}$ as $L\to\infty$, and investigate the fluctuations near the boundary of the region affected by this conditioning.

There have been several studies of upper-tail conditioning. In the two-time setting, the regime in which the earlier height is large and positive was analyzed using the replica Bethe ansatz in \cite{de_Nardis-Le_Doussal17}. The resulting tail formula was subsequently shown to agree with the corresponding asymptotics of the rigorous two-time formula \cite{de_Nardis-Le_Doussal18}. At the macroscopic scale, sharp upper-tail estimates and conditional spatial profiles were obtained in \cite{Ganguly-Hegde22}, and multipoint upper-tail large deviations and space-time limit shapes for the narrow-wedge KPZ equation were studied in \cite{Lin-Tsai25}. For the directed landscape, an upper-tail large deviation principle was established in \cite{Das-Dauvergne-Virag26}, and narrow-wedge conditional limit shapes were obtained in \cite{Das-Tsai24}.

The conditional fluctuations depend on the space-time scale. Near the segment joining the origin to the conditioned point, the limiting height fluctuations are described by two independent Brownian bridges \cite{Liu-Wang24}. Related limits hold for the periodic KPZ fixed point \cite{Baik-Liu24}. In a smaller neighborhood of the conditioned point, a different limiting field, called the upper tail field, was obtained in \cite{Liu-Zhang25}. At a later macroscopic time, a GUE Tracy--Widom one-point limit and an Edgeworth-type expansion were proved in \cite{Nissim-Zhang26}.

Upper-tail conditioning also changes the geometry of geodesics. In exponential last-passage percolation, it was proved in \cite{Basu-Ganguly23} that upper-tail conditioning changes the transversal fluctuation exponent from $2/3$ to $1/2$. Further results for conditioned geodesics include rigidity and Gaussian one-point fluctuations \cite{Liu22b}, Brownian-bridge limits for geodesics and polymer path measures \cite{Ganguly-Hegde-Zhang23}, and a joint one-point limit near the endpoints \cite{Liu-Ma-Tripathi26}.

For exponential directed last-passage percolation, a conditional law of large numbers and fluctuation limits in certain regimes were obtained in \cite{Baik-Cordaro-Tripathi25}. More recently, one-point fluctuation laws were determined in two of the three spatial regions and at their transition boundaries, giving Gaussian, GUE Tracy--Widom, or one-spike BBP limits \cite{Baik-Tripathi26}.

\begin{figure}[htbp]
\centering
\begin{tikzpicture}[x=2.05cm,y=3.15cm,font=\small,>=Stealth]
  \definecolor{utblue}{RGB}{40,112,170}
  \definecolor{utpurple}{RGB}{127,70,155}
  \definecolor{utcyan}{RGB}{0,125,130}
  \definecolor{utedge}{RGB}{210,85,25}
  \fill[black!5] (-1.8,0) rectangle (1.8,1);
  \fill[black!17] (0,0)--(-1,1)--(1,1)--cycle;
  \draw[black!45,dashed] (0,0)--(-1,1) (0,0)--(1,1);
  \draw[black!40,densely dashed] (-1.8,1)--(1.8,1);
  \fill[utblue!35] (-.085,.08) rectangle (.085,.93);
  \draw[utblue,thick] (-.085,.08)--(-.085,.93) (.085,.08)--(.085,.93);
  \fill[utedge!25] (.07,.12)--(.17,.12)--(1.70,1.65)--(1.60,1.65)--cycle;
  \draw[utedge,very thick] (.12,.12)--(1.65,1.65);
  \draw[utedge,thick,dashed] (-.12,.12)--(-1.65,1.65);
  \draw[utcyan,densely dotted] (0,1.12)--(0,1.66);
  \foreach \yy in {1.23,1.43,1.63}
    \draw[utcyan,fill=white,thick] (0,\yy) circle[radius=2pt];
  \filldraw[fill=utpurple!25,draw=utpurple,thick] (0,1) ellipse[x radius=.22cm,y radius=.16cm];
  \draw[->,black!65] (-1.85,0)--(1.94,0);
  \node[below] at (1.8,-.15) {$X=x/\sqrt L$};
  \draw[->,black!65] (0,-.02)--(0,1.79) node[above] {$t$};
  \fill (0,0) circle[radius=1.6pt];
  \node[below left] at (0,0) {$(0,0)$};
  \fill (0,1) circle[radius=2.6pt];
  \node[anchor=east,fill=white,inner sep=2pt] at (-.17,1.08) {conditioning point $(0,1)$};
  \node[anchor=north] at (-1,0) {$-1$};
  \node[anchor=north] at (1,0) {$1$};
  \node[black!65,align=center,font=\footnotesize] at (-.43,.72) {affected\\region};
  \node[black!60,align=center,font=\footnotesize] at (-1.38,.47) {unchanged\\leading shape};
  \node[utedge,rotate=56,fill=white,inner sep=1.5pt,font=\footnotesize] at (1.27,1.16) {$X=t$};
  \node[anchor=west,font=\small\bfseries] at (2.15,1.68) {Known results in different regions};
  \fill[black!17] (2.15,1.46) rectangle (2.31,1.51);
  \node[anchor=west,text width=5.5cm,font=\footnotesize,align=left] at (2.4,1.485)
    {Macroscopic limit shape: $0<t\le1$.\ Dark gray: $|X|\le t$ is affected.};
  \filldraw[fill=utblue!35,draw=utblue] (2.15,1.19) rectangle (2.31,1.24);
  \node[anchor=west,text width=5.5cm,font=\footnotesize,align=left] at (2.4,1.215)
    {Central height fluctuations: $x=O(L^{-1/4})$, $0<t<1$.\ Two-Brownian-bridge limit.};
  \filldraw[fill=utpurple!25,draw=utpurple] (2.23,.94) ellipse[x radius=.15cm,y radius=.08cm];
  \node[anchor=west,text width=5.5cm,font=\footnotesize,align=left] at (2.4,.94)
    {Microscopic upper-tail field:\ $x=O(L^{-1})$, $t-1=O(L^{-3/2})$.};
  \draw[utcyan,fill=white,thick] (2.23,.67) circle[radius=2pt];
  \node[anchor=west,text width=5.5cm,font=\footnotesize,align=left] at (2.4,.67)
    {Later-time one-point limits:\ fixed $x=O(1)$ and $t>1$.};
  \draw[utedge,very thick] (2.15,.405)--(2.31,.405);
  \node[anchor=west,text width=5.5cm,font=\footnotesize,align=left] at (2.4,.405)
    {This paper: $x=t\sqrt L+O(1)$, $t>0$.\ Joint space-time edge fluctuations;\ dashed left edge follows by reflection.};
\end{tikzpicture}
\caption{Space-time locations of shape and fluctuation results for the KPZ fixed point $\hKPZ(x,t)$ under the upper-tail conditioning $\{\hKPZ(0,1)\ge L\}$. The horizontal coordinate is $X=x/\sqrt L$.}
\label{fig:upper-tail-spacetime}
\end{figure}
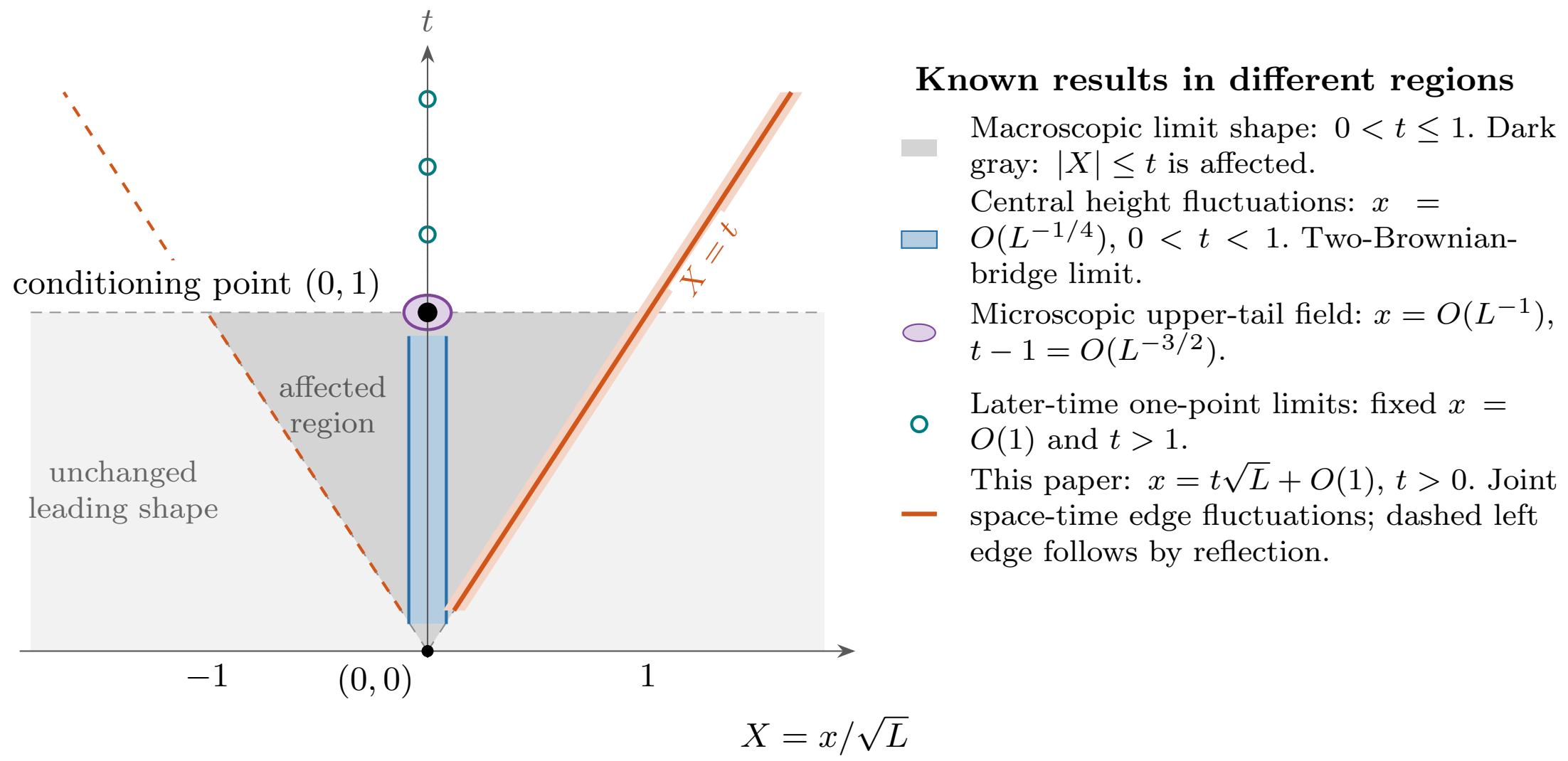

Our main contribution is a joint space-time limit at the edge of the region affected by the conditioning. We consider $O(1)$ spatial distances from the right edge $x=t\sqrt L$, where the conditional macroscopic profile meets the unaffected leading parabola of the KPZ fixed point with the narrow-wedge initial profile. The limiting field has the same finite-dimensional distributions as the KPZ fixed point with half-Brownian initial profile. Spatial reflection gives the corresponding left-edge limit. Figure~\ref{fig:upper-tail-spacetime} compares this window with the previously studied regimes. Before stating the main theorem, we explain heuristically the edge location and the centering.

Write $X=x/\sqrt L$. The leading profile of the rescaled height $L^{-1}\hKPZ(\sqrt L\,X,t)$ without conditioning is $-X^2/t$. When a large height is prescribed at $(0,1)$, the narrow-wedge limit-shape calculation in \cite[Example 1.3 and Corollary 1.2]{Das-Tsai24} gives
\begin{equation}\label{eq:conditional-profile-before-time-one}
\mathcal{F}_t(X)=
\begin{cases}
t-2|X|,& |X|\le t,\\
-X^2/t,& |X|\ge t,
\end{cases}
\qquad 0<t\le1.
\end{equation}
Thus, for $0<t\leq1$, the tent meets the parabola tangentially at $X=\pm t$. The case $t=1$ was also discussed in \cite{Ganguly-Hegde22}.

To understand the edge after time $1$, we use the directed landscape variational formula
\begin{equation*}
\hKPZ(x,t)=\sup_{y\in\realR}\{\hKPZ(y,1)+\dL(y,1;x,t)\},\qquad t>1.
\end{equation*}
Replacing the two terms by their macroscopic profiles leads to the following heuristic variational problem:
\begin{equation*}
\mathcal{F}_t(X)=\sup_{Y\in\realR}\left\{\mathcal{F}_1(Y)-\frac{(X-Y)^2}{t-1}\right\},\qquad t>1.
\end{equation*}
Evaluating the supremum yields
\begin{equation}\label{eq:conditional-profile-after-time-one}
\mathcal{F}_t(X)=
\begin{cases}
1-X^2/(t-1),& |X|\le t-1,\\
t-2|X|,& t-1\le |X|\le t,\\
-X^2/t,& |X|\ge t,
\end{cases}
\qquad t>1.
\end{equation}

The profiles \eqref{eq:conditional-profile-before-time-one} and \eqref{eq:conditional-profile-after-time-one} predict that, for every $t>0$, the conditioning changes the unconditioned parabola $-X^2/t$ precisely on $|X|<t$. The outer linear branch and the parabola have the same value and slope at $X=\pm t$. Near the right edge, the difference between their extensions is
\begin{equation*}
(t-2X)-\left(-\frac{X^2}{t}\right)=\frac{(X-t)^2}{t}.
\end{equation*}
This suggests the edge location $x=t\sqrt L$ for all positive times, including $t>1$. It also explains the spatial scale in our theorem. At $x=t\sqrt L+\alpha$, with $\alpha$ fixed, the difference in the original height scale is $L(X-t)^2/t=\alpha^2/t$, which is of order one.

\medskip
For $L>0$, define the rescaled field
\begin{equation}\label{eq:rescaled-edge-field}
\mathcal H_L^{\mathrm{edge}}(\alpha,\tau)
:=\hKPZ(\tau\sqrt L+\alpha,\tau)+\tau L+2\alpha\sqrt L,
\qquad (\alpha,\tau)\in\realR\times(0,\infty).
\end{equation}
Here we subtract the common tangent height $-\tau L-2\alpha\sqrt L$ of the linear and parabolic profiles. Our main result is the following. Throughout the paper, $\xrightarrow{\fdd}$ denotes convergence of finite-dimensional distributions.

\begin{thm}\label{thm:main}
As $L\to\infty$,
\begin{equation*}
\operatorname{Law}\left(
\left\{\mathcal H_L^{\mathrm{edge}}(\alpha,\tau)\right\}_{\alpha\in\realR,\,\tau>0}
\,\middle|\,\hKPZ(0,1)\ge L
\right)
\xrightarrow{\fdd}
\operatorname{Law}\left(
\left\{\hKPZhB\right\}_{\alpha\in\realR,\,\tau>0}
\right),
\end{equation*}
where $\hKPZhB$ is the KPZ fixed point started from the half-Brownian (one-sided Brownian) initial profile
\begin{equation*}
\mathfrak{h}_{\mathrm{hB}}(\alpha)=
\begin{cases}
\rB(-\alpha), & \alpha\leq0,\\
-\infty, & \alpha>0.
\end{cases}
\end{equation*}
Here $\rB$ is a Brownian motion on $[0,\infty)$ with $\rB(0)=0$ and $\Var(\rB(s))=2s$, independent of the directed landscape driving the  KPZ fixed point.
\end{thm}

\begin{rmk}\label{rmk:left-edge}
The narrow-wedge KPZ fixed point and the conditioning event $\{\hKPZ(0,1)\ge L\}$ are invariant in law under spatial reflection. Consequently, Theorem~\ref{thm:main} also describes fluctuations at the left edge $x=-\tau\sqrt L$. With the spatial displacement measured in the original coordinate, define
\begin{equation*}
\mathcal H_L^{\mathrm{edge;left}}(\alpha,\tau)
:=\hKPZ(-\tau\sqrt L+\alpha,\tau)+\tau L-2\alpha\sqrt L.
\end{equation*}
Under the conditioning, this field has the same finite-dimensional distributions as $\mathcal H_L^{\mathrm{edge}}(-\alpha,\tau)$. Its limit is therefore the reflected field $\hKPZ(-\alpha,\tau;\mathfrak{h}_{\mathrm{hB}})$.
\end{rmk}

\subsection{Multipoint distribution formulas}\label{sec:half-brownian-distribution-formulas}

Multipoint space-time formulas for the KPZ fixed point are known for narrow-wedge \cite{Johansson-Rahman21,Liu22}, flat \cite{Liu22}, and compactly supported initial profiles \cite{Liao-Liu25}. A two-time formula for a two-sided Brownian initial profile and a simplified multitime formula for the narrow-wedge initial profile were obtained in \cite{Rahman26}.

For the half-Brownian initial profile, a two-time formula was obtained in \cite{Johansson-Rahman22}. Here we record the multipoint cumulative distribution formula for $\hKPZhB$, obtained by specializing the limiting formulas for inhomogeneous last-passage percolation \cite{Liao21} and identifying the limit through convergence to the directed landscape \cite{Dauvergne-Virag21}. We also derive a joint upper-tail representation suited to the asymptotic analysis. These formulas will be used to identify the limit in Theorem~\ref{thm:main}; their proofs are given in Sections~\ref{sec:half-brownian-cdf-proof} and~\ref{sec:half-brownian-upper-tail-proof}, respectively.

Fix $\mm\geq1$. We first specify the contours. The $2\mm-1$ left contours are pairwise disjoint and ordered from left to right as
\[
\Sigma_{\mm,\rL}^{\inn},\ldots,\Sigma_{2,\rL}^{\inn},
\Sigma_{1,\rL},
\Sigma_{2,\rL}^{\out},\ldots,\Sigma_{\mm,\rL}^{\out}.
\]
Each contour crosses the real axis to the right of zero, has straight tails in directions $\pm2\pi/3$, and is oriented from the lower tail to the upper tail.

\begin{figure}[htbp]
\centering
\begin{tikzpicture}[
    x=1cm,y=1cm,line cap=round,line join=round,
    contour/.style={line width=1.15pt,rounded corners=5pt,
      postaction={decorate},
      decoration={markings,mark=at position 0.72 with {\arrow{>}}}},
    contour label/.style={inner sep=2pt,font=\small}
]
\draw[->,gray!70,line width=0.55pt] (-3.8,0)--(8.9,0) node[right,black] {$\Re$};
\draw[->,gray!70,line width=0.55pt] (-3.2,-1.65)--(-3.2,1.9) node[above,black] {$\Im$};
\fill (-3.2,0) circle[radius=2pt];
\node[below left,inner sep=3pt] at (-3.2,0) {$0$};
\foreach \x in {-1.6,-0.4,0.4,1.2,2.4}{
  \begin{scope}[xshift=\x cm]
    \draw[blue,contour,line width=1.35pt] (-120:1.50)--(0,0)--(120:1.50);
  \end{scope}
}
\foreach \x in {3,4.3,5,5.7,7}{
  \begin{scope}[xshift=\x cm]
    \draw[red,contour,line width=1.35pt] (-36:2.2)--(0,0)--(36:2.2);
  \end{scope}
}
\node[contour label,anchor=south] at (-2.3,1.34) {$\Sigma_{\mm,\rL}^{\inn}$};
\node[contour label,anchor=south] at (-1.2,1.34) {$\Sigma_{2,\rL}^{\inn}$};
\node[contour label,anchor=south] at (-0.35,1.34) {$\Sigma_{1,\rL}$};
\node[contour label,anchor=south] at (0.50,1.34) {$\Sigma_{2,\rL}^{\out}$};
\node[contour label,anchor=south] at (1.6,1.34) {$\Sigma_{\mm,\rL}^{\out}$};
\node[contour label,anchor=south] at (4.5,1.34) {$\Sigma_{\mm,\rR}^{\out}$};
\node[contour label,anchor=south] at (5.8,1.34) {$\Sigma_{2,\rR}^{\out}$};
\node[contour label,anchor=south] at (6.8,1.34) {$\Sigma_{1,\rR}$};
\node[contour label,anchor=south] at (7.7,1.34) {$\Sigma_{2,\rR}^{\inn}$};
\node[contour label,anchor=south] at (8.82,1.34) {$\Sigma_{\mm,\rR}^{\inn}$};
\foreach \x in {-1.3,1.5,4.3,7.1}
  \node[contour label] at (\x,0.5) {$\boldsymbol{\cdots}$};
\end{tikzpicture}
\caption{The $\Sigma$-contours and their orientations. All contours cross the real axis to the right of zero.}
\label{fig:half-brownian-contours}
\end{figure}
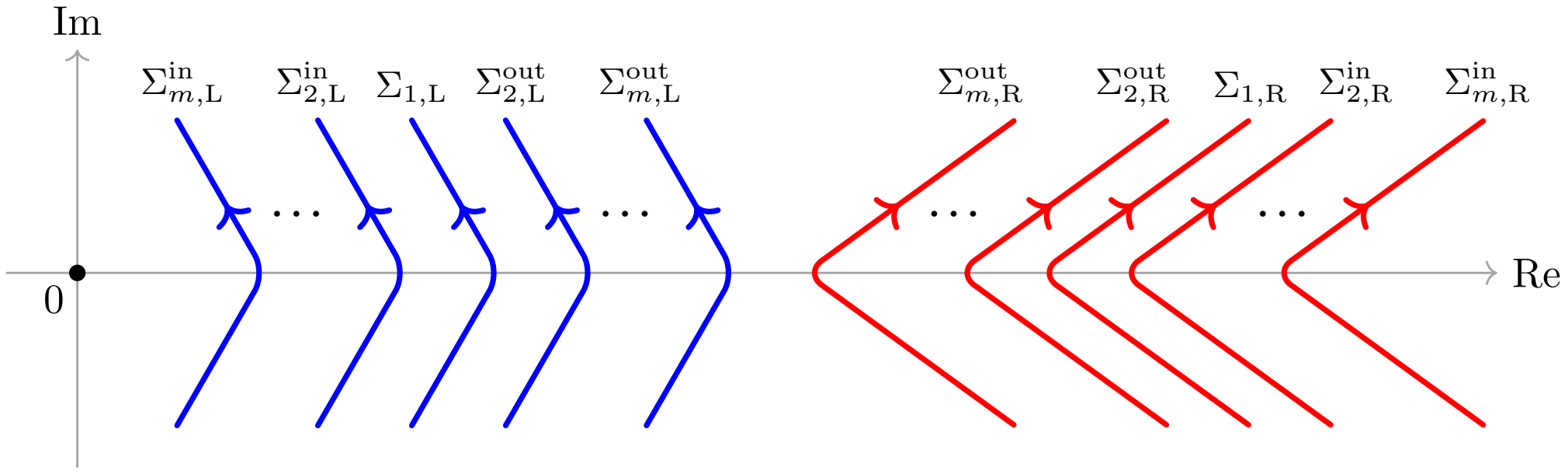

The $2\mm-1$ right contours are pairwise disjoint and ordered from left to right as
\[
\Sigma_{\mm,\rR}^{\out},\ldots,\Sigma_{2,\rR}^{\out},
\Sigma_{1,\rR},
\Sigma_{2,\rR}^{\inn},\ldots,\Sigma_{\mm,\rR}^{\inn}.
\]
Each right contour has straight tails in directions $\pm\pi/5$ and is oriented from the lower half plane to the upper half plane. Every left contour lies to the left of every right contour, and the distance between any two distinct contours is bounded below by a fixed constant $c>0$. For the left family, the $\inn$-contours are nested to the left of $\Sigma_{1,\rL}$ and the $\out$-contours to its right. For the right family, the $\inn$-contours are nested to the right of $\Sigma_{1,\rR}$ and the $\out$-contours to its left. For $\mm=1$, only $\Sigma_{1,\rL}$ and $\Sigma_{1,\rR}$ are present. Figure~\ref{fig:half-brownian-contours} illustrates the nesting.

We use the following order on space-time points, comparing time first and then space:
\begin{equation*}
(\alpha,\tau)\prec(\alpha',\tau')
\quad\Longleftrightarrow\quad
 \tau<\tau'\quad\text{or}\quad(\tau=\tau'\ \text{and}\ \alpha<\alpha').
\end{equation*}

For half-Brownian space-time points $(\alpha_j,\tau_j)\in\realR\times(0,\infty)$ satisfying $(\alpha_1,\tau_1)\prec\cdots\prec(\alpha_\mm,\tau_\mm)$ and height parameters $\beta_1,\ldots,\beta_\mm\in\realR$, define the exponential factors and ratios of consecutive factors by
\begin{equation}\begin{split}\label{eq:half-brownian-phases-and-increments}
g_j(\zeta)&:=\exp\left\{-\frac{\tau_j}{3}\zeta^3+\alpha_j\zeta^2+\beta_j\zeta\right\},
\qquad 1\leq j\leq \mm,\\
G_1(\zeta)&:=g_1(\zeta),
\qquad
G_j(\zeta):=\frac{g_j(\zeta)}{g_{j-1}(\zeta)},
\qquad 2\leq j\leq \mm.
\end{split}\end{equation}

For two vectors 
\begin{equation*}
  \bzeta=(\zeta_1,\ldots,\zeta_n)\in\complexC^n,\quad \bzeta'=(\zeta'_1,\ldots,\zeta'_{n'})\in\complexC^{n'},
\end{equation*}
let $\bzeta\sqcup\bzeta'=(\zeta_1,\ldots,\zeta_n,\zeta'_1,\ldots,\zeta'_{n'})$ denote their concatenation. When $n=n'$ and $\zeta_i\ne\zeta'_j$ for all $i,j$, define the Cauchy determinant
\begin{equation*}
  \cC(\bzeta;\bzeta')= \det\left[ \frac{1}{\zeta_i-\zeta'_j}\right]_{i,j=1}^n = (-1)^{\binom{n}{2}}\frac{\prod_{i<j}(\zeta_j-\zeta_i)(\zeta'_j-\zeta'_i)}{\prod_i\prod_j (\zeta_i-\zeta'_j)}.
\end{equation*}
The following formula is a half-Brownian specialization of \cite[Theorem~3.2.1(i)]{Liao21}, with the identification and extension to coincident times justified in Section~\ref{sec:half-brownian-cdf-proof}.
\begin{thm}\label{thm:half-brownian-cdf}
Fix $\mm\geq1$ and space-time points $(\alpha_\ell,\tau_\ell)\in\realR\times(0,\infty)$, $1\leq\ell\leq\mm$, satisfying $(\alpha_1,\tau_1)\prec\cdots\prec(\alpha_\mm,\tau_\mm)$. For any height parameters $\bbeta=(\beta_1,\ldots,\beta_\mm)\in\realR^\mm$, the joint cumulative distribution function is given by
\begin{equation}
\begin{split}
 \label{eq:half-brownian-cdf}
 &\prob\left(\bigcap_{\ell=1}^{\mm}\left\{\hKPZ(\alpha_\ell,\tau_\ell;\mathfrak{h}_{\mathrm{hB}})\leq\beta_\ell\right\}\right)\\
  &\quad =\sum_{\nn_1,\ldots,\nn_\mm\geq0}
\frac{1}{(\nn_1!\cdots\nn_\mm!)^2}
\oint_{<1}\cdots\oint_{<1}\rK_{\bn}(z_1,\ldots,z_{\mm-1})
\prod_{\ell=1}^{\mm-1}\frac{\rd z_\ell}{2\pi\rmi z_\ell(1-z_\ell)},
  \end{split}
\end{equation}
where each $z_\ell$-contour is a counterclockwise circle centered at the origin with an arbitrary fixed radius in $(0,1)$. For $\bn=(\nn_1,\ldots,\nn_\mm)$, the function $\rK_\bn$ is defined by
\begin{equation}
\begin{split}
 \label{eq:half-brownian-coefficient}
 \rK_\bn(z_1,\ldots,z_{\mm-1})
 &= \rK_\bn(z_1,\ldots,z_{\mm-1};\bbeta; (\alpha_1,\tau_1),\ldots,(\alpha_\mm,\tau_\mm))\\
&:= \prod_{\ell=1}^{\mm-1}(1-z_\ell)^{\nn_\ell}(1-z_\ell^{-1})^{\nn_{\ell+1}}\\
&\quad\cdot\prod_{i_1=1}^{\nn_1}\int_{\Sigma_{1,\rL}}\frac{\rd\xi_{i_1}^{(1)}}{2\pi\rmi}
\cdot\prod_{\ell=2}^\mm\prod_{i_\ell=1}^{\nn_\ell}\left(
\frac{1}{1-z_{\ell-1}}\int_{\Sigma_{\ell,\rL}^{\inn}}\frac{\rd\xi_{i_\ell}^{(\ell)}}{2\pi\rmi}
-\frac{z_{\ell-1}}{1-z_{\ell-1}}\int_{\Sigma_{\ell,\rL}^{\out}}\frac{\rd\xi_{i_\ell}^{(\ell)}}{2\pi\rmi}\right)\\
&\quad\cdot\prod_{i_1=1}^{\nn_1}\int_{\Sigma_{1,\rR}}\frac{\rd\eta_{i_1}^{(1)}}{2\pi\rmi}
\cdot\prod_{\ell=2}^\mm\prod_{i_\ell=1}^{\nn_\ell}\left(
\frac{1}{1-z_{\ell-1}}\int_{\Sigma_{\ell,\rR}^{\inn}}\frac{\rd\eta_{i_\ell}^{(\ell)}}{2\pi\rmi}
-\frac{z_{\ell-1}}{1-z_{\ell-1}}\int_{\Sigma_{\ell,\rR}^{\out}}\frac{\rd\eta_{i_\ell}^{(\ell)}}{2\pi\rmi}\right)\\
&\quad \cdot  \prod_{\ell=0}^\mm \cC(\bxi^{(\ell)}\sqcup \vecteta^{(\ell+1)}; \vecteta^{(\ell)}\sqcup \bxi^{(\ell+1)}) \cdot \prod_{i_1=1}^{\nn_1} \frac{\eta_{i_1}^{(1)}}{\xi_{i_1}^{(1)}} \cdot \prod_{\ell=1}^\mm\prod_{i_\ell=1}^{\nn_\ell} \frac{G_\ell(\xi_{i_\ell}^{(\ell)})}{G_\ell (\eta_{i_\ell}^{(\ell)})}.
\end{split}
\end{equation}
Here, for $1\leq\ell\leq\mm$, we write
\[
\bxi^{(\ell)}=(\xi_1^{(\ell)},\ldots,\xi_{\nn_\ell}^{(\ell)}),
\qquad
\vecteta^{(\ell)}=(\eta_1^{(\ell)},\ldots,\eta_{\nn_\ell}^{(\ell)}),
\]
with both vectors empty when $\nn_\ell=0$. At the auxiliary indices $0$ and $\mm+1$, which carry no integration variables, set $\bxi^{(0)}=\bxi^{(\mm+1)}=\vecteta^{(0)}=\vecteta^{(\mm+1)}=\emptyset$. Throughout the paper, empty products and the Cauchy determinant of two empty vectors are taken to be $1$. When $\mm=1$, the circle integrals are omitted.
\end{thm}

\begin{rmk}
The series \eqref{eq:half-brownian-cdf} is absolutely convergent. Indeed, the first-layer weights have cubic exponential decay since $\tau_1>0$. For $\ell\geq2$, the factors $G_\ell$ on the left contours and $G_\ell^{-1}$ on the right contours have cubic exponential decay if $\tau_\ell>\tau_{\ell-1}$, and Gaussian decay if $\tau_\ell=\tau_{\ell-1}$, since then $\alpha_\ell>\alpha_{\ell-1}$. Contour separation and Hadamard's inequality give a bound that is summable over the half-Brownian series indices $\nn_1,\ldots,\nn_\mm$. The details are given in Lemma~\ref{lm:half-brownian-absolute-convergence}; see also \cite[Lemma~3.3]{Liu22c} and \cite{Liu-Zhang25} for similar estimates.
\end{rmk}

For the asymptotic analysis, we need a formula for joint upper-tail probabilities. Adapting the circle-contour deformation and inclusion--exclusion arguments of \cite{Liu22,Liu-Wang24,Liu-Zhang25} to \eqref{eq:half-brownian-cdf} gives the following representation, in which the $z$-contours have radii greater than $1$ and the sum is restricted to $\nn_1,\ldots,\nn_\mm\geq1$.
 
\begin{prop}\label{prop:half-brownian-upper-tail}
Under the assumptions and notation of Theorem~\ref{thm:half-brownian-cdf},
\begin{equation}
\begin{split}
 \label{eq:half-brownian-upper-tail}
 &\prob\left(\bigcap_{\ell=1}^{\mm}\left\{\hKPZ(\alpha_\ell,\tau_\ell;\mathfrak{h}_{\mathrm{hB}})\geq\beta_\ell\right\}\right)\\
  &\quad =(-1)^{\mm}\sum_{\nn_1,\ldots,\nn_\mm\geq1}
\frac{1}{(\nn_1!\cdots\nn_\mm!)^2}
\oint_{>1}\cdots\oint_{>1}\rK_{\bn}(z_1,\ldots,z_{\mm-1})
\prod_{\ell=1}^{\mm-1}\frac{\rd z_\ell}{2\pi\rmi z_\ell(1-z_\ell)},
  \end{split}
\end{equation}
where each $z_\ell$-contour is a counterclockwise circle centered at the origin with an arbitrary fixed radius greater than $1$. Here $\bn=(\nn_1,\ldots,\nn_\mm)$, and the contour integrals $\rK_\bn$ are defined by \eqref{eq:half-brownian-coefficient}. The series converges absolutely, and the circle integrals are omitted when $\mm=1$.
\end{prop}

\subsection{Continuity, scaling, and fixed-time laws}\label{sec:half-brownian-properties}

We recall the continuity, scaling, and fixed-time laws of $\hKPZhB$ needed in the proof of Theorem~\ref{thm:main}. These also imply continuity of its joint CDF. We write $\xlongequal{\fdd}$ for equality of finite-dimensional distributions.

\begin{prop}[Sample-path continuity and scaling invariance]\label{prop:half-brownian-continuity-scaling}
The field $\hKPZhB$ has a version with jointly continuous sample paths on $\realR\times(0,\infty)$ and is therefore stochastically continuous. Moreover, for every $\epsilon>0$,
\begin{equation}\label{eq:half-brownian-scaling}
\left\{\epsilon^{-1}\hKPZ(\epsilon^2\alpha,\epsilon^3\tau;\mathfrak{h}_{\mathrm{hB}})\right\}_{\alpha\in\realR,\,\tau>0}
\xlongequal{\fdd}
\left\{\hKPZhB\right\}_{\alpha\in\realR,\,\tau>0}.
\end{equation}
\end{prop}

\begin{proof}
Almost surely, $\mathfrak{h}_{\mathrm{hB}}(\alpha)\leq C(1+|\alpha|)$ for all $\alpha\in\realR$ and some finite random constant $C$. Thus the directed landscape construction gives a version of the KPZ fixed point that is jointly continuous at positive times \cite[Section~16.1]{Dauvergne-Virag21}. This also implies stochastic continuity.

Brownian scaling gives
\[
\epsilon^{-1}\mathfrak{h}_{\mathrm{hB}}(\epsilon^2\,\cdot)
\stackrel{\mathrm d}{=}\mathfrak{h}_{\mathrm{hB}}.
\]
Together with the $1:2:3$ scaling of the KPZ fixed point, this proves \eqref{eq:half-brownian-scaling}.
\end{proof}

The fixed-time laws are expressed in terms of the Brownian-to-Airy$_2$ transition process introduced in \cite{Imamura-Sasamoto04}, following earlier work on one-point transition distributions for polynuclear growth with external sources \cite{Baik-Rains00}. We denote this process by $\Airy_{\mathrm{BM}\to2}$ and use the normalization in \cite[Definition~2.13]{Corwin-Ferrari-Peche10}. For its one-point laws, we write $F_{\mathrm{BBP},1,x}$ for the rank-one Baik--Ben Arous--P\'ech\'e (BBP) distribution function \cite{Baik-Ben_Arous-Peche05} with deformation parameter $x$, normalized as in \cite[Definition~1.6]{Borodin-Corwin-Ferrari14}.

\begin{prop}[Fixed-time spatial and one-point laws]\label{prop:half-brownian-fixed-time-laws}
For every $\tau>0$, $\mm\geq1$, and $\alpha_1<\cdots<\alpha_\mm$,
\begin{equation}\label{eq:half-brownian-fixed-time-law}
\left(\hKPZ(\alpha_\ell,\tau;\mathfrak{h}_{\mathrm{hB}})\right)_{\ell=1}^{\mm}
\stackrel{\mathrm d}{=}
\left(\tau^{1/3}\Airy_{\mathrm{BM}\to2}(\alpha_\ell\tau^{-2/3})-\frac{\alpha_\ell^2}{\tau}\right)_{\ell=1}^{\mm}.
\end{equation}
Consequently, for every $(\alpha,\beta,\tau)\in\realR^2\times(0,\infty)$,
\begin{equation}\label{eq:half-brownian-one-point-bbp}
\prob\left(\hKPZhB\leq\beta\right)
=F_{\mathrm{BBP},1,-\alpha\tau^{-2/3}}
\left(\tau^{-1/3}\beta+\alpha^2\tau^{-4/3}\right).
\end{equation}
\end{prop}

\begin{proof}
Taking $\epsilon=\tau^{1/3}$ in the scaling relation of Proposition~\ref{prop:half-brownian-continuity-scaling} gives
\[
\left\{\hKPZhB\right\}_{\alpha\in\realR}
\xlongequal{\fdd}
\left\{\tau^{1/3}\hKPZ(\alpha\tau^{-2/3},1;\mathfrak{h}_{\mathrm{hB}})\right\}_{\alpha\in\realR}.
\]
At time $1$, the KPZ fixed point with this initial profile has the spatial law $\{\Airy_{\mathrm{BM}\to2}(x)-x^2\}_{x\in\realR}$; see \cite{Imamura-Sasamoto04,Corwin-Ferrari-Peche10}. Combining these identities proves \eqref{eq:half-brownian-fixed-time-law}. In the chosen normalization,
\[
\prob\bigl(\Airy_{\mathrm{BM}\to2}(x)\leq s\bigr)
=F_{\mathrm{BBP},1,-x}(s).
\]
Taking $x=\alpha\tau^{-2/3}$ and $s=\tau^{-1/3}\beta+\alpha^2\tau^{-4/3}$ yields \eqref{eq:half-brownian-one-point-bbp}.
\end{proof}

Thus the normalized height $\tau^{-1/3}(\hKPZhB+\alpha^2/\tau)$ has the rank-one BBP law with deformation parameter $-\alpha\tau^{-2/3}$. Rank-one BBP laws also occur at the transition boundaries in exponential directed last-passage percolation under upper-tail conditioning \cite{Baik-Tripathi26}.

\begin{rmk}[Borodin--P\'ech\'e and Airy-wanderer descriptions]
The process $\Airy_{\mathrm{BM}\to2}$ can equivalently be described as the rank-one Borodin--P\'ech\'e process with deformation sets $\emptyset$ and $\{0\}$ \cite{Borodin-Peche08}, or as a one-sided rank-one specialization of the Airy processes with wanderers \cite{Adler-Ferrari-van_Moerbeke10}.
\end{rmk}

Combining sample-path continuity with the atomlessness of the one-point laws gives the following consequence.
\begin{cor}[Joint CDF continuity]\label{cor:half-brownian-joint-cdf-continuity}
For every $\mm\geq1$, the map
\begin{equation}\label{eq:half-brownian-joint-cdf-continuity}
\bigl((\alpha_\ell,\tau_\ell,\beta_\ell)\bigr)_{\ell=1}^{\mm}\longmapsto
\prob\bigl(\hKPZ(\alpha_\ell,\tau_\ell;\mathfrak{h}_{\mathrm{hB}})\leq\beta_\ell,\ 1\leq \ell\leq\mm\bigr)
\end{equation}
is continuous on $(\realR\times(0,\infty)\times\realR)^\mm$.
\end{cor}

\begin{proof}
By Proposition~\ref{prop:half-brownian-fixed-time-laws}, the one-point marginals are atomless. For any convergent sequence of space-time points and thresholds, the jointly continuous version from Proposition~\ref{prop:half-brownian-continuity-scaling} gives almost sure convergence of the corresponding event indicators, since each limiting height equals its threshold with probability zero. The bounded convergence theorem proves continuity of \eqref{eq:half-brownian-joint-cdf-continuity}.
\end{proof}

\subsection{Proof strategy and organization}\label{sec:proof-strategy}

The proof combines multipoint distribution formulas for the half-Brownian KPZ fixed point with asymptotics of the conditional narrow-wedge probabilities.

In Section~\ref{sec:half-brownian-formula-proofs}, we identify the limiting critical-boundary last-passage formula of \cite{Liao21} as the multipoint distribution of the half-Brownian KPZ fixed point. We then convert the resulting CDF into the joint upper-tail representation by deforming the circle contours.

In Section~\ref{sec:contour-decompositions}, we introduce the narrow-wedge multipoint formula from \cite{Liu22,Liu-Zhang25} and derive separate residue decompositions for the contour integrals in the narrow-wedge and half-Brownian joint upper-tail series.

In Section~\ref{sec:main-theorem-proof}, we insert the conditioning point into the ordered list of space-time points and translate all contours by $\sqrt L$. Through further contour deformations and estimates, we obtain two parallel results: Proposition~\ref{prop:narrow-wedge-fixed-order-asymptotics} and Lemma~\ref{lm:terminal-fixed-order-asymptotics}. Both identify the leading contributions in the narrow-wedge series and determine their limits after normalization. After establishing a summable bound uniform in $L$, we pass to the limit under the sum. Reindexing the resulting series then recovers the half-Brownian residue expansion, thereby proving convergence at distinct positive times $\tau_1,\ldots,\tau_\mm$. 

At the end of Section~\ref{sec:main-theorem-proof}, we use joint CDF continuity of the limiting field and \cite[Lemma~3.6]{Liu-Wang24} to extend the convergence to coincident times, completing the proof of Theorem~\ref{thm:main}.

\section*{Acknowledgement}

This project was initiated and carried out in part during H.Y.'s visit to the University of Kansas from spring 2025 to spring 2026, when Z.L. was based there. Z.L. thanks Duncan Dauvergne for suggesting possible connections between the limiting field and the Airy process with wanderers during the conference ``The Kardar--Parisi--Zhang Universality Class \& Related Topics'', held at the BRIN Mathematics Research Center in April 2026. Z.L.'s research is supported in part by the Spofford Harris Kimball Endowed Fund at the University of Maine and by NSF grants DMS-2552334/DMS-2639619 and DMS-2246683/DMS-2639616. Y.L. and H.Y. are supported in part by the NSFC grant 12471148.

Our initial asymptotic analysis treated times before, at, and after $1$ separately. At the final stage, ChatGPT 6 was used to assist us in developing a unified treatment of these cases and in writing the manuscript.

\section{Multipoint formulas for the half-Brownian initial profile}\label{sec:half-brownian-formula-proofs}

This section supplies the analytic details for the formulas in Section~\ref{sec:half-brownian-distribution-formulas}. After establishing absolute convergence of the series, we identify the limiting critical-boundary formula of \cite{Liao21} as the joint CDF of the half-Brownian KPZ fixed point. Deforming the circle contours then gives the joint upper-tail representation.

We also allow some indices $\nn_\ell$ in the half-Brownian series to be zero. If $\nn_\ell=0$, then $\bxi^{(\ell)}=\vecteta^{(\ell)}=\emptyset$ in \eqref{eq:half-brownian-coefficient}, and the integrals over the layer-$\ell$ variables are omitted.

Set
\begin{equation}\label{eq:half-brownian-integrated-coefficient}
\hatK_{\bn}
:=(-1)^\mm
\left[\prod_{\ell=1}^{\mm-1}\oint_{>1}
\frac{\rd z_\ell}{2\pi\rmi z_\ell(1-z_\ell)}\right]
\rK_{\bn}(z_1,\ldots,z_{\mm-1}).
\end{equation}
For each fixed choice of the half-Brownian series indices $\bn$, the integrals converge absolutely by the decay of the exponential factors and contour separation described in Section~\ref{sec:half-brownian-distribution-formulas}.

For later use, write
\begin{equation}\label{eq:half-brownian-upper-tail-series}
\rT_\mm(\bbeta;(\alpha_1,\tau_1),\ldots,(\alpha_\mm,\tau_\mm))
:=\sum_{\nn_1,\ldots,\nn_\mm\geq1}
\frac{\hatK_{\bn}}{(\nn_1!\cdots\nn_\mm!)^2}.
\end{equation}
With this notation, Proposition~\ref{prop:half-brownian-upper-tail} states that $\rT_\mm$ is the joint upper-tail probability for the KPZ fixed point with initial profile $\mathfrak h_{\mathrm{hB}}$.

\subsection{Absolute convergence}\label{sec:half-brownian-absolute-convergence}

We choose the contours to consist of two straight rays with the directions and nesting specified in Section~\ref{sec:half-brownian-distribution-formulas}. They stay away from zero, and any contours appearing together in a Cauchy denominator have positive separation. Let $a_1>0$ be the real-axis vertex of $\Sigma_{1,\rL}$. We use the convention $0^0=1$.

If $X=(x_1,\ldots,x_n)$ and $Y=(y_1,\ldots,y_n)$ satisfy
$|x_p-y_q|\geq d>0$, then Hadamard's inequality gives
\begin{equation}\label{eq:cauchy-hadamard-bound}
|\cC(X;Y)|
\leq \prod_{p=1}^n
\left(\sum_{q=1}^n|x_p-y_q|^{-2}\right)^{1/2}
\leq n^{n/2}d^{-n}.
\end{equation}

\begin{lm}\label{lm:half-brownian-absolute-convergence}
For fixed half-Brownian parameters, contours, and circle radii, there is
a constant $C>0$, independent of the half-Brownian series indices $\bn$, such that
\begin{equation}\label{eq:half-brownian-order-majorant}
\frac{|\hatK_{(\nn_1,\ldots,\nn_\mm)}|}
{(\nn_1!\cdots\nn_\mm!)^2}
\leq C\prod_{\ell=1}^{\mm}
\frac{C^{\nn_\ell}\nn_\ell^{\nn_\ell}}{(\nn_\ell!)^2},
\qquad \nn_\ell\geq0.
\end{equation}
The same bounds hold for the sums of absolute integrals obtained
by expanding the in/out contour functionals. They remain valid if
the $z$-circles are replaced by circles centered at zero with fixed
radii in $(0,1)\cup(1,\infty)$, uniformly when the radii range over
compact subsets of this set.

Consequently, the series in \eqref{eq:half-brownian-cdf} and \eqref{eq:half-brownian-upper-tail-series} converge absolutely. The bounds are locally uniform in the height parameters. They are also uniform along small perturbations of the times that preserve the ordering of the space-time points, including approximations of equal times by strictly increasing times.
\end{lm}

\begin{proof}
We expand the weighted contour integrals and take absolute values. Since the $z$-circles are fixed, their weights and the in/out coefficients contribute at most $C^{1+\nn_1+\cdots+\nn_\mm}$. This bound is uniform for radii in compact subsets of $(0,1)\cup(1,\infty)$.

For $\ell\geq2$, the factors $G_\ell(\xi)$ and $G_\ell(\eta)^{-1}$ decay super-exponentially on their respective contours. The cubic terms give this decay when $\tau_\ell>\tau_{\ell-1}$. When $\tau_\ell=\tau_{\ell-1}$, the quadratic terms give Gaussian decay since $\alpha_\ell>\alpha_{\ell-1}$. The first-layer factors $G_1(\xi)/\xi$ and $\eta/G_1(\eta)$ are also integrable, since $\tau_1>0$ and the contours avoid zero. Hence each one-variable absolute integral is bounded by a constant.

For $\hatK_{\bn}$, the first and last Cauchy determinants have
orders $\nn_1$ and $\nn_\mm$, and the determinant between layers
$\ell$ and $\ell+1$ has order $\nn_\ell+\nn_{\ell+1}$.
Combining \eqref{eq:cauchy-hadamard-bound} with the preceding
bounds gives
\begin{equation}\label{eq:half-brownian-determinant-bound}
|\hatK_{(\nn_1,\ldots,\nn_\mm)}|
\leq C^{1+\nn_1+\cdots+\nn_\mm}
\nn_1^{\nn_1/2}
\prod_{\ell=1}^{\mm-1}
(\nn_\ell+\nn_{\ell+1})^{(\nn_\ell+\nn_{\ell+1})/2}
\nn_\mm^{\nn_\mm/2}.
\end{equation}
Since
\begin{equation}\label{eq:two-order-splitting}
(p+q)^{(p+q)/2}
\leq 2^{(p+q)/2}p^{p/2}q^{q/2},
\qquad p,q\geq0,
\end{equation}
we obtain \eqref{eq:half-brownian-order-majorant} after increasing $C$.

Finally, Stirling's bound gives
\begin{equation}\label{eq:factorial-series-summability}
\sum_{\nn\geq0}\frac{C^\nn\nn^\nn}{(\nn!)^2}
\leq 1+\sum_{\nn\geq1}\left(\frac{Ce^2}{\nn}\right)^\nn
<\infty,
\end{equation}
which proves absolute summability. The same argument bounds the sum of absolute integrals, and the bounds remain valid on truncated contours.

For fixed ordered points, the bounds are uniform when the height parameters range over a compact set. Continuity in these parameters follows by dominated convergence. We can also keep the same contours under small ordered perturbations of the times. A positive time increment still gives cubic decay, while a zero time increment has a strictly positive spatial increment and hence uniform Gaussian decay.
\end{proof}

\subsection{Proof of Theorem~\ref{thm:half-brownian-cdf}}\label{sec:half-brownian-cdf-proof}

We identify the limiting formula in \cite{Liao21} with the KPZ fixed point started from $\mathfrak h_{\mathrm{hB}}$. The argument uses exponential directed last-passage percolation (DLPP) with one critical first column.

\subsubsection{Convergence to the KPZ fixed point}\label{sec:critical-boundary-kpz-convergence}
Let $L(i,j)$ be the last passage value from $(1,1)$ to $(i,j)$ for independent exponential weights $\omega_{i,j}$ of rate $1/2$ on the first column $\{(i,1):i\geq1\}$ and rate $1$ elsewhere. The rescaled inverse height $Z_T$ is characterized, for fixed $(\alpha,\tau,\beta)$ and large $T$, by
\begin{equation}\begin{split}\label{eq:critical-boundary-inverse-event}
\{Z_T(\alpha,\tau)\leq\beta\}
&=\{L(k_T,a_T+k_T)\leq2\tau T\},\\
a_T=\lfloor2\alpha T^{2/3}\rfloor,\qquad
k_T&=\left\lceil\frac{\tau T-a_T-\beta T^{1/3}}2\right\rceil.
\end{split}\end{equation}

Decomposing a path at its exit from the first column identifies the cumulative column weights as a boundary profile for the homogeneous rate-$1$ LPP environment outside the first column. This environment is independent of the first-column weights. By Donsker's theorem, after linear interpolation,
\[
\frac{1}{2T^{1/3}}\sum_{i=1}^{\lfloor2sT^{2/3}\rfloor}(\omega_{i,1}-2)
\Longrightarrow\rB(s),\qquad s\geq0,\qquad \Var(\rB(s))=2s,
\]
locally uniformly in $s$. In the scaling limit, the first column becomes the negative spatial half-line at time zero, so these Brownian fluctuations supply precisely the initial profile $\mathfrak h_{\mathrm{hB}}$ for the subsequent evolution.

This boundary-profile argument appears in \cite[Corollary~2.11(a) and Section~4.1]{Corwin-Liu-Wang16} for geometric LPP. To obtain joint convergence, first restrict the rescaled exit location to $0\leq s\leq R$. Theorem~13.8(2) of \cite{Dauvergne-Virag21} gives convergence of the rescaled homogeneous exponential LPP field to a single directed landscape, uniformly on compact sets of space-time endpoint pairs. Together with independence and the preceding Brownian convergence, this gives joint convergence of the truncated variational maxima at finitely many space-time points, with the same $\rB$ and $\mathcal L$ at every point.

Let $E_T$ be the largest $i$ for which the maximizing path to $(k_T,a_T+k_T)$ visits $(i,1)$. Then
\begin{equation}\label{eq:critical-boundary-exit-tightness}
\lim_{R\to\infty}\limsup_{T\to\infty}
\prob(E_T>RT^{2/3})=0.
\end{equation}
This follows from \cite[Theorem~2.5 and Proposition~2.7(a)]{Emrah-Janjigian-Seppalainen20} by a monotone comparison with a slightly tilted stationary exponential LPP model.

A finite union bound in \eqref{eq:critical-boundary-exit-tightness} removes the cutoff at $2RT^{2/3}$ simultaneously for any finite collection of fixed triples $(\alpha_j,\tau_j,\beta_j)$. The limiting cutoff can also be removed almost surely, since the parabolic decay of $\mathcal L$ dominates the Brownian growth. Applying the converging-together argument to the truncated maxima and using \eqref{eq:critical-boundary-inverse-event}, we obtain, jointly at any finite collection of space-time points,
\begin{equation}\label{eq:critical-boundary-kpz-limit}
Z_T(\alpha,\tau)\Longrightarrow\hKPZhB =\sup_{s\geq0}\{\rB(s)+\mathcal L(-s,0;\alpha,\tau)\}.
\end{equation}

\subsubsection{The joint CDF on full space-time}
Suppose first that $0<\tau_1<\cdots<\tau_\mm$. We apply \cite[Theorem~3.2.1(i)]{Liao21}\footnote{There is a sign typo in the cited theorem. With both limiting contour families oriented upward, the denominator in \cite[Equation (3.22)]{Liao21} should be $\zeta-\zeta'$.} with $r=0$, $s=1$, and $\mu_1=0$. The rates are $\pi_i=0$, $\widehat\pi_1=1/2$, and $\widehat\pi_j=1$ for $j\geq2$. Up to integer roundings, \eqref{eq:critical-boundary-inverse-event} is the scaling in \cite[(3.13)]{Liao21}. By Section~\ref{sec:critical-boundary-kpz-convergence}, the resulting formula gives the joint CDF of $\hKPZhB$.

Under this specialization, the functions in \cite[(3.24)]{Liao21} are $w^{-1}g_\ell(w)$. Their successive ratios are $G_\ell$ for $\ell\geq2$, and the first layer contributes $\prod_i\eta_i^{(1)}/\xi_i^{(1)}$. Expanding the limiting Fredholm determinant as in \cite[(3.121)]{Liao21} gives \eqref{eq:half-brownian-cdf}, with $\rK_\bn$ given by \eqref{eq:half-brownian-coefficient}.

To include equal times, we use the contours in Section~\ref{sec:half-brownian-distribution-formulas} and approximate the times $\tau_1,\ldots,\tau_\mm$ by strictly increasing ones. If $\tau_\ell=\tau_{\ell-1}$, then $\alpha_\ell>\alpha_{\ell-1}$, so the quadratic terms give Gaussian decay. Lemma~\ref{lm:half-brownian-absolute-convergence} gives a uniform summable bound along the approximation. We may therefore take the limit in the formula by dominated convergence. Continuity of $\hKPZhB$ and atomlessness of its marginals identify the result with the joint CDF at the original space-time points. Thus the formula gives the joint CDF of $\hKPZhB$ at all ordered space-time points. Together with \eqref{eq:critical-boundary-kpz-limit}, this proves Theorem~\ref{thm:half-brownian-cdf}.

\subsection{Proof of Proposition~\ref{prop:half-brownian-upper-tail}}\label{sec:half-brownian-upper-tail-proof}

\begin{lm}\label{lm:half-brownian-empty-layer}
Suppose $\nn_\mm\geq1$.  If $\nn_\ell=0$ for some $1\leq \ell<\mm$, then
$\hatK_{(\nn_1,\ldots,\nn_\mm)}=0$.
\end{lm}

\begin{proof}
This lemma follows from the proof of \cite[Lemma~3.1]{Liu-Wang24}, taking $\ell$ to be the largest index with $\nn_\ell=0$. Note that the extra factor $\prod_i\eta_i^{(1)}/\xi_i^{(1)}$ is independent of $z_\ell$ and the layer-$(\ell+1)$ variables, so the contour cancellation in the original proof is unchanged.
\end{proof}

We follow the proof of \cite[Proposition~3.1]{Liu-Zhang25}. The required half-Brownian analogue of the mixed upper/lower-tail formula \cite[(2.4)]{Liu-Wang24} is
\begin{equation}\label{eq:half-brownian-mixed-tail}
\prob\left(
\bigcap_{\ell=1}^{\mm-1}\left\{\hKPZ(\alpha_\ell,\tau_\ell;\mathfrak{h}_{\mathrm{hB}})>\beta_\ell\right\}
\cap\left\{\hKPZ(\alpha_\mm,\tau_\mm;\mathfrak{h}_{\mathrm{hB}})\leq\beta_\mm\right\}
\right)
=-\sum_{\nn_1,\ldots,\nn_\mm\geq0}
\frac{\hatK_{\bn}}{(\nn_1!\cdots\nn_\mm!)^2}.
\end{equation}

Apply the circle-deletion argument of \cite[Proposition~2.17 and Section~8.7]{Liu22} to Theorem~\ref{thm:half-brownian-cdf} and then use inclusion--exclusion. Only variables in layers $2,\ldots,\mm$ are deformed, so the extra first-layer factor introduces no new poles. The identity $(g_1(w)/w)(g_2(w)/g_1(w))=g_2(w)/w$ preserves the half-Brownian form when the first point is deleted. Lemma~\ref{lm:half-brownian-absolute-convergence} justifies interchanging sums and integrals; at equal times, strict spatial ordering supplies the Gaussian decay needed for the deformations.

For $\mm\geq2$, split \eqref{eq:half-brownian-mixed-tail} as in \cite{Liu-Zhang25}. Lemma~\ref{lm:half-brownian-empty-layer} gives the contribution $-\rT_\mm$ from $\nn_\mm\geq1$. The terms with $\nn_\mm=0$ contribute $\rT_{\mm-1}$, with the last space-time point and its height parameter omitted, since
\[
\oint_{>1}\frac{(1-z)^{\nn_{\mm-1}}}{z(1-z)}\frac{\rd z}{2\pi\rmi}
=\identity_{\{\nn_{\mm-1}\geq1\}},
\]
and Lemma~\ref{lm:half-brownian-empty-layer} applies to the remaining layers. The zero-index term $1$ in the one-point CDF gives the base case. Induction and marginal atomlessness from Proposition~\ref{prop:half-brownian-fixed-time-laws} then yield \eqref{eq:half-brownian-upper-tail}.

\section{Contour decompositions}\label{sec:contour-decompositions}

We derive residue decompositions of the narrow-wedge joint upper-tail formula \cite{Liu-Zhang25} and the half-Brownian formula in Proposition~\ref{prop:half-brownian-upper-tail}. The argument uses contour deformations and residue computations similar to those in \cite{Liu22,Liu_Ortiz25,Nissim-Zhang26,Baik-Cordaro-Tripathi25,Baik-Tripathi26}. The resulting decompositions allow us to compare the two series in the edge asymptotic analysis.

\subsection{The narrow-wedge joint upper-tail formula}\label{sec:narrow-wedge-upper-tail-formula}

We use the formula in \cite{Liu-Zhang25} for the joint upper-tail probabilities of the narrow-wedge KPZ fixed point at $\mm+1$ space-time points. Fix $\mm\geq1$ and $(x_\ell,t_\ell)\in\realR\times(0,\infty)$, $1\leq \ell\leq\mm+1$, satisfying
\begin{equation*}
(x_1,t_1)\prec\cdots\prec(x_{\mm+1},t_{\mm+1}),
\end{equation*}
together with the height parameters $h_1,\ldots,h_{\mm+1}\in\realR$.
Define the exponential factors and ratios of consecutive factors by
\begin{equation}\begin{split}\label{eq:narrow-wedge-phases-and-increments}
f_\ell(w)&:=\exp\left\{-\frac{t_\ell}{3}w^3+x_\ell w^2+h_\ell w\right\},
\qquad 1\leq \ell\leq\mm+1,\\
F_1(w)&:=f_1(w),
\qquad
F_\ell(w):=\frac{f_\ell(w)}{f_{\ell-1}(w)},
\qquad 2\leq \ell\leq\mm+1.
\end{split}\end{equation}

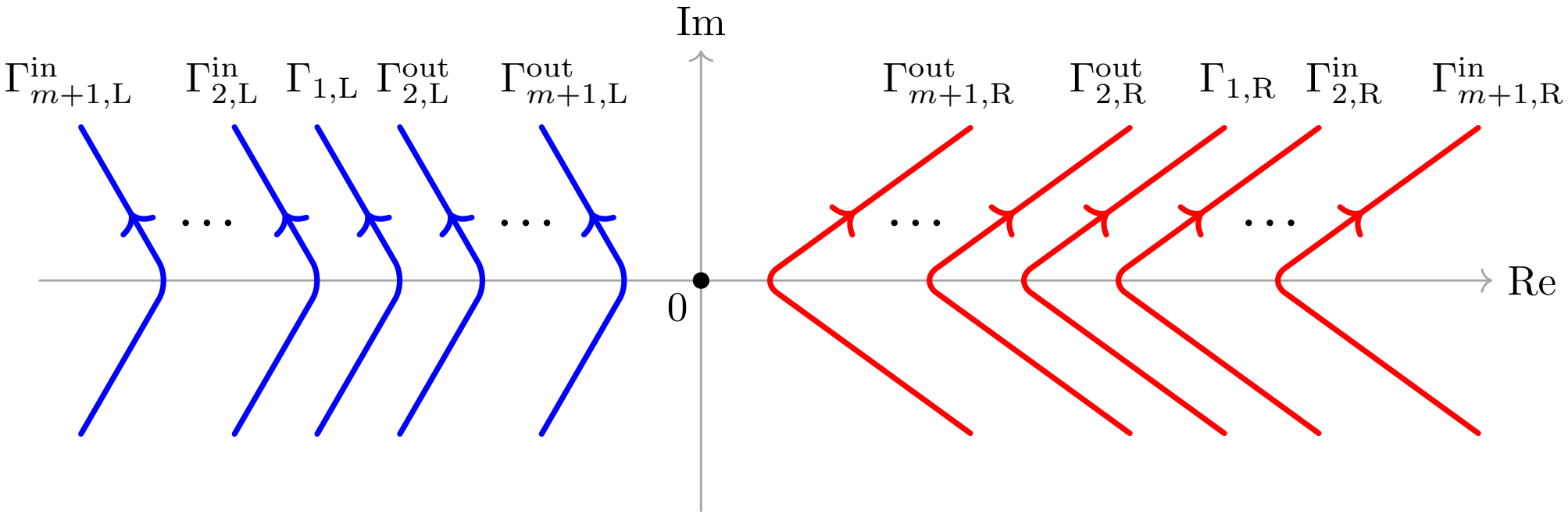
\begin{figure}[htbp]
\centering
\begin{tikzpicture}[
    x=1cm,
    y=1cm,
    line cap=round,
    line join=round,
    contour/.style={
        line width=0.85pt,
        rounded corners=5pt,
        postaction={decorate},
        decoration={
            markings,
            mark=at position 0.72 with {\arrow{>}}
        }
    },
    contour label/.style={anchor=base,inner sep=1pt}
]
    \draw[->,gray!70,line width=0.55pt]
        (-5.55,0) -- (6.75,0) node[right,black] {$\Re$};
    \draw[->,gray!70,line width=0.55pt]
        (0.05,-1.95) -- (0.05,1.95) node[above,black] {$\Im$};
    \node[below left,inner sep=3pt] at (0.05,0) {$0$};
    \fill (0.05,0) circle[radius=2pt];
    \foreach \x in {0.55,1.9,2.7,3.5,4.85} {
        \begin{scope}[xshift=\x cm]
            \draw[red,contour,line width=1.35pt]
                (-36:2.2) -- (0,0) -- (36:2.2);
        \end{scope}
    }
    \foreach \x in {-4.45,-3.15,-2.45,-1.75,-0.55} {
        \begin{scope}[xshift=\x cm]
            \draw[blue,contour,line width=1.35pt]
                (-120:1.50) -- (0,0) -- (120:1.50);
        \end{scope}
    }
    \foreach \x in {-1.4,1.9,4.9,-4.1} {
        \node[contour label] at (\x,0.4) {$\boldsymbol{\cdots}$};
    }
    \foreach \x/\y/\lab in {
        -5.3/1.60/{\Gamma^{\inn}_{\mm+1,\rL}},
        -4/1.60/{\Gamma^{\inn}_{2,\rL}},
        -3.15/1.60/{\Gamma_{1,\rL}},
        -2.38/1.60/{\Gamma^{\out}_{2,\rL}},
        -1.1/1.60/{\Gamma^{\out}_{\mm+1,\rL}}
    } {
        \node[contour label] at (\x,\y) {$\lab$};
    }
    \foreach \x/\y/\lab in {
        2.15/1.60/{\Gamma^{\out}_{\mm+1,\rR}},
        3.5/1.60/{\Gamma^{\out}_{2,\rR}},
        4.6/1.60/{\Gamma_{1,\rR}},
        5.5/1.60/{\Gamma^{\inn}_{2,\rR}},
        6.8/1.60/{\Gamma^{\inn}_{\mm+1,\rR}}
    } {
        \node[contour label] at (\x,\y) {$\lab$};
    }
\end{tikzpicture}
\caption{Finite-$L$ $\Gamma$-contours for $D_{(\kk_1,\ldots,\kk_{\mm+1})}$.}
\label{fig:narrow-wedge-contours}
\end{figure}

For $\kk_1,\ldots,\kk_{\mm+1}\geq1$, write
\begin{equation*}
U^{(\ell)}=(u^{(\ell)}_1,\ldots,u^{(\ell)}_{\kk_\ell}),
\qquad
V^{(\ell)}=(v^{(\ell)}_1,\ldots,v^{(\ell)}_{\kk_\ell}),
\qquad 1\leq \ell\leq\mm+1,
\end{equation*}
and set $U^{(\mm+2)}=V^{(\mm+2)}=\emptyset$ at the auxiliary index $\mm+2$, beyond the last actual layer.
\begin{defn}
For $z_1,\ldots,z_\mm\in\complexC\setminus\{0,1\}$, define
\begin{equation*}\begin{split}
D_{(\kk_1,\ldots,\kk_{\mm+1})}(z_1,\ldots,z_\mm)&:=\prod_{\ell=1}^{\mm}(1-z_\ell)^{\kk_\ell}(1-z_\ell^{-1})^{\kk_{\ell+1}}\\
&\quad\cdot\prod_{i_1=1}^{\kk_1}
\int_{\Gamma_{1,\rL}}\frac{\rd u^{(1)}_{i_1}}{2\pi\rmi}
\prod_{\ell=2}^{\mm+1}\prod_{i_\ell=1}^{\kk_\ell}
\left(
\frac{1}{1-z_{\ell-1}}
\int_{\Gamma^{\inn}_{\ell,\rL}}
\frac{\rd u^{(\ell)}_{i_\ell}}{2\pi\rmi}
-
\frac{z_{\ell-1}}{1-z_{\ell-1}}
\int_{\Gamma^{\out}_{\ell,\rL}}
\frac{\rd u^{(\ell)}_{i_\ell}}{2\pi\rmi}
\right)\\
&\quad\cdot\prod_{i_1=1}^{\kk_1}
\int_{\Gamma_{1,\rR}}\frac{\rd v^{(1)}_{i_1}}{2\pi\rmi}
\prod_{\ell=2}^{\mm+1}\prod_{i_\ell=1}^{\kk_\ell}
\left(
\frac{1}{1-z_{\ell-1}}
\int_{\Gamma^{\inn}_{\ell,\rR}}
\frac{\rd v^{(\ell)}_{i_\ell}}{2\pi\rmi}
-
\frac{z_{\ell-1}}{1-z_{\ell-1}}
\int_{\Gamma^{\out}_{\ell,\rR}}
\frac{\rd v^{(\ell)}_{i_\ell}}{2\pi\rmi}
\right)\\
&\quad\cdot \cC(V^{(1)};U^{(1)}) \cdot
\prod_{\ell=1}^{\mm+1}
\cC(U^{(\ell)}\sqcup V^{(\ell+1)};
V^{(\ell)}\sqcup U^{(\ell+1)}) \cdot
\prod_{\ell=1}^{\mm+1}\prod_{i_\ell=1}^{\kk_\ell}
\frac{F_\ell(u^{(\ell)}_{i_\ell})}{F_\ell(v^{(\ell)}_{i_\ell})}.
\end{split}\end{equation*}
The $\Gamma$-contours have the same nesting, positive separation, orientations, and tail directions as the $\Sigma$-contours in Section~\ref{sec:half-brownian-distribution-formulas}, with $\mm$ replaced by $\mm+1$. Here the left and right contours lie in the left and right half-planes, respectively; see Figure~\ref{fig:narrow-wedge-contours}.
\end{defn}

Set
\begin{equation}\label{eq:narrow-wedge-integrated-coefficient}
\hatD_{(\kk_1,\ldots,\kk_{\mm+1})}
:=(-1)^{\mm+1}
\left[\prod_{\ell=1}^{\mm}\oint_{>1}
\frac{\rd z_\ell}{2\pi\rmi z_\ell(1-z_\ell)}\right]
D_{(\kk_1,\ldots,\kk_{\mm+1})}(z_1,\ldots,z_\mm).
\end{equation}
We will use the following upper-tail formula from \cite{Liu-Zhang25}.
\begin{prop}\cite[Proposition~3.1]{Liu-Zhang25}
For all $(x_1,t_1)\prec\cdots\prec(x_{\mm+1},t_{\mm+1}),
t_\ell>0,$
\begin{equation}\label{eq:narrow-wedge-upper-tail}
\prob\left(\bigcap_{\ell=1}^{\mm+1}
\{\hKPZ(x_\ell,t_\ell)\geq h_\ell\}\right)
=\sum_{\kk_1,\ldots,\kk_{\mm+1}\geq1}
\frac{\hatD_{(\kk_1,\ldots,\kk_{\mm+1})}}
{(\kk_1!\cdots \kk_{\mm+1}!)^2}.
\end{equation}
\end{prop}
The contour integrals and the series in
\eqref{eq:narrow-wedge-upper-tail} are absolutely convergent.

\subsection{Decomposition of the narrow-wedge contour integrals}
We use the notation
\begin{equation*}
U^{(\ell)}_{[a,b]}:=
\begin{cases}
(u_a^{(\ell)},\ldots,u_b^{(\ell)}),&a\leq b,\\
\emptyset,&a>b.
\end{cases}
\end{equation*}
The notation $V^{(\ell)}_{[a,b]}$ is defined similarly. For a tuple $W=(w_1,\ldots,w_q)$, write
\[
F_\ell(W):=\prod_{i=1}^q F_\ell(w_i),
\qquad
G_\ell(W):=\prod_{i=1}^q G_\ell(w_i).
\]
For $A\in\{\rL,\rR\}$, abbreviate weighted contour integration by
\begin{equation}\begin{split}\label{eq:weighted-contour-integration}
\int \varphi(w)\,\rd\mu_{1,A}(w)
&:=\int_{\Gamma_{1,A}}\varphi(w)\frac{\rd w}{2\pi\rmi},\\
\int \varphi(w)\,\rd\mu_{\ell,A}(w)
&:=\frac{1}{1-z_{\ell-1}}\int_{\Gamma_{\ell,A}^{\inn}}\varphi(w)\frac{\rd w}{2\pi\rmi}
-\frac{z_{\ell-1}}{1-z_{\ell-1}}\int_{\Gamma_{\ell,A}^{\out}}\varphi(w)\frac{\rd w}{2\pi\rmi},
\qquad \ell\geq2.
\end{split}\end{equation}
For tuples, $\rd\mu_{\ell,A}(W):=\prod_{i=1}^q\rd\mu_{\ell,A}(w_i)$ denotes product integration. If $W=\emptyset$, no integration is performed and the integrand is unchanged. For $\ell\geq2$, we suppress the dependence of $\rd\mu_{\ell,A}$ on $z_{\ell-1}$ in the notation. Each integration operator acts on all factors to its right.

\begin{defn}\label{defn:r-coefficients}
Let $\bk=(\kk_1,\ldots,\kk_{\mm+1})$, $1\leq\rr\leq\mm$,
$0\leq s\leq\min(\kk_\rr-1,\kk_{\rr+1})$, and
$0\leq s'\leq\min(\kk_\rr,\kk_{\rr+1})$. Set
\[
\widetilde U^{(\rr+1)}
:=U^{(\rr)}_{[1,s]}\sqcup U^{(\rr+1)}_{[s+1,\kk_{\rr+1}]},
\qquad
\widetilde V^{(\rr+1)}
:=V^{(\rr)}_{[1,s']}\sqcup V^{(\rr+1)}_{[s'+1,\kk_{\rr+1}]},
\]
and leave all other tuples unchanged:
$\widetilde U^{(\ell)}=U^{(\ell)}$ and $\widetilde V^{(\ell)}=V^{(\ell)}$ for $\ell\ne\rr+1$.
Define the scalar coefficient
\begin{equation}\label{eq:r-prefactor}
c_{(\bk;\,\rr)}(s,s')
:=(-1)^{\mm+1+\kk_\rr(s+s')+ss'+s'}
\binom{\kk_\rr}s\binom{\kk_{\rr+1}}s s!\cdot
\binom{\kk_\rr}{s'}\binom{\kk_{\rr+1}}{s'}s'!\cdot
\binom{\kk_\rr+\kk_{\rr+1}-s-s'-1}{\kk_{\rr+1}-s'}.
\end{equation}
Then set
\begin{equation}\begin{split}\label{eq:r-coefficient}
R_{(\bk;\,\rr)}(s,s')
&:=c_{(\bk;\,\rr)}(s,s')\cdot
\left[\prod_{\substack{1\leq\ell\leq\mm\\\ell\ne\rr}}
\oint_{>1}\frac{\rd z_\ell}{2\pi\rmi z_\ell(1-z_\ell)}\right]
\prod_{\substack{1\leq\ell\leq\mm\\\ell\ne\rr}}
(1-z_\ell)^{\kk_\ell}(1-z_\ell^{-1})^{\kk_{\ell+1}}\\
&\quad\cdot
\left[\prod_{\substack{1\leq\ell\leq\mm+1\\\ell\ne\rr+1}}
\int \rd\mu_{\ell,\rL}(U^{(\ell)})\,\rd\mu_{\ell,\rR}(V^{(\ell)})\right]
\prod_{i=s+1}^{\kk_{\rr+1}}\int_{\Gamma_{\rr+1,\rL}^{\out}}\frac{\rd u_i^{(\rr+1)}}{2\pi\rmi}
\prod_{i=s'+1}^{\kk_{\rr+1}}\int_{\Gamma_{\rr+1,\rR}^{\inn}}\frac{\rd v_i^{(\rr+1)}}{2\pi\rmi}\\
&\quad\cdot
\cC(V^{(1)};U^{(1)})
\cC\left(
U^{(\rr)}_{[s+1,\kk_\rr]}\sqcup V^{(\rr+1)}_{[s'+1,\kk_{\rr+1}]};
V^{(\rr)}_{[s'+1,\kk_\rr]}\sqcup U^{(\rr+1)}_{[s+1,\kk_{\rr+1}]}
\right)\\
&\quad\cdot
\prod_{\substack{1\leq\ell\leq\mm+1\\\ell\ne\rr}}
\cC\left(\widetilde U^{(\ell)}\sqcup\widetilde V^{(\ell+1)};
\widetilde V^{(\ell)}\sqcup\widetilde U^{(\ell+1)}\right)
\prod_{\ell=1}^{\mm+1}
\frac{F_\ell(\widetilde U^{(\ell)})}{F_\ell(\widetilde V^{(\ell)})}.
\end{split}\end{equation}
The entries of $U^{(\rr)}_{[1,s]}$ and $V^{(\rr)}_{[1,s']}$ inserted into $\widetilde U^{(\rr+1)}$ and $\widetilde V^{(\rr+1)}$ retain their original layer-$\rr$ contour integrations.
\end{defn}

We call the Cauchy factor coupling layers $\rr$ and $\rr+1$ the \emph{central determinant}. In \eqref{eq:r-coefficient}, it is reduced to the determinant involving only the unmatched variables in these two layers:
\[
\cC\left(
U^{(\rr)}_{[s+1,\kk_\rr]}\sqcup V^{(\rr+1)}_{[s'+1,\kk_{\rr+1}]};
V^{(\rr)}_{[s'+1,\kk_\rr]}\sqcup U^{(\rr+1)}_{[s+1,\kk_{\rr+1}]}
\right).
\]
The other Cauchy factors are $\cC(V^{(1)};U^{(1)})$ and the factors indexed by $\ell\ne\rr$ in the Cauchy product in \eqref{eq:r-coefficient}. We call the minors obtained from these other factors by deleting rows and columns in subsequent residue calculations the \emph{outer minors}, including the unchanged factors when no rows or columns are deleted.

The exponential factors in \eqref{eq:r-coefficient} are $F_\ell$ and $F_\rr F_{\rr+1}$ for the $u$-variables, and their reciprocals for the $v$-variables. The ordering of the narrow-wedge space-time points and \eqref{eq:narrow-wedge-phases-and-increments} imply super-exponential decay along the contours. Since the contours appearing in a common Cauchy denominator are separated, the integrals defining $R_{(\bk;\,\rr)}(s,s')$ converge absolutely.

The residue terms $R$ give the following residue decomposition of $\hatD$.

\begin{lm}\label{lm:narrow-wedge-residue-decomposition}
For $1\leq\rr\leq\mm$ and $\kk_1,\ldots,\kk_{\mm+1}\geq1$, we have
\begin{equation}\label{eq:narrow-wedge-residue-decomposition}
\hatD_{(\kk_1,\ldots,\kk_{\mm+1})}
=
\sum_{s=0}^{\min(\kk_{\rr}-1,\kk_{\rr+1})}
\sum_{s'=0}^{\min(\kk_{\rr},\kk_{\rr+1})}
R_{(\bk;\,\rr)}(s,s').
\end{equation}
\end{lm}

\begin{proof}
In \eqref{eq:narrow-wedge-integrated-coefficient}, we deform $\Gamma_{\rr+1,\rL}^{\inn}$ to $\Gamma_{\rr+1,\rL}^{\out}$ and $\Gamma_{\rr+1,\rR}^{\out}$ to $\Gamma_{\rr+1,\rR}^{\inn}$; see Figure~\ref{fig:narrow-wedge-residue-deformation}. The left deformation crosses only the Cauchy poles at the $U^{(\rr)}$-variables, which lie on $\Gamma_{\rr,\rL}^{\inn}$ or $\Gamma_{\rr,\rL}^{\out}$, or on $\Gamma_{1,\rL}$ when $\rr=1$. The right deformation crosses only the poles $v_j^{(\rr+1)}=v_i^{(\rr)}$.

Let $s$ and $s'$ be the numbers of residues taken in the $U^{(\rr+1)}$- and $V^{(\rr+1)}$-variables, respectively. We keep the matched variables in layer $\rr$ on their original contours and include the remaining integrals in $R_{(\bk;\,\rr)}(s,s')$. In particular, these retained variables are not moved across any additional poles created by the residues.

\begin{figure}[htpb]
\centering
\resizebox{0.90\textwidth}{!}{\begin{tikzpicture}[
    x=1cm,
    y=1cm,
    line cap=round,
    line join=round,
    contour/.style={
        line width=0.85pt,
        rounded corners=5pt,
        postaction={decorate},
        decoration={
            markings,
            mark=at position 0.72 with {\arrow{>}}
        }
    },
    oldcontour/.style={blue,dashed,contour},
    deform/.style={
        -{Stealth[length=2.4mm,width=1.5mm]},
        gray,
        line width=0.9pt
    },
    contour label/.style={anchor=base,inner sep=1pt}
]
    \draw[->,gray!70,line width=0.55pt]
        (-6.85,0) -- (8.25,0) node[right,black] {$\Re$};
    \draw[->,gray!70,line width=0.55pt]
        (0,-1.95) -- (0,1.95) node[above,black] {$\Im$};
    \node[below left,inner sep=3pt] at (0,0) {$0$};
    \fill (0,0) circle[radius=2pt];
    \begin{scope}[xshift=2cm]
        \draw[red,dashed,contour,line width=1.35pt]
            (-36:2.2) -- (0,0) -- (36:2.2);
    \end{scope}
    \foreach \x in {0.5,3,4.6,5.6,7.2} {
        \begin{scope}[xshift=\x cm]
            \draw[red,contour,line width=1.35pt]
                (-36:2.2) -- (0,0) -- (36:2.2);
        \end{scope}
    }
    \begin{scope}[xshift=-4.4cm]
        \draw[oldcontour,line width=1.35pt]
            (-120:1.50) -- (0,0) -- (120:1.50);
    \end{scope}
    \foreach \x in {-0.6,-1.9,-2.5,-3.8,-5.7} {
        \begin{scope}[xshift=\x cm]
            \draw[blue,contour,line width=1.35pt]
                (-120:1.50) -- (0,0) -- (120:1.50);
        \end{scope}
    }
    \foreach \x in {-1.5,-3.4,1.9,4.5,7.1,-5.3} {
        \node[contour label] at (\x,0.4) {$\boldsymbol{\cdots}$};
    }
    \foreach \x/\y/\lab in {
        -6.6/1.60/{\Gamma^{\inn}_{\mm+1,\rL}},
        -5.30/1.60/{\Gamma^{\inn}_{\rr+1,\rL}},
        -4.4/1.60/{\Gamma^{\inn}_{\rr,\rL}},
        -3.30/1.60/{\Gamma^{\out}_{\rr,\rL}},
        -2.30/1.60/{\Gamma^{\out}_{\rr+1,\rL}},
        -1/1.60/{\Gamma^{\out}_{\mm+1,\rL}}
    } {
        \node[contour label] at (\x,\y) {$\lab$};
    }
    \foreach \x/\y/\lab in {
        2.1/1.60/{\Gamma^{\out}_{\mm+1,\rR}},
        3.8/1.60/{\Gamma^{\out}_{\rr+1,\rR}},
        5/1.60/{\Gamma^{\out}_{\rr,\rR}},
        6.40/1.60/{\Gamma^{\inn}_{\rr,\rR}},
        7.60/1.60/{\Gamma^{\inn}_{\rr+1,\rR}},
        9.2/1.60/{\Gamma^{\inn}_{\mm+1,\rR}},
    } {
        \node[contour label] at (\x,\y) {$\lab$};
    }
    \draw[deform,bend left=20]
        (-4.6,0.63) to (-2.3,0.60);
    \draw[deform,bend left=20]
        (2.95,0.60) to (6.35,0.68);
\end{tikzpicture}
}
\caption{Deformation of $\Gamma_{\rr+1,\rL}^{\inn}$ to $\Gamma_{\rr+1,\rL}^{\out}$ and of $\Gamma_{\rr+1,\rR}^{\out}$ to $\Gamma_{\rr+1,\rR}^{\inn}$, crossing the layer-$\rr$ contours.}
\label{fig:narrow-wedge-residue-deformation}
\end{figure}
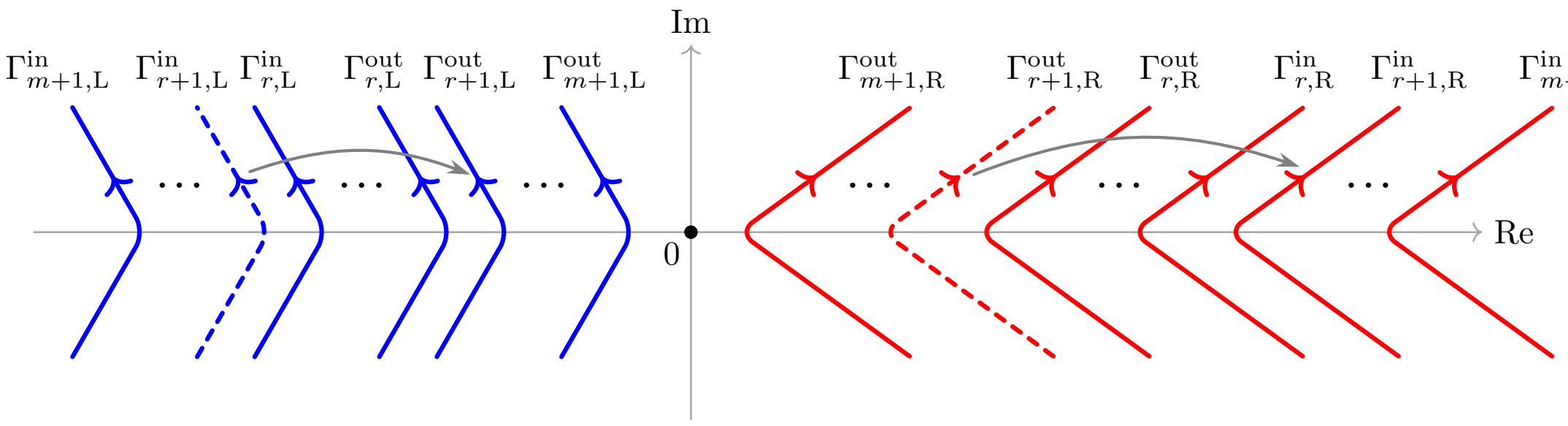

The decay of the exponential factors makes the connecting integrals vanish as the contour truncation tends to infinity. This justifies the deformations.

For $1\leq \ell\leq\kk_{\rr+1}$, we apply the residue theorem
in each of the variables $u_\ell^{(\rr+1)}$ and $v_\ell^{(\rr+1)}$,
keeping the others fixed.
Since all contours are oriented upward, we obtain
\begin{equation*}\begin{split}
&\int_{\Gamma_{\rr+1,\rL}^{\inn}}
\cdots\frac{\rd u_\ell^{(\rr+1)}}{2\pi\rmi}=\int_{\Gamma_{\rr+1,\rL}^{\out}}\cdots\frac{\rd u_\ell^{(\rr+1)}}{2\pi\rmi}-\sum_{i=1}^{\kk_{\rr}}\Res_{u_\ell^{(\rr+1)}=u_i^{(\rr)}}(\cdots),\\
&\int_{\Gamma_{\rr+1,\rR}^{\out}}
\cdots\frac{\rd v_\ell^{(\rr+1)}}{2\pi\rmi}=\int_{\Gamma_{\rr+1,\rR}^{\inn}}
\cdots\frac{\rd v_\ell^{(\rr+1)}}{2\pi\rmi}-\sum_{i=1}^{\kk_{\rr}}
\Res_{v_\ell^{(\rr+1)}=v_i^{(\rr)}}(\cdots).
\end{split}\end{equation*}
Each residue deletes one row and one column from the central Cauchy determinant. Hence each variable in $U^{(\rr)}$ or $V^{(\rr)}$ can be matched at most once in a nonzero term. By symmetry within the four tuples, all matchings with the same $(s,s')$ give equal contributions after integration. The number of matchings is
\begin{equation*}
\binom{\kk_{\rr}}s\binom{\kk_{\rr+1}}s s!
\binom{\kk_{\rr}}{s'}\binom{\kk_{\rr+1}}{s'}s'!.
\end{equation*}
It remains to evaluate the canonical matching
\begin{equation*}
\begin{split}
u_i^{(\rr+1)}&=u_i^{(\rr)},\qquad 1\leq i\leq s,\\
v_i^{(\rr+1)}&=v_i^{(\rr)},\qquad 1\leq i\leq s'.
\end{split}
\end{equation*}

For this matching, a direct calculation gives
\begin{equation*}\begin{split}
&\left(\prod_{i=1}^{s}
\Res_{u_i^{(\rr+1)}=u_i^{(\rr)}}\right)
\left(\prod_{i=1}^{s'}
\Res_{v_i^{(\rr+1)}=v_i^{(\rr)}}\right)
\cC(U^{(\rr)}\sqcup V^{(\rr+1)};
V^{(\rr)}\sqcup U^{(\rr+1)})\\
&\qquad=
(-1)^{\kk_{\rr}(s+s')+ss'+s}
\cC\left(
U^{(\rr)}_{[s+1,\kk_{\rr}]}\sqcup
V^{(\rr+1)}_{[s'+1,\kk_{\rr+1}]};
V^{(\rr)}_{[s'+1,\kk_{\rr}]}\sqcup
U^{(\rr+1)}_{[s+1,\kk_{\rr+1}]}
\right).
\end{split}\end{equation*}
Moving the matched rows and columns to the front gives
$(-1)^{\kk_{\rr}(s+s')}$, and interchanging the two matched
column blocks gives $(-1)^{ss'}$.
The $s$ residues at $u_i^{(\rr+1)}=u_i^{(\rr)}$
contribute $(-1)^s$ through the factors
$(u_i^{(\rr)}-u_i^{(\rr+1)})^{-1}$,
while the residues at $v_i^{(\rr+1)}=v_i^{(\rr)}$
contribute no additional sign.

The $z_{\rr}$-dependence is now confined to
\begin{equation*}\begin{split}
&\oint_{>1}\frac{\rd z_{\rr}}
{2\pi\rmi z_{\rr}(1-z_{\rr})}
(1-z_{\rr})^{\kk_{\rr}}
(1-z_{\rr}^{-1})^{\kk_{\rr+1}}
\left(-\frac{1}{1-z_{\rr}}\right)^s
\left(\frac{z_{\rr}}{1-z_{\rr}}\right)^{s'}\\
&\qquad=\begin{cases}(-1)^{s+s'}
\displaystyle\binom{\kk_{\rr}+\kk_{\rr+1}-s-s'-1}{\kk_{\rr+1}-s'},
&0\leq s\leq \kk_{\rr}-1,\\[2mm]
0,&s=\kk_{\rr}.
\end{cases}
\end{split}\end{equation*}
For $s<\kk_\rr$, the only pole is at $z_\rr=0$, and the binomial expansion gives the stated value. For $s=\kk_\rr$, the integrand is $O(z_\rr^{-2})$ at infinity, so enlarging the circle gives zero.

The unmatched variables in layer $\rr+1$ remain on $\Gamma_{\rr+1,\rL}^{\out}$ and $\Gamma_{\rr+1,\rR}^{\inn}$. Substituting the matched variables gives the tuples $\widetilde U^{(\rr+1)}$ and $\widetilde V^{(\rr+1)}$ in the next Cauchy determinant and the exponential factors, as in Definition~\ref{defn:r-coefficients}. The total sign is
\begin{equation*}
(-1)^{\mm+1}(-1)^{\kk_\rr(s+s')+ss'+s}(-1)^{s+s'}=(-1)^{\mm+1+\kk_\rr(s+s')+ss'+s'}.
\end{equation*}
Multiplying by the number of matchings gives $R_{(\bk;\,\rr)}(s,s')$. Summing over the admissible pairs $(s,s')$ proves \eqref{eq:narrow-wedge-residue-decomposition}.
\end{proof}

\subsection{Decomposition of the half-Brownian contour integrals}
\label{sec:half-brownian-coefficient-decomposition}

Recall the contour integrals $\hatK_{\bn}$ from \eqref{eq:half-brownian-integrated-coefficient} and the empty-layer identity in Lemma~\ref{lm:half-brownian-empty-layer}.

We deform the contours in $\hatK_{\bn}$ from the $\Sigma$-contours to the $\Gamma$-contours. As shown in the proof of Proposition~\ref{prop:half-brownian-residue-decomposition}, every nonzero term produced by this deformation has exactly one residue at zero in each of layers $1,\ldots,p$ and no residues in layers $p+1,\ldots,\mm$, for some $0\leq p\leq\mm$. We denote the total contribution of the terms with this residue pattern by $Q_{\bn}^{(p)}$; in particular, $Q_{\bn}^{(0)}$ is the term with no residues.

\begin{defn}\label{defn:q-coefficients}
Let the $\Gamma$-contours satisfy the conditions in
Section~\ref{sec:narrow-wedge-upper-tail-formula} for $\mm$ layers; see
Figure~\ref{fig:half-brownian-left-contour-deformation}.
Use the convention $G_\ell(W)=\prod_{i=1}^qG_\ell(w_i)$ for a tuple $W=(w_1,\ldots,w_q)$ and the weighted integrals in \eqref{eq:weighted-contour-integration}, with
$G_\ell$ as in \eqref{eq:half-brownian-phases-and-increments}.
For $\bn=(\nn_1,\ldots,\nn_\mm)$ with $\nn_\ell\geq1$, define
\begin{equation}\begin{split}\label{eq:q-coefficient-no-residue}
Q^{(0)}_{\bn}
&:=(-1)^\mm
\left[\prod_{\ell=1}^{\mm-1}\oint_{>1}
\frac{\rd z_\ell}{2\pi\rmi z_\ell(1-z_\ell)}\right]
\prod_{\ell=1}^{\mm-1}(1-z_\ell)^{\nn_\ell}(1-z_\ell^{-1})^{\nn_{\ell+1}}\left[\prod_{\ell=1}^{\mm}\int
\rd\mu_{\ell,\rL}(\bxi^{(\ell)})\,\rd\mu_{\ell,\rR}(\vecteta^{(\ell)})\right]\\
&\quad\cdot
\cC(\vecteta^{(1)};\bxi^{(1)})
\prod_{\ell=1}^{\mm}\cC(\bxi^{(\ell)}\sqcup\vecteta^{(\ell+1)};
\vecteta^{(\ell)}\sqcup\bxi^{(\ell+1)})\cdot
\prod_{i=1}^{\nn_1}\frac{\eta_i^{(1)}}{\xi_i^{(1)}}
\prod_{\ell=1}^{\mm}
\frac{G_\ell(\bxi^{(\ell)})}{G_\ell(\vecteta^{(\ell)})}.
\end{split}\end{equation}

For each $0\leq p\leq\mm$, write
\[
\widehat{\bxi}^{(\ell)}
:=
\begin{cases}
\bxi^{(\ell)}_{[1,\nn_\ell-1]},&1\leq\ell\leq p,\\
\bxi^{(\ell)},&p<\ell\leq\mm+1.
\end{cases}
\]
The dependence of $\widehat{\bxi}^{(\ell)}$ on $p$ is suppressed. We retain the convention $\widehat{\bxi}^{(\mm+1)}=\vecteta^{(\mm+1)}=\emptyset$ at the auxiliary index $\mm+1$. Also, $\widehat{\bxi}^{(\ell)}=\emptyset$ when $\ell\leq p$ and $\nn_\ell=1$.
For $1\leq p\leq\mm$, write $0$ for the one-entry tuple $(0)$ and define
\begin{equation}\begin{split}\label{eq:q-coefficient-with-residues}
Q^{(p)}_{\bn}
&:=(-1)^{\mm+\sum_{\ell=2}^{p}\nn_\ell}
\left(\prod_{\ell=1}^{p}\nn_\ell\right)
\left[\prod_{\ell=1}^{\mm-1}\oint_{>1}
\frac{\rd z_\ell}{2\pi\rmi z_\ell(1-z_\ell)}\right]\\
&\quad\cdot
\prod_{\ell=1}^{p-1}\frac{z_\ell}{1-z_\ell}
\prod_{\ell=1}^{\mm-1}(1-z_\ell)^{\nn_\ell}(1-z_\ell^{-1})^{\nn_{\ell+1}}
\left[\prod_{\ell=1}^{\mm}\int
\rd\mu_{\ell,\rL}(\widehat{\bxi}^{(\ell)})\,\rd\mu_{\ell,\rR}(\vecteta^{(\ell)})\right]\\
&\quad\cdot
\cC(\vecteta^{(1)};\widehat{\bxi}^{(1)}\sqcup0)
\cC(\widehat{\bxi}^{(p)}\sqcup0\sqcup\vecteta^{(p+1)};
\vecteta^{(p)}\sqcup\widehat{\bxi}^{(p+1)})\\
&\quad\cdot
\prod_{\substack{1\leq\ell\leq\mm\\\ell\ne p}}
\cC(\widehat{\bxi}^{(\ell)}\sqcup\vecteta^{(\ell+1)};
\vecteta^{(\ell)}\sqcup\widehat{\bxi}^{(\ell+1)})
\frac{\prod_{i=1}^{\nn_1}\eta_i^{(1)}}{\prod_{i=1}^{\nn_1-1}\xi_i^{(1)}}
\prod_{\ell=1}^{\mm}\frac{G_\ell(\widehat{\bxi}^{(\ell)})}{G_\ell(\vecteta^{(\ell)})}.
\end{split}\end{equation}
For $0\leq p\leq\mm$, set $Q^{(p)}_{\bn}:=0$ if some $\nn_\ell\leq0$.
\end{defn}
The integrals defining $Q^{(p)}_{\bn}$ converge absolutely. Indeed, no Cauchy denominator in \eqref{eq:q-coefficient-with-residues} contains two fixed zeros. We choose the $\Gamma$-contours at positive distance from one another and from zero. The Cauchy factors are then bounded for each fixed $\bn$, and \eqref{eq:half-brownian-phases-and-increments} gives the required decay of the exponential factors.

\begin{prop}\label{prop:half-brownian-residue-decomposition}
For every $\nn_1,\ldots,\nn_\mm\geq1$, we have
\begin{equation}\label{eq:half-brownian-residue-decomposition}
\hatK_{(\nn_1,\ldots,\nn_\mm)}
=\sum_{p=0}^{\mm}Q^{(p)}_{\bn}.
\end{equation}
\end{prop}

\begin{proof}
The decay of the exponential factors justifies the contour deformations below.

Starting from \eqref{eq:half-brownian-integrated-coefficient}, we move the left $\Sigma$-contours to the $\Gamma$-contours in the left half-plane. We first move $\Sigma_{\ell,\rL}^{\inn}$ to $\Gamma_{\ell,\rL}^{\inn}$ in the order $\ell=\mm,\ldots,2$. This preserves the nesting, and the integrand is regular at zero in these variables, so no pole is crossed. We also move the right contours to their $\Gamma$-contours without crossing a pole. Next, we move $\Sigma_{1,\rL}$ to $\Gamma_{1,\rL}$, and then $\Sigma_{\ell,\rL}^{\out}$ to $\Gamma_{\ell,\rL}^{\out}$ in the order $\ell=2,\ldots,\mm$. See Figures~\ref{fig:half-brownian-contours} and~\ref{fig:half-brownian-left-contour-deformation}. The only possible poles crossed in these steps are at zero. At each layer, we separate the term with no residue from the residue terms.

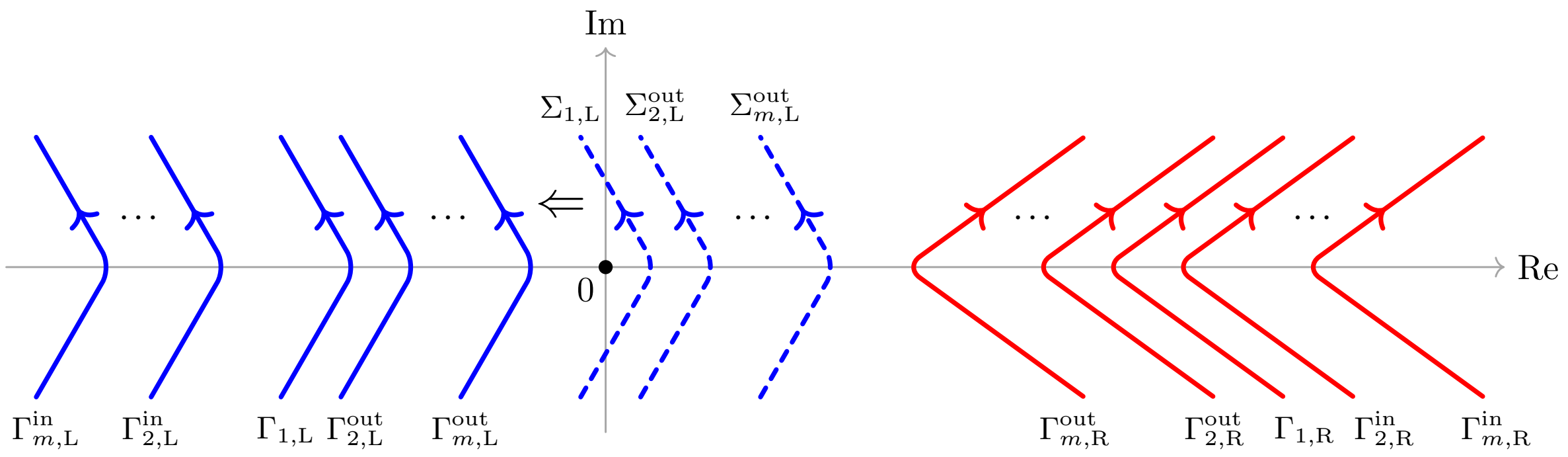
\begin{figure}[htbp]
\centering
\begin{tikzpicture}[
    x=1cm,y=1cm,line cap=round,line join=round,
    contour/.style={line width=1.15pt,rounded corners=5pt,
      postaction={decorate},
      decoration={markings,mark=at position 0.72 with {\arrow{>}}}},
    contour label/.style={inner sep=2pt,font=\small},
    deform/.style={
        -{Stealth[length=2.4mm,width=1.5mm]},
        gray,
        line width=0.9pt
    },
]
\draw[->,gray!70,line width=0.55pt] (-6,0)--(9,0) node[right,black] {$\Re$};
\draw[->,gray!70,line width=0.55pt] (0,-1.65)--(0,2.2) node[above,black] {$\Im$};
\fill (0,0) circle[radius=2pt];
\node[below left,inner sep=3pt] at (0,0) {$0$};
\foreach \x in {-4.95,-3.8,-2.5,-1.9,-0.7}{
  \begin{scope}[xshift=\x cm]
    \draw[blue,contour,line width=1.35pt] (-120:1.50)--(0,0)--(120:1.50);
  \end{scope}
}
\foreach \x in {0.5,1.1,2.3}{
  \begin{scope}[xshift=\x cm]
    \draw[blue,dashed,contour,line width=1.35pt] (-120:1.50)--(0,0)--(120:1.50);
  \end{scope}
}
\foreach \x in {3,4.3,5,5.7,7}{
  \begin{scope}[xshift=\x cm]
    \draw[red,contour,line width=1.35pt] (-36:2.2)--(0,0)--(36:2.2);
  \end{scope}
}
\foreach \x in {-4.65,-1.55,1.5,4.3,7.1}
  \node[font=\scriptsize] at (\x,0.5) {$\boldsymbol{\cdots}$};
\node[contour label,anchor=south] at (-0.35,1.34) {$\Sigma_{1,\rL}$};
\node[contour label,anchor=south] at (0.50,1.34) {$\Sigma_{2,\rL}^{\out}$};
\node[contour label,anchor=south] at (1.6,1.34) {$\Sigma_{\mm,\rL}^{\out}$};
\node[contour label,anchor=south] at (-5.60,-1.9) 
{$\Gamma_{\mm,\rL}^{\inn}$};
\node[contour label,anchor=south] at (-4.55,-1.9) {$\Gamma_{2,\rL}^{\inn}$};
\node[contour label,anchor=south] at (-3.2,-1.87) {$\Gamma_{1,\rL}$};
\node[contour label,anchor=south] at (-2.50,-1.9) {$\Gamma_{2,\rL}^{\out}$};
\node[contour label,anchor=south] at (-1.4,-1.9) {$\Gamma_{\mm,\rL}^{\out}$};
\node[contour label,anchor=south] at (4.7,-1.9) {$\Gamma_{\mm,\rR}^{\out}$};
\node[contour label,anchor=south] at (6.1,-1.9) {$\Gamma_{2,\rR}^{\out}$};
\node[contour label,anchor=south] at (7,-1.87) {$\Gamma_{1,\rR}$};
\node[contour label,anchor=south] at (7.8,-1.9) {$\Gamma_{2,\rR}^{\inn}$};
\node[contour label,anchor=south] at (8.92,-1.9) {$\Gamma_{\mm,\rR}^{\inn}$};
\node[font=\Large,inner sep=0pt]
    at (-0.45,0.60) {$\Leftarrow$};
\end{tikzpicture}
\caption{Deformation of $\Sigma_{1,\rL},\Sigma_{2,\rL}^{\out},\ldots,\Sigma_{\mm,\rL}^{\out}$ (dashed), in this order, to their $\Gamma$-contours (solid). The inner left contours and all right contours have already been deformed without crossing any pole.}
\label{fig:half-brownian-left-contour-deformation}
\end{figure}

In the first layer, the poles at zero come from $(\xi_i^{(1)})^{-1}$. The squared Vandermonde factor in $\bxi^{(1)}$ makes every term with two or more zero residues vanish. If no residue is taken in the first layer, the contours $\Sigma_{\ell,\rL}^{\out}$, $2\leq\ell\leq\mm$, can be moved to $\Gamma_{\ell,\rL}^{\out}$ without crossing a pole, giving $Q^{(0)}_{\bn}$. By symmetry, the contribution with one residue is $\nn_1$ times the term obtained by taking the residue at $\xi_{\nn_1}^{(1)}=0$.

We continue with the contribution in which a residue at zero was taken in the first layer. Suppose one residue has been taken in each of the first $\ell-1$ layers, with $\xi_{\nn_q}^{(q)}=0$ for $1\leq q<\ell$, where $2\leq \ell\leq\mm$. We move $\Sigma_{\ell,\rL}^{\out}$ to $\Gamma_{\ell,\rL}^{\out}$. The fixed zero in layer $\ell-1$ gives poles at $\xi_i^{(\ell)}=0$ in the Cauchy determinant joining the two layers. As before, terms with two or more residues at zero vanish.

If no residue is taken in layer $\ell<\mm$, the Cauchy factor joining layers $\ell$ and $\ell+1$ has no fixed zero and hence no pole at zero in the layer-$(\ell+1)$ variables. The exponential factors are entire, so all later contours can be moved without taking another residue. This gives $Q^{(\ell-1)}_{\bn}$, also when $\ell=\mm$ and there are no later contours to move.

For the single-residue contribution, symmetry allows us to take
the residue at $\xi_{\nn_\ell}^{(\ell)}=0$ and multiply by $\nn_\ell$.
The pole is the entry $-1/\xi_{\nn_\ell}^{(\ell)}$ in row
$\nn_{\ell-1}$ and column $\nn_{\ell-1}+\nn_\ell$ of the adjacent
Cauchy determinant. Hence
\begin{equation}\label{eq:q-cauchy-residue}
\Res_{\xi_{\nn_\ell}^{(\ell)}=0}
\cC\left(
\bxi^{(\ell-1)}_{[1,\nn_{\ell-1}-1]}\sqcup0\sqcup\vecteta^{(\ell)};
\vecteta^{(\ell-1)}\sqcup\bxi^{(\ell)}
\right)=
(-1)^{\nn_\ell+1}
\cC\left(
\bxi^{(\ell-1)}_{[1,\nn_{\ell-1}-1]}\sqcup\vecteta^{(\ell)};
\vecteta^{(\ell-1)}\sqcup\bxi^{(\ell)}_{[1,\nn_\ell-1]}
\right).
\end{equation}
The residue theorem gives this term the contour coefficient
$-z_{\ell-1}/(1-z_{\ell-1})$.
Combining this contour coefficient with the cofactor sign $(-1)^{\nn_\ell+1}$ in \eqref{eq:q-cauchy-residue} and the $\nn_\ell$ choices gives
the factor $(-1)^{\nn_\ell}\nn_\ell z_{\ell-1}/(1-z_{\ell-1})$.
We then proceed to the next layer in this residue term.

Repeating the argument, we obtain one contribution for each $0\leq p\leq\mm$, with a residue in the first $p$ layers and none thereafter. The case $p=\mm$ ends with a residue in the last layer. For $1\leq p\leq\mm$, the coefficient, including the original sign $(-1)^\mm$, is
\begin{equation}\label{eq:q-prefix-coefficient}
(-1)^\mm\nn_1
\prod_{\ell=2}^{p}
\left(
(-1)^{\nn_\ell}\nn_\ell\frac{z_{\ell-1}}{1-z_{\ell-1}}
\right)
=
(-1)^{\mm+\sum_{\ell=2}^{p}\nn_\ell}
\left(\prod_{\ell=1}^{p}\nn_\ell\right)
\prod_{\ell=1}^{p-1}\frac{z_\ell}{1-z_\ell}.
\end{equation}
Using the shortened tuples $\widehat{\bxi}^{(\ell)}$ from Definition~\ref{defn:q-coefficients}, repeated application of \eqref{eq:q-cauchy-residue} leaves the Cauchy factors
\begin{equation}\label{eq:q-cauchy-product}
\cC(\vecteta^{(1)};\widehat{\bxi}^{(1)}\sqcup0)
\cC\left(\widehat{\bxi}^{(p)}\sqcup0\sqcup\vecteta^{(p+1)};
\vecteta^{(p)}\sqcup\widehat{\bxi}^{(p+1)}\right)\cdot
\prod_{\substack{1\leq\ell\leq\mm\\\ell\ne p}}
\cC\left(
\widehat{\bxi}^{(\ell)}\sqcup\vecteta^{(\ell+1)};
\vecteta^{(\ell)}\sqcup\widehat{\bxi}^{(\ell+1)}
\right).
\end{equation}
The first residue removes $(\xi_{\nn_1}^{(1)})^{-1}$, and $G_\ell(0)=1$ for every $\ell$. The remaining integration operators are $\rd\mu_{\ell,\rL}(\widehat{\bxi}^{(\ell)})\,\rd\mu_{\ell,\rR}(\vecteta^{(\ell)})$, and the exponential factors are $\prod_{\ell=1}^{\mm}G_\ell(\widehat{\bxi}^{(\ell)})/G_\ell(\vecteta^{(\ell)})$. Together with \eqref{eq:q-prefix-coefficient}, this gives exactly \eqref{eq:q-coefficient-with-residues}, so the contribution is $Q^{(p)}_{\bn}$. Summing these contributions over $p=1,\ldots,\mm$ and adding $Q^{(0)}_{\bn}$ proves \eqref{eq:half-brownian-residue-decomposition}.
\end{proof}

\subsection{Absolute summability of the residue terms}

\begin{lm}\label{lm:q-coefficient-summability}
For fixed half-Brownian parameters, contours, and circle radii, there is a constant $C>0$, independent of the half-Brownian series indices $\bn$, such that
\begin{equation}\label{eq:q-order-majorant}
\frac{\textstyle\sum_{p=0}^{\mm}|Q^{(p)}_{\bn}|}
{(\nn_1!\cdots\nn_\mm!)^2}
\leq C\prod_{\ell=1}^{\mm}
\frac{C^{\nn_\ell}\nn_\ell^{\nn_\ell}}{(\nn_\ell!)^2},
\qquad \nn_\ell\geq1.
\end{equation}
The same bound holds for the sums of absolute integrals after expanding the contour functionals. In particular,
\begin{equation*}
\sum_{\nn_1,\ldots,\nn_\mm\geq1}
\frac{\textstyle\sum_{p=0}^{\mm}|Q^{(p)}_{\bn}|}
{(\nn_1!\cdots\nn_\mm!)^2}<\infty.
\end{equation*}
\end{lm}

\begin{proof}
We apply the bounds on the exponential factors and Hadamard's inequality as in the proof of Lemma~\ref{lm:half-brownian-absolute-convergence}. For $Q^{(p)}_{\bn}$, the first and last determinants have orders $\nn_1$ and $\nn_\mm$, and each intermediate determinant has order $\nn_\ell+\nn_{\ell+1}$ or $\nn_\ell+\nn_{\ell+1}-1$. The fixed zeros do not change the bounds on the remaining entries. The extra factor $\prod_{\ell=1}^{p}\nn_\ell$ is absorbed into $C^{\nn_1+\cdots+\nn_\mm}$ using $\nn_\ell\leq2^{\nn_\ell}$. Thus \eqref{eq:half-brownian-determinant-bound} also bounds $|Q^{(p)}_{\bn}|$, uniformly for $0\leq p\leq\mm$. Summing over $p$ and using \eqref{eq:two-order-splitting} proves \eqref{eq:q-order-majorant}.

Summability follows from \eqref{eq:factorial-series-summability}.
\end{proof}

\section{Asymptotic analysis and proof of Theorem~\ref{thm:main}}
\label{sec:main-theorem-proof}

We first prove convergence at distinct times. Fix $\mm\geq1$, $\alpha_\ell,\beta_\ell\in\realR$ for $1\leq \ell\leq\mm$, and $0<\tau_1<\cdots<\tau_\mm$. Define
\begin{equation*}
\rr:=\min\bigl(\{\mm+1\}\cup\{\ell:1\leq \ell\leq\mm,\ \tau_\ell\geq1\}\bigr).
\end{equation*}
Then
\begin{equation*}
0<\tau_1<\cdots<\tau_{\rr-1}<1\leq\tau_\rr<\cdots<\tau_\mm,
\end{equation*}
where the inequalities to the left of $1$ are omitted when $\rr=1$, and those to the right are omitted when $\rr=\mm+1$.

We place the conditioning point at layer $\rr$ in the narrow-wedge $(\mm+1)$-point formula \eqref{eq:narrow-wedge-upper-tail} and choose
\begin{equation}\label{eq:narrow-wedge-observation-parameters}
(x_\ell,h_\ell,t_\ell)
=
\begin{cases}
(\tau_\ell\sqrt L+\alpha_\ell,\ -\tau_\ell L-2\alpha_\ell\sqrt L+\beta_\ell,\ \tau_\ell),
&1\leq \ell\leq\rr-1,\\
(0,L,1),&\ell=\rr,\\
(\tau_{\ell-1}\sqrt L+\alpha_{\ell-1},\ -\tau_{\ell-1}L-2\alpha_{\ell-1}\sqrt L+\beta_{\ell-1},\ \tau_{\ell-1}),
&\rr+1\leq \ell\leq\mm+1.
\end{cases}
\end{equation}
The case $\tau_\rr=1$ is allowed. Then $t_\rr=t_{\rr+1}=1$, but $(x_\rr,t_\rr)=(0,1)\prec(\sqrt{L}+\alpha_\rr,1)=(x_{\rr+1},t_{\rr+1})$ for sufficiently large $L$. We use the narrow-wedge parameters \eqref{eq:narrow-wedge-observation-parameters} throughout this section and write $\bk=(\kk_1,\ldots,\kk_{\mm+1})$ for the indices of the narrow-wedge series.

By \eqref{eq:rescaled-edge-field} and \eqref{eq:narrow-wedge-observation-parameters},
\begin{equation*}
\bigcap_{\ell=1}^{\mm}
\left\{\mathcal H_L^{\mathrm{edge}}(\alpha_\ell,\tau_\ell)\geq\beta_\ell\right\}
=
\bigcap_{\substack{1\leq \ell\leq\mm+1\\ \ell\neq\rr}}
\left\{\hKPZ(x_\ell,t_\ell)\geq h_\ell\right\}.
\end{equation*}
Set
\begin{equation*}
\delta_L:=\prob\bigl(\hKPZ(0,1)\geq L\bigr).
\end{equation*}
Applying \eqref{eq:narrow-wedge-upper-tail}, we obtain
\begin{equation}\label{eq:conditional-upper-tail-series}
\prob\left(
\bigcap_{\ell=1}^{\mm}
\left\{\mathcal H_L^{\mathrm{edge}}(\alpha_\ell,\tau_\ell)\geq\beta_\ell\right\}
\,\middle|\,
\hKPZ(0,1)\geq L
\right)=
\frac{1}{\delta_L}
\sum_{\kk_1,\ldots,\kk_{\mm+1}\geq1}
\frac{\hatD_{\bk}}{(\kk_1!\cdots\kk_{\mm+1}!)^2}.
\end{equation}
The random variable $\hKPZ(0,1)$ has the GUE Tracy--Widom distribution. The representation in \cite{Tracy-Widom94} yields the upper-tail estimate below; see \cite{Baik-Buckingham-DiFranco08} for the explicit asymptotic expansion:
\begin{equation}\label{eq:gue-tracy-widom-upper-tail}
\delta_L
=
\frac{1}{16\pi L^{3/2}}
\exp\left\{-\frac{4}{3}L^{3/2}\right\}\bigl(1+O(L^{-3/2})\bigr).
\end{equation}

We use the contours
\begin{equation*}
\Gamma_{\rL}:=\left\{\rho e^{\pm2\pi\rmi/3}:\rho\geq0\right\},
\qquad
\Gamma_{\rR}:=\left\{\rho e^{\pm\pi\rmi/3}:\rho\geq0\right\}.
\end{equation*}
Both contours are oriented from infinity on the lower ray to the origin and then to infinity on the upper ray.

We first compute the limits for fixed narrow-wedge series indices $\bk$, treating $1\leq\rr\leq\mm$ and $\rr=\mm+1$ separately. We then prove a summable bound, uniform in $L$, and pass to the limit in \eqref{eq:conditional-upper-tail-series}. Reindexing the surviving terms by the half-Brownian series indices $\bn$ and summing over the residue index $p$ will give $\rT_\mm(\bbeta;(\alpha_1,\tau_1),\ldots,(\alpha_\mm,\tau_\mm))$.

\subsection{Asymptotics for fixed narrow-wedge series indices}\label{sec:fixed-order-asymptotics}

We first consider $1\leq\rr\leq\mm$. Recall from \eqref{eq:narrow-wedge-residue-decomposition} that
\begin{equation*}
\hatD_{\bk}
=\sum_{s=0}^{\min(\kk_\rr-1,\kk_{\rr+1})}
\sum_{s'=0}^{\min(\kk_\rr,\kk_{\rr+1})}
R_{(\bk;\,\rr)}(s,s').
\end{equation*}

Fix a term $R_{(\bk;\,\rr)}(s,s')$. We first translate all the $u$- and $v$-contours in \eqref{eq:r-coefficient} by $\sqrt L$ to the right. Moving them together preserves their relative positions, so no pole is crossed. The connecting integrals vanish by cubic exponential decay, or by quadratic decay for the equal-time increment when $\tau_\rr=1$.

We change variables accordingly:
\begin{equation}\label{eq:narrow-wedge-shifted-variables}
u_i^{(\ell)}=\sqrt L+\xi_i^{(\ell)},
\qquad
v_i^{(\ell)}=\sqrt L+\eta_i^{(\ell)}.
\end{equation}
The $\xi_i^{(\ell)}$- and $\eta_i^{(\ell)}$-contours are now the original $u_i^{(\ell)}$- and $v_i^{(\ell)}$-contours, respectively. All differences between integration variables are unchanged, so the Cauchy determinants retain their form.

Recall that the exponential weights in the narrow-wedge and half-Brownian integrands are $F_\ell=f_\ell/f_{\ell-1}$ and $G_\ell=g_\ell/g_{\ell-1}$, respectively, with $f_0=g_0=1$. We first describe the translation of the functions $f_\ell$ entering these ratios. Multiplicative constants independent of the integration variable cancel between the numerator and denominator weights in each layer. By the choice of narrow-wedge parameters in \eqref{eq:narrow-wedge-observation-parameters}, we have
\begin{equation}\label{eq:translated-narrow-wedge-phase-ratios}
\frac{f_\ell(\xi+\sqrt L)}{f_\ell(\eta+\sqrt L)}
=
\begin{cases}
\displaystyle\frac{g_\ell(\xi)}{g_\ell(\eta)},
&1\leq \ell\leq\rr-1,\\[8pt]
\displaystyle
\exp\left\{-\sqrt L(\xi^2-\eta^2)-\frac{\xi^3-\eta^3}{3}\right\},
&\ell=\rr,\\[6pt]
\displaystyle\frac{g_{\ell-1}(\xi)}{g_{\ell-1}(\eta)},
&\rr+1\leq \ell\leq\mm+1,
\end{cases}
\end{equation}
Here $f_\ell$ has the narrow-wedge parameters $(x_\ell,h_\ell,t_\ell)$, while $g_\ell$ has the half-Brownian parameters $(\alpha_\ell,\beta_\ell,\tau_\ell)$. Thus the translated integral retains the form of \eqref{eq:r-coefficient}, with $u,v$ replaced by $\xi,\eta$. Inside the ratios defining $F_\ell$, we replace $f_\ell(\,\cdot+\sqrt L)$ by $g_\ell$ for $\ell<\rr$ or $g_{\ell-1}$ for $\ell>\rr$, up to the cancelling constants just described. The factors involving $f_\rr$ become
\begin{equation}\label{eq:translated-conditioning-phase}
\frac{
\textstyle\prod_{i=s+1}^{\kk_\rr}f_\rr(\xi_i^{(\rr)}+\sqrt L)
\prod_{i=s'+1}^{\kk_{\rr+1}}f_\rr(\eta_i^{(\rr+1)}+\sqrt L)
}{
\textstyle\prod_{i=s'+1}^{\kk_\rr}f_\rr(\eta_i^{(\rr)}+\sqrt L)
\prod_{i=s+1}^{\kk_{\rr+1}}f_\rr(\xi_i^{(\rr+1)}+\sqrt L)
}, \qquad \text{where }f_\rr(w)=\exp\{-\frac{1}{3}w^3+Lw\}.
\end{equation}

The next two lemmas identify the leading terms after division by $\delta_L$ and compute their limits.

\begin{lm}\label{lm:nonleading-r-coefficients}
Fix integers $\kk_1,\ldots,\kk_{\mm+1}\geq1$ and $s,s'$ satisfying
\begin{equation*}
0\leq s\leq\min(\kk_\rr-1,\kk_{\rr+1}),
\qquad
0\leq s'\leq\min(\kk_\rr,\kk_{\rr+1}).
\end{equation*}
Then
\begin{equation*}
\lim_{L\to\infty}\frac{R_{(\bk;\,\rr)}(s,s')}{\delta_L}=0,
\end{equation*}
unless
\begin{equation}\label{eq:leading-r-configurations}
(\kk_\rr,\kk_{\rr+1};s,s')
=(\kk+1,\kk;\kk,\kk)
\quad\text{or}\quad
(\kk,\kk;\kk-1,\kk)
\end{equation}
for some integer $\kk\geq1$.
\end{lm}

Recall $Q^{(p)}_\bn$ from Definition~\ref{defn:q-coefficients}, with the convention $Q^{(p)}_\bn=0$ whenever some $\nn_\ell\leq0$.

\begin{lm}\label{lm:leading-r-asymptotics}
Fix integers $\kk_1,\ldots,\kk_{\mm+1}\geq1$.
\begin{enumerate}[(i)]
\item If $(\kk_\rr,\kk_{\rr+1})=(\kk+1,\kk)$ for some integer $\kk\geq1$, then
\begin{equation}\label{eq:leading-r-first-limit}
\lim_{L\to\infty}\frac{R_{(\bk;\,\rr)}(\kk,\kk)}{\delta_L((\kk+1)!)^2}=\sum_{p=0}^{\rr-1}
\left(\prod_{\ell=p+1}^{\rr-1}\kk_\ell^2\right)
Q^{(p)}_{(\kk_1,\ldots,\kk_p,\kk_{p+1}-1,\ldots,\kk_{\rr-1}-1,\kk,\kk_{\rr+2},\ldots,\kk_{\mm+1})}.
\end{equation}

\item If $\kk_\rr=\kk_{\rr+1}=\kk\geq1$, then
\begin{equation}\label{eq:leading-r-second-limit}
\lim_{L\to\infty}\frac{R_{(\bk;\,\rr)}(\kk-1,\kk)}{\delta_L(\kk!)^2}
=\sum_{p=\rr}^{\mm}
Q^{(p)}_{(\kk_1,\ldots,\kk_{\rr-1},\kk,\kk_{\rr+2},\ldots,\kk_{\mm+1})}.
\end{equation}
\end{enumerate}
\end{lm}

Combining Lemmas~\ref{lm:narrow-wedge-residue-decomposition}, \ref{lm:nonleading-r-coefficients}, and~\ref{lm:leading-r-asymptotics} gives the following proposition.
\begin{prop}\label{prop:narrow-wedge-fixed-order-asymptotics}
For $1\leq\rr\leq\mm$ and fixed integers $\kk_1,\ldots,\kk_{\mm+1}\geq1$, we have, as $L\to\infty$,
\begin{equation*}
\begin{split}&\frac{\hatD_{(\kk_1,\ldots,\kk_{\mm+1})}}{\delta_L}\\
&\quad=\begin{cases}\displaystyle(\kk_\rr!)^2\sum_{p=0}^{\rr-1}\left(\prod_{j=p+1}^{\rr-1}\kk_j^2\right)Q^{(p)}_{(\kk_1,\ldots,\kk_p,\kk_{p+1}-1,\ldots,\kk_{\rr-1}-1,\kk_{\rr+1},\kk_{\rr+2},\ldots,\kk_{\mm+1})}+o(1),&\kk_\rr=\kk_{\rr+1}+1,\\
\displaystyle(\kk_\rr!)^2\sum_{p=\rr}^{\mm}Q^{(p)}_{(\kk_1,\ldots,\kk_{\rr-1},\kk_\rr,\kk_{\rr+2},\ldots,\kk_{\mm+1})}+o(1),&\kk_\rr=\kk_{\rr+1},\\
o(1),&\text{otherwise}.
\end{cases}
\end{split}\end{equation*}
\end{prop}

\subsubsection{Proof of Lemma~\ref{lm:nonleading-r-coefficients}}

\begin{figure}[htpb]
\centering
\resizebox{0.90\textwidth}{!}{\begin{tikzpicture}[
    x=1cm,
    y=1cm,
    line cap=round,
    line join=round,
    contour/.style={
        line width=0.85pt,
        rounded corners=5pt,
        postaction={decorate},
        decoration={
            markings,
            mark=at position 0.72 with {\arrow{>}}
        }
    },
    oldcontour/.style={blue,dashed,contour},
    deform/.style={
        -{Stealth[length=2.4mm,width=1.5mm]},
        gray,
        line width=0.9pt
    },
    contour label/.style={anchor=base,inner sep=1pt}
]
    \draw[->,gray!70,line width=0.55pt]
        (-12,0) -- (2,0) node[right,black] {$\Re$};
    \draw[->,gray!70,line width=0.55pt]
        (1.35,-1.95) -- (1.35,1.95) node[above,black] {$\Im$};
    \node[below right,inner sep=1.5pt] at (1.4,0) {$0$};
    
    \begin{scope}[xshift=-2.8cm]
        \draw[oldcontour,line width=1.35pt]
            (-120:1.50) -- (0,0) -- (120:1.50);
    \end{scope}
    \begin{scope}[xshift=-6.4cm]
        \draw[oldcontour,line width=1.35pt]
            (-120:1.50) -- (0,0) -- (120:1.50);
    \end{scope}
    \begin{scope}[xshift=-2cm]
        \draw[oldcontour,line width=1.35pt]
            (-120:1.50) -- (0,0) -- (120:1.50);
    \end{scope}
    \begin{scope}[xshift=-11 cm]
            \draw[black,contour,line width=1.35pt]
                (-120:1.50) -- (0,0) -- (120:1.50);
        \end{scope}
    \fill (-11.05,0) circle[radius=2pt];
    \begin{scope}[xshift=1.4 cm]
            \draw[black,contour,line width=1.35pt]
                (-120:1.50) -- (0,0) -- (120:1.50);
        \end{scope}
    \fill (1.35,0) circle[radius=2pt];
    \node[below right,inner sep=2.5pt] at (-11.05,0) {$\omega_{\rr-1}$};
    \foreach \x in {0.2,-1.2,-3.6,-5.6,-7.6,-9} {
        \begin{scope}[xshift=\x cm]
            \draw[blue,contour,line width=1.35pt]
                (-120:1.50) -- (0,0) -- (120:1.50);
        \end{scope}
    }
    \foreach \x in {-0.7,-4.8,-8.6} {
        \node[contour label] at (\x,0.4) {$\boldsymbol{\cdots}$};
    }
    \foreach \x/\y/\lab in {
        -11.8/1.60/{\omega_{\rr-1}+\Gamma_{\rL}},
        -9.6/1.60/{\Gamma^{\inn}_{\mm+1,\rL}},
        -8.20/1.60/{\Gamma^{\inn}_{\rr+2,\rL}},
        -7.1/1.60/{\Gamma^{\inn}_{\rr,\rL}},
        -6/1.60/{\Gamma^{\inn}_{\rr-1,\rL}},
        -4.50/1.60/{\Gamma^{\out}_{\rr-1,\rL}},
        -3.50/1.60/{\Gamma^{\out}_{\rr,\rL}},
        -2.55/1.60/{\Gamma^{\out}_{\rr+1,\rL}},
        -1.50/1.60/{\Gamma^{\out}_{\rr+2,\rL}},
        -0.2/1.60/{\Gamma^{\out}_{\mm+1,\rL}},
        0.7/1.60/{\Gamma_{\rL}}
    } {
        \node[contour label] at (\x,\y) {$\lab$};
    }
    \draw[deform,bend right=12]
        (-6.85,0.63) to (-11.3,0.63);
    \draw[deform,bend right=10]
        (-3.25,0.63) to (-11.3,0.63);
    \draw[deform,bend left=15]
        (-2.25,0.63) to (0.95,0.63);
\end{tikzpicture}
}
\caption{Initial deformation of the unmatched $\xi^{(\rr)}$- and
$\xi^{(\rr+1)}$-contours to $\omega_{\rr-1}+\Gamma_{\rL}$ and $\Gamma_{\rL}$.}
\label{fig:unmatched-left-contour-deformation}
\end{figure}
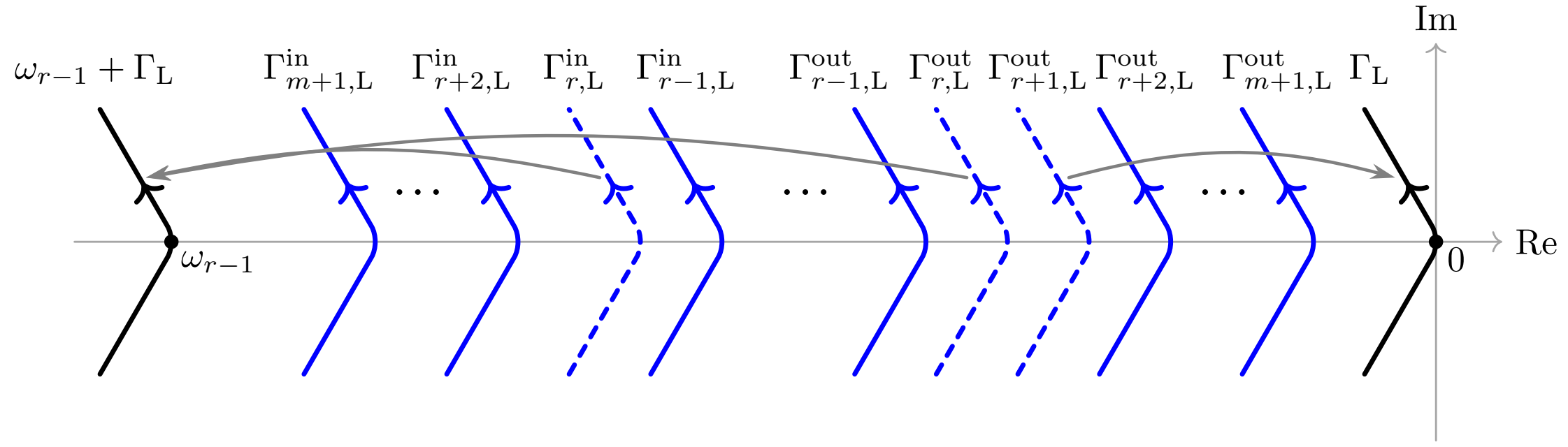

We choose the initial left contours to be horizontal translates of $\Gamma_{\rL}$. We first deform the unmatched $\xi^{(\rr)}$-contours to the left, and then the unmatched $\xi^{(\rr+1)}$-contours to $\Gamma_{\rL}$. Figure~\ref{fig:unmatched-left-contour-deformation} shows the initial step for $\rr\geq2$, with $\omega_{\rr-1}=-2(\sqrt L+\alpha_{\rr-1})/(1-\tau_{\rr-1})$.

For $\xi_i^{(\rr)}$, $s+1\leq i\leq\kk_\rr$, the deformation to $-2\sqrt L+\Gamma_{\rL}$ crosses no pole when $\rr=1$. For $\rr\geq2$, the in-contour crosses no pole, while the out-contour crosses only the poles at $\xi_i^{(\rr)}=\xi_a^{(\rr-1)}$, $1\leq a\leq\kk_{\rr-1}$.

A residue at $\xi_i^{(\rr)}=\xi_a^{(\rr-1)}$ multiplies the two exponential factors and leaves $\xi_a^{(\rr-1)}$ on its original contour. We then deform the contour of this variable. We call a sequence of residues that successively identify variables in adjacent layers a residue chain; its endpoint is the variable that remains to be integrated. We divide each $f_\rr$-factor in \eqref{eq:translated-conditioning-phase} by $f_\rr(\sqrt L)$. Since the numerator and denominator contain equally many such factors, this normalization leaves the ratio unchanged. At a general backward step, let $q+1$ be the layer of the surviving variable in the narrow-wedge formula, where $1\leq q\leq\rr-1$; the initial step corresponds to $q=\rr-1$. The residues identify the variables in layers $q+1,\ldots,\rr$, so their increment factors telescope:
\[
\prod_{\ell=q+1}^{\rr}F_\ell(w)=\frac{f_\rr(w)}{f_q(w)}.
\]
After setting $w=\xi+\sqrt L$ and removing the variable-independent constants as above, \eqref{eq:translated-narrow-wedge-phase-ratios} gives the normalized exponential factor
\begin{equation}\label{eq:backward-left-chain-phase}
\frac{f_\rr(\xi+\sqrt L)}{f_\rr(\sqrt L)g_q(\xi)}
=\exp\left\{-\frac{1-\tau_q}{3}\xi^3-(\sqrt L+\alpha_q)\xi^2-\beta_q\xi\right\},
\end{equation}

Ignoring the linear term in the exponent of \eqref{eq:backward-left-chain-phase}, the critical points are $0$ and
\begin{equation}\label{eq:backward-left-saddle-location}
\omega_q:=-\frac{2(\sqrt L+\alpha_q)}{1-\tau_q}.
\end{equation}
The left critical point is $\omega_q+O(L^{-1/2})$. At layer $q+1$ in the narrow-wedge formula, we deform to $\omega_q+\Gamma_{\rL}$, which lies to the left of all the initial left contours for sufficiently large $L$.

The residue theorem gives
\begin{equation}\label{eq:backward-left-contour-residues}
\frac{1}{1-z_q}\int_{\Gamma_{q+1,\rL}^{\inn}}\cdots\frac{\rd\xi}{2\pi\rmi}
-\frac{z_q}{1-z_q}\int_{\Gamma_{q+1,\rL}^{\out}}\cdots\frac{\rd\xi}{2\pi\rmi}=\int_{\omega_q+\Gamma_{\rL}}\cdots\frac{\rd\xi}{2\pi\rmi}
-\frac{z_q}{1-z_q}\sum_{a=1}^{\kk_q}\Res_{\xi=\xi_a^{(q)}}(\cdots),
\end{equation}
where the dots denote the integrand at the current step. Since $\tau_q<1$, the exponential factor in \eqref{eq:backward-left-chain-phase} has cubic decay and the connecting integrals vanish. Each residue deletes a Cauchy row and column. We may therefore repeat the argument in the preceding layer, excluding labels already assigned to another residue chain. We call the labels still available for selection unused labels.

We iterate \eqref{eq:backward-left-contour-residues} from $q=\rr-1$ down to $q=1$. The integral term at step $q$ ends the chain on $\omega_q+\Gamma_{\rL}$ with exponential factor \eqref{eq:backward-left-chain-phase}. If a residue is taken at every step, the chain reaches the first layer and has exponential factor $f_\rr(\xi+\sqrt L)/f_\rr(\sqrt L)$. Its contour can then be moved to $-2\sqrt L+\Gamma_{\rL}$ without crossing a pole. We use the following bounds for these two cases.

\begin{lm}\label{lm:backward-left-chain-phase-bounds}
For all sufficiently large $L$, we have
\begin{equation}\label{eq:stopped-left-chain-phase-bound}
\int_{\omega_q+\Gamma_{\rL}}
\left|\frac{f_\rr(\xi+\sqrt L)}{f_\rr(\sqrt L)g_q(\xi)}\right|\,|\rd\xi|
\leq CL^{-1/4}\exp\left\{-\frac{4}{3(1-\tau_q)^2}L^{3/2}+CL\right\},
\end{equation}
for $1\leq q\leq\rr-1$, with $\omega_q$ given by \eqref{eq:backward-left-saddle-location}, and
\begin{equation}\label{eq:complete-left-chain-phase-bound}
\int_{-2\sqrt L+\Gamma_{\rL}}
\left|\frac{f_\rr(\xi+\sqrt L)}{f_\rr(\sqrt L)}\right|\,|\rd\xi|
\leq CL^{-1/4}\exp\left\{-\frac43L^{3/2}\right\}.
\end{equation}
The constant $C$ is independent of $L$ and $q$.
\end{lm}

\begin{proof}
For $\xi=\omega_q+\rho e^{\pm2\pi\rmi/3}$, $\rho\geq0$, direct substitution into \eqref{eq:backward-left-chain-phase} gives
\begin{equation*}
\log\left|\frac{f_\rr(\xi+\sqrt L)}{f_\rr(\sqrt L)g_q(\xi)}\right|
=-\frac{4(\sqrt L+\alpha_q)^3}{3(1-\tau_q)^2}
+\frac{2\beta_q(\sqrt L+\alpha_q)}{1-\tau_q}+\frac{\beta_q}{2}\rho
-\frac{\sqrt L+\alpha_q}{2}\rho^2
-\frac{1-\tau_q}{3}\rho^3.
\end{equation*}
For sufficiently large $L$, the terms involving $\rho$ are bounded above by $-\frac14(\sqrt L+\alpha_q)\rho^2+CL^{-1/2}$. Their exponential has integral at most $CL^{-1/4}$ over the two rays. The constant terms are $-4L^{3/2}/(3(1-\tau_q)^2)+O(L)$, which proves \eqref{eq:stopped-left-chain-phase-bound}.

For $\xi=-2\sqrt L+\rho e^{\pm2\pi\rmi/3}$, we similarly obtain
\begin{equation}\label{eq:complete-left-chain-phase-expansion}
\log\left|\frac{f_\rr(\xi+\sqrt L)}{f_\rr(\sqrt L)}\right|
=-\frac43L^{3/2}-\frac{\sqrt L}{2}\rho^2-\frac13\rho^3.
\end{equation}
Integration gives \eqref{eq:complete-left-chain-phase-bound}. The constants are uniform over the finitely many values of $q$.
\end{proof}

We repeat the deformation for all $\xi_i^{(\rr)}$, $s+1\leq i\leq\kk_\rr$, omitting labels already used. By Lemma~\ref{lm:backward-left-chain-phase-bounds} and $0<\tau_q<1$, the absolute integral of each exponential factor is at most $CL^{-1/4}e^{-4L^{3/2}/3}$ for sufficiently large $L$. These contours remain fixed in the subsequent deformations.

We next deform the unmatched $\xi_i^{(\rr+1)}$-contours, $s+1\leq i\leq\kk_{\rr+1}$, to $\Gamma_{\rL}$, following the rightward arrow in Figure~\ref{fig:unmatched-left-contour-deformation}. At layer $\ell$ in the narrow-wedge formula, $\rr+1\leq \ell\leq\mm$, the only crossed poles are at unused $\xi_a^{(\ell+1)}$ on $\Gamma_{\ell+1,\rL}^{\out}$. A residue continues the chain in layer $\ell+1$. At layer $\mm+1$, no further pole is crossed. If the chain ends in layer $p+1$ in the narrow-wedge formula, $\rr\leq p\leq\mm$, its exponential factor is $g_p(\xi)f_\rr(\sqrt L)/f_\rr(\xi+\sqrt L)$. Substituting $\xi=\rho e^{\pm2\pi\rmi/3}$ gives
\begin{equation}\label{eq:forward-left-phase-bound}
\log\left|\frac{g_p(\xi)f_\rr(\sqrt L)}{f_\rr(\xi+\sqrt L)}\right|
=-\frac{\tau_p-1}{3}\rho^3-\frac{\sqrt L+\alpha_p}{2}\rho^2-\frac{\beta_p}{2}\rho.
\end{equation}
The absolute integral of this exponential factor is at most $CL^{-1/4}$, uniformly in $p$. The connecting integrals vanish by cubic decay when $\tau_p>1$. When $p=\rr$ and $\tau_\rr=1$, the negative quadratic term gives the same bound and justifies the deformation.

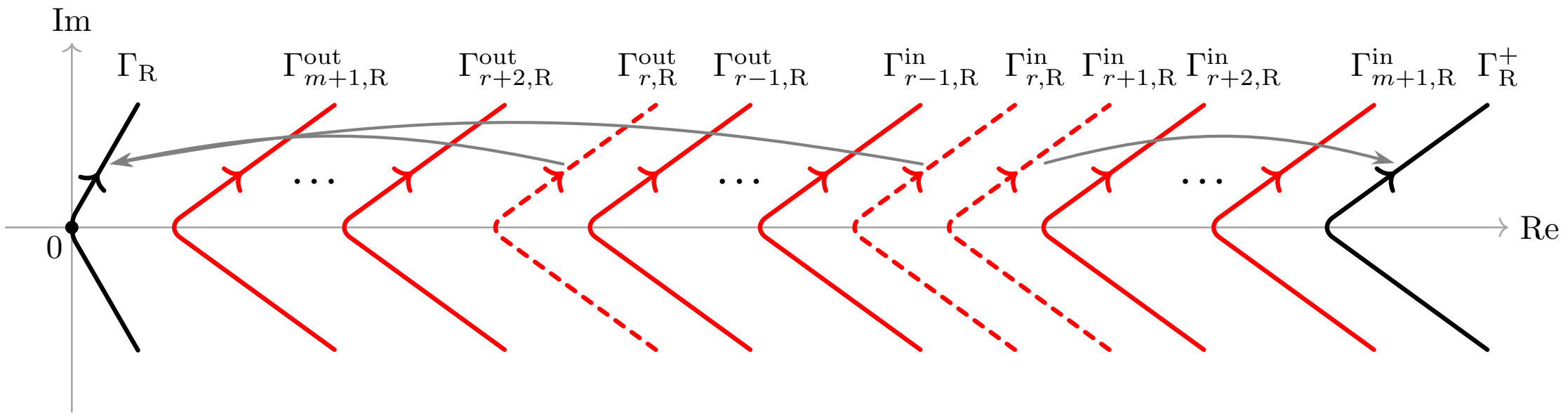
\begin{figure}[htpb]
\centering
\resizebox{0.90\textwidth}{!}{\begin{tikzpicture}[
    x=1cm,
    y=1cm,
    line cap=round,
    line join=round,
    contour/.style={
        line width=0.85pt,
        rounded corners=5pt,
        postaction={decorate},
        decoration={
            markings,
            mark=at position 0.72 with {\arrow{>}}
        }
    },
    oldcontour/.style={blue,dashed,contour},
    deform/.style={
        -{Stealth[length=2.4mm,width=1.5mm]},
        gray,
        line width=0.9pt
    },
    contour label/.style={anchor=base,inner sep=1pt}
]
    \draw[->,gray!70,line width=0.55pt]
        (-0.7,0) -- (15.2,0) node[right,black] {$\Re$};
    \draw[->,gray!70,line width=0.55pt]
        (0,-1.95) -- (0,1.95) node[above,black] {$\Im$};
    \node[below left,inner sep=2.5pt] at (0,0) {$0$};
    \fill (0,0) circle[radius=2pt];
    \begin{scope}[xshift=-0.05 cm]
            \draw[black,contour,line width=1.35pt]
                (-60:1.50) -- (0,0) -- (60:1.50);
    \end{scope}
    \begin{scope}[xshift=13.2 cm]
            \draw[black,contour,line width=1.35pt]
                (-36:2.2) -- (0,0) -- (36:2.2);
    \end{scope}
    \begin{scope}[xshift=4.4cm]
        \draw[red,dashed,contour,line width=1.35pt]
            (-36:2.2) -- (0,0) -- (36:2.2);
    \end{scope}
    \begin{scope}[xshift=8.2cm]
        \draw[red,dashed,contour,line width=1.35pt]
            (-36:2.2) -- (0,0) -- (36:2.2);
    \end{scope}
    \begin{scope}[xshift=9.2cm]
        \draw[red,dashed,contour,line width=1.35pt]
            (-36:2.2) -- (0,0) -- (36:2.2);
    \end{scope}
    \foreach \x in {1,2.8,5.4,7.2,10.2,12} {
        \begin{scope}[xshift=\x cm]
            \draw[red,contour,line width=1.35pt]
                (-36:2.2) -- (0,0) -- (36:2.2);
        \end{scope}
    }
    
    \foreach \x in {2.6,7.1,12} {
        \node[contour label] at (\x,0.4) {$\boldsymbol{\cdots}$};
    }
    \foreach \x/\y/\lab in {
        0.7/1.60/{\Gamma_\rR},
        2.8/1.60/{\Gamma^{\out}_{\mm+1,\rR}},
        4.6/1.60/{\Gamma^{\out}_{\rr+2,\rR}},
        6.1/1.60/{\Gamma^{\out}_{\rr,\rR}},
        7.3/1.60/{\Gamma^{\out}_{\rr-1,\rR}},
        9.10/1.60/{\Gamma^{\inn}_{\rr-1,\rR}},
        10.20/1.60/{\Gamma^{\inn}_{\rr,\rR}},
        11.20/1.60/{\Gamma^{\inn}_{\rr+1,\rR}},
        12.30/1.60/{\Gamma^{\inn}_{\rr+2,\rR}},
        14.1/1.60/{\Gamma^{\inn}_{\mm+1,\rR}},
        15.1/1.60/{\Gamma^{+}_{\rR}},
    } {
        \node[contour label] at (\x,\y) {$\lab$};
    }

    \draw[deform,bend right=12]
        (5.2,0.67) to (0.4,0.67);
    \draw[deform,bend right=10]
        (9,0.67) to (0.4,0.67);
    \draw[deform,bend left=15]
        (10.3,0.68) to (14,0.68);
\end{tikzpicture}
}
\caption{Deformation of the unmatched $\eta^{(\rr)}$- and
$\eta^{(\rr+1)}$-contours to $\Gamma_{\rR}$ and
$\Gamma_{\rR}^{+}$.}
\label{fig:unmatched-right-contour-deformation}
\end{figure}

We take the initial right contours to be horizontal translates of rays with arguments $\pm\pi/5$. Define
\begin{equation}\label{eq:forward-right-deformed-rays}
\Gamma_{\rR}^{+}:=\left\{\frac{\sqrt L}{2\tau_\mm}+\rho e^{\pm\pi\rmi/5}:\rho\geq0\right\}.
\end{equation}
The contour is oriented from the lower half plane to the upper half plane. As shown in Figure~\ref{fig:unmatched-right-contour-deformation}, we deform the unmatched $\eta^{(\rr)}$-contours to $\Gamma_{\rR}$ and the unmatched $\eta^{(\rr+1)}$-contours to $\Gamma_{\rR}^{+}$.

For $\eta_i^{(\rr)}$, $s'+1\leq i\leq\kk_\rr$, the deformation to $\Gamma_{\rR}$ crosses no pole when $\rr=1$. For $\rr\geq2$, the out-contour crosses no pole, while the in-contour crosses only the poles at $\eta_a^{(\rr-1)}$. We continue through the preceding layers as for the $\xi^{(\rr)}$-variables, using
\begin{equation}\label{eq:backward-right-contour-residues}
\frac{1}{1-z_q}\int_{\Gamma_{q+1,\rR}^{\inn}}\cdots\frac{\rd\eta}{2\pi\rmi}
-\frac{z_q}{1-z_q}\int_{\Gamma_{q+1,\rR}^{\out}}\cdots\frac{\rd\eta}{2\pi\rmi}=\int_{\Gamma_{\rR}}\cdots\frac{\rd\eta}{2\pi\rmi}
+\frac{1}{1-z_q}\sum_a\Res_{\eta=\eta_a^{(q)}}(\cdots),
\end{equation}
for $1\leq q\leq\rr-1$, with the sum over unused labels. Each integral term ends the chain on $\Gamma_{\rR}$. If this occurs at layer $q+1$, then for $\eta=\rho e^{\pm\pi\rmi/3}$, $\rho\geq0$,
\begin{equation}\label{eq:backward-right-phase-bound}
\log\left|\frac{g_q(\eta)f_\rr(\sqrt L)}{f_\rr(\eta+\sqrt L)}\right|
=-\frac{1-\tau_q}{3}\rho^3-\frac{\sqrt L+\alpha_q}{2}\rho^2+\frac{\beta_q}{2}\rho\leq-\frac{\sqrt L+\alpha_q}{4}\rho^2+CL^{-1/2}.
\end{equation}
For a chain reaching the first layer, the exponential factor is $f_\rr(\sqrt L)/f_\rr(\eta+\sqrt L)$, with logarithmic modulus $-\rho^3/3-\sqrt L\,\rho^2/2$. In either case, the absolute integral of the exponential factor is at most $CL^{-1/4}$. Cubic decay in the sectors between the initial and final ray directions justifies the deformation and the rotation of the tails.

We now deform the unmatched $\eta_i^{(\rr+1)}$-contours, $s'+1\leq i\leq\kk_{\rr+1}$, to $\Gamma_{\rR}^{+}$. At layer $\ell$ in the narrow-wedge formula, $\rr+1\leq \ell\leq\mm$, the only crossed poles are at unused $\eta_a^{(\ell+1)}$ on $\Gamma_{\ell+1,\rR}^{\inn}$. We continue the deformations as for the $\xi^{(\rr+1)}$-variables. For a chain ending in layer $p+1$ in the narrow-wedge formula, $\rr\leq p\leq\mm$, the parametrization \eqref{eq:forward-right-deformed-rays} gives
\begin{equation}\begin{split}\label{eq:forward-right-phase-bound}
\log\left|\frac{f_\rr(\eta+\sqrt L)}{f_\rr(\sqrt L)g_p(\eta)}\right|
&=\Re\left\{\frac{\tau_p-1}{3}\eta^3-(\sqrt L+\alpha_p)\eta^2-\beta_p\eta\right\}\\
&\leq-\frac{L^{3/2}}{8\tau_\mm^2}
-\frac{\cos(\pi/5)}{4\tau_\mm}L\rho
-\frac{\cos(2\pi/5)}{4}\sqrt L\,\rho^2
\end{split}\end{equation}
for sufficiently large $L$. To obtain the inequality, expand in $\rho$ and use $0\leq(\tau_p-1)/(2\tau_\mm)<1/2$ and $\cos(3\pi/5)<0$. Integration gives the bound $CL^{-1}e^{-L^{3/2}/(8\tau_\mm^2)}$, uniformly in $p$. The connecting integrals vanish by cubic decay when $\tau_p>1$ and by quadratic decay when $\tau_p=1$.

After the cofactor deletions, each chain endpoint occurs in the central Cauchy determinant and one outer minor, once as a row variable and once as a column variable. The endpoints on $\Gamma_{\rL}$ and $\Gamma_{\rR}$ occupy columns in the central determinant and rows in the outer minors, so their common vertex gives no pole. All contours coupled by a remaining Cauchy denominator are uniformly separated for large $L$. The separation is of order $\sqrt L$ for $\omega_q+\Gamma_{\rL}$, $-2\sqrt L+\Gamma_{\rL}$, and $\Gamma_{\rR}^{+}$. Hadamard's inequality therefore bounds the Cauchy factors uniformly in the integration variables and $L$, for fixed $\bk$.

The remaining variables have half-Brownian exponential factors on fixed contours, so the absolute integrals of these factors are bounded independently of $L$. For fixed $\bk$, there are finitely many residue terms, and their scalar factors and $z_\ell$-weights are bounded on the fixed circles. There are $\kk_\rr-s$ chains starting from $\xi^{(\rr)}$ and $\kk_{\rr+1}-s'$ starting from $\eta^{(\rr+1)}$. For chains starting from $\eta^{(\rr)}$ or $\xi^{(\rr+1)}$, the absolute integrals of the exponential factors are bounded by $CL^{-1/4}$. Combining these estimates gives
\begin{equation}\label{eq:r-fixed-order-exponential-bound}
\left|R_{(\bk;\,\rr)}(s,s')\right|
\leq C\exp\left\{-\frac43(\kk_\rr-s)L^{3/2}
-\frac{\kk_{\rr+1}-s'}{8\tau_\mm^2}L^{3/2}\right\},
\end{equation}
where $C$ may depend on $\bk$ and the parameters $(\alpha_\ell,\tau_\ell,\beta_\ell)_{\ell=1}^{\mm}$, but not on $L$. If a $\xi^{(\rr)}$-chain stops before reaching the first layer, we use the sharper bound \eqref{eq:stopped-left-chain-phase-bound}. The same argument then shows that its contribution is $o(\delta_L)$.

Since $\kk_\rr-s\geq1$, the bound \eqref{eq:r-fixed-order-exponential-bound} and the upper-tail asymptotic \eqref{eq:gue-tracy-widom-upper-tail} imply $R_{(\bk;\,\rr)}(s,s')=o(\delta_L)$ unless $s=\kk_\rr-1$ and $s'=\kk_{\rr+1}$. In that case, the admissible ranges give
\begin{equation*}
\kk_\rr-1\leq\kk_{\rr+1}\leq\kk_\rr.
\end{equation*}
These are precisely the two configurations in \eqref{eq:leading-r-configurations}.\qed

\subsubsection{Proof of Lemma~\ref{lm:leading-r-asymptotics}}

In both cases, $\kk_\rr=s+1$ and $\kk_{\rr+1}=s'$. In \eqref{eq:r-coefficient}, all $\eta^{(\rr+1)}$-variables are matched with $\eta_1^{(\rr)},\ldots,\eta_\kk^{(\rr)}$, whose contour integrals remain unchanged. Moreover, there is one unmatched left variable $\xi_{\kk_\rr}^{(\rr)}$. We apply the backward left-contour deformation from the preceding proof. For $\rr\geq2$, we first deform to $\omega_{\rr-1}+\Gamma_{\rL}$, with $\omega_{\rr-1}$ given by \eqref{eq:backward-left-saddle-location}.

\begin{figure}[htpb]
\centering
\resizebox{1\textwidth}{!}{\begin{tikzpicture}[
    x=1cm,
    y=1cm,
    line cap=round,
    line join=round,
    contour/.style={
        line width=0.85pt,
        rounded corners=5pt,
        postaction={decorate},
        decoration={
            markings,
            mark=at position 0.72 with {\arrow{>}}
        }
    },
    oldcontour/.style={blue,dashed,contour},
    deform/.style={
        -{Stealth[length=2.4mm,width=1.5mm]},
        gray,
        line width=0.9pt
    },
    contour label/.style={anchor=base,inner sep=1pt}
]
    \draw[->,gray!70,line width=0.55pt]
        (-8.85,0) -- (11.5,0) node[right,black] {$\Re$};
    \draw[->,gray!70,line width=0.55pt]
        (0,-2) -- (0,2.2) node[above,black] {$\Im$};
    \node[below left,inner sep=1.5pt] at (0,0) {$0$};
    \begin{scope}[xshift=3.5cm]
        \draw[red,dashed,contour,line width=1.35pt]
            (-36:2.2) -- (0,0) -- (36:2.2);
    \end{scope}
    \begin{scope}[xshift=7cm]
        \draw[red,dashed,contour,line width=1.35pt]
            (-36:2.2) -- (0,0) -- (36:2.2);
    \end{scope}
    \foreach \x in {1,2.5,4.5,6,8,9.5} {
        \begin{scope}[xshift=\x cm]
            \draw[red,contour,line width=1.35pt]
                (-36:2.2) -- (0,0) -- (36:2.2);
        \end{scope}
    }
    \begin{scope}[xshift=-0.05 cm]
            \draw[black,contour,line width=1.35pt]
                (-60:1.50) -- (0,0) -- (60:1.50);
    \end{scope}
    \fill (0,0) circle[radius=2pt];
    \foreach \x in {-0.6,-1.9,-2.5,-3.8,-4.4,-5.7} {
        \begin{scope}[xshift=\x cm]
            \draw[blue,contour,line width=1.35pt]
                (-120:1.50) -- (0,0) -- (120:1.50);
        \end{scope}
    }
    \begin{scope}[xshift=-7.5 cm]
            \draw[black,contour,line width=1.35pt]
                (-120:1.50) -- (0,0) -- (120:1.50);
        \end{scope}
    \foreach \x in {-1.5,-3.4,2.5,6,9.5,-5.3} {
        \node[contour label] at (\x,0.4) {$\boldsymbol{\cdots}$};
    }
    \foreach \x/\y/\lab in {
        -8.4/1.60/{-2\sqrt{L}+\Gamma_{\rL}},
        -6.6/1.60/{\Gamma^{\inn}_{\mm+1,\rL}},
        -5.30/1.60/{\Gamma^{\inn}_{\rr+2,\rL}},
        -4.4/1.60/{\Gamma^{\inn}_{\rr,\rL}},
        -3.30/1.60/{\Gamma^{\out}_{\rr,\rL}},
        -2.30/1.60/{\Gamma^{\out}_{\rr+2,\rL}},
        -1/1.60/{\Gamma^{\out}_{\mm+1,\rL}}
    } {
        \node[contour label] at (\x,\y) {$\lab$};
    }
    \foreach \x/\y/\lab in {
        0.7/1.60/{\Gamma_{\rR}},
        2.8/1.60/{\Gamma^{\out}_{\mm+1,\rR}},
        4.3/1.60/{\Gamma^{\out}_{\rr+2,\rR}},
        5.4/1.60/{\Gamma^{\out}_{\rr,\rR}},
        6.5/1.60/{\Gamma^{\out}_{\rr-1,\rR}},
        8/1.60/{\Gamma^{\inn}_{\rr-1,\rR}},
        9/1.60/{\Gamma^{\inn}_{\rr,\rR}},
        10.1/1.60/{\Gamma^{\inn}_{\rr+2,\rR}},
        11.6/1.60/{\Gamma^{\inn}_{\mm+1,\rR}},
    } {
        \node[contour label] at (\x,\y) {$\lab$};
    }

    \node[contour label] at (1.4,1.1) {$\eta_{k+1}^{(r)}$};
    \draw[deform,bend right=12]
        (4.3,0.67) to (0.4,0.67);
    \draw[deform,bend right=10]
        (7.8,0.67) to (0.4,0.67);
\end{tikzpicture}
}
\caption{Backward deformation of $\eta_{\kk+1}^{(\rr)}$ for $(\kk_\rr,\kk_{\rr+1};s,s')=(\kk+1,\kk;\kk,\kk)$. The contours are schematic: only existing layers are intended, and contours carrying no remaining variables are to be omitted.}
\label{fig:leading-backward-right-deformation}
\end{figure}
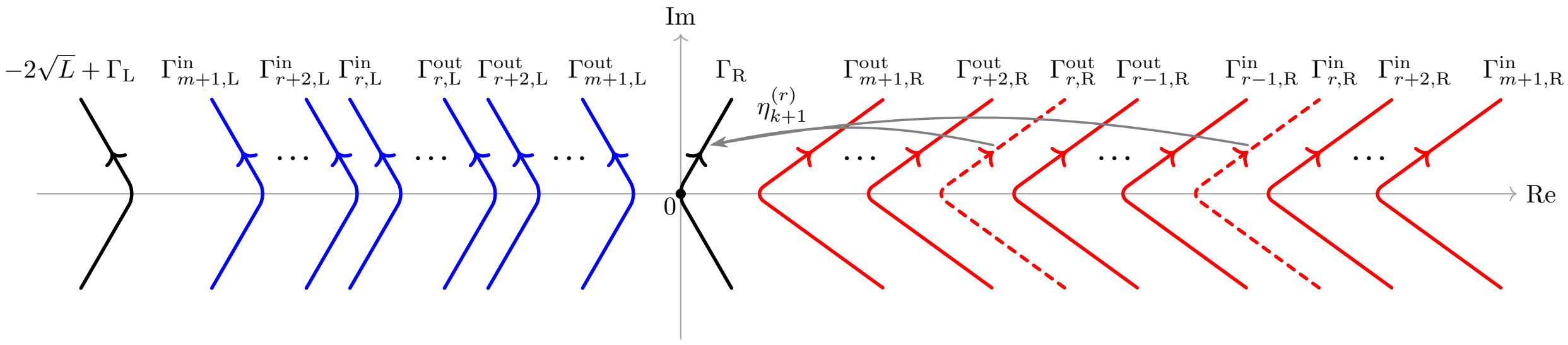

\begin{figure}[htpb]
\centering
\resizebox{1\textwidth}{!}{\begin{tikzpicture}[
    x=1cm,
    y=1cm,
    line cap=round,
    line join=round,
    contour/.style={
        line width=0.85pt,
        rounded corners=5pt,
        postaction={decorate},
        decoration={
            markings,
            mark=at position 0.72 with {\arrow{>}}
        }
    },
    oldcontour/.style={blue,dashed,contour},
    deform/.style={
        -{Stealth[length=2.4mm,width=1.5mm]},
        gray,
        line width=0.9pt
    },
    contour label/.style={anchor=base,inner sep=1pt}
]
    \draw[->,gray!70,line width=0.55pt]
        (-8.85,0) -- (11.5,0) node[right,black] {$\Re$};
    \draw[->,gray!70,line width=0.55pt]
        (1.95,-2) -- (1.95,2.2) node[above,black] {$\Im$};
    \node[below right,inner sep=1.5pt] at (2,0) {$0$};
    \begin{scope}[xshift=-1.7cm]
        \draw[blue,dashed,contour,line width=1.35pt]
            (-120:1.50) -- (0,0) -- (120:1.50);
    \end{scope}
    \foreach \x in {2.5,4.3,5.3,7.1,8.1,9.9} {
        \begin{scope}[xshift=\x cm]
            \draw[red,contour,line width=1.35pt]
                (-36:2.2) -- (0,0) -- (36:2.2);
        \end{scope}
    }
    \begin{scope}[xshift=2 cm]
            \draw[black,contour,line width=1.35pt]
                 (-120:1.50) -- (0,0) -- (120:1.50);
    \end{scope}
    \fill (1.95,0) circle[radius=2pt];
    \foreach \x in {0.5,-0.9,-2.5,-3.8,-4.4,-5.7} {
        \begin{scope}[xshift=\x cm]
            \draw[blue,contour,line width=1.35pt]
                (-120:1.50) -- (0,0) -- (120:1.50);
        \end{scope}
    }
    \begin{scope}[xshift=-7.5 cm]
            \draw[black,contour,line width=1.35pt]
                (-120:1.50) -- (0,0) -- (120:1.50);
        \end{scope}
    \foreach \x in {-0.5,-3.4,4,7,9.6,-5.3} {
        \node[contour label] at (\x,0.4) {$\boldsymbol{\cdots}$};
    }
    \foreach \x/\y/\lab in {
        -8.4/1.60/{-2\sqrt{L}+\Gamma_{\rL}},
        -6.6/1.60/{\Gamma^{\inn}_{\mm+1,\rL}},
        -5.30/1.60/{\Gamma^{\inn}_{\rr+2,\rL}},
        -4.4/1.60/{\Gamma^{\inn}_{\rr,\rL}},
        -3.40/1.60/{\Gamma^{\out}_{\rr,\rL}},
        -2.35/1.60/{\Gamma^{\out}_{\rr+1,\rL}},
        -1.2/1.60/{\Gamma^{\out}_{\rr+2,\rL}},
        0.1/1.60/{\Gamma^{\out}_{\mm+1,\rL}}
    } {
        \node[contour label] at (\x,\y) {$\lab$};
    }
    \foreach \x/\y/\lab in {
        1.3/1.60/{\Gamma_{\rL}},
        4.3/1.60/{\Gamma^{\out}_{\mm+1,\rR}},
        6.1/1.60/{\Gamma^{\out}_{\rr+2,\rR}},
        7.1/1.60/{\Gamma^{\out}_{\rr,\rR}},
        9/1.60/{\Gamma^{\inn}_{\rr,\rR}},
        10.1/1.60/{\Gamma^{\inn}_{\rr+2,\rR}},
        11.6/1.60/{\Gamma^{\inn}_{\mm+1,\rR}},
    } {
        \node[contour label] at (\x,\y) {$\lab$};
    }

     \draw[deform,bend left=20]
         (-1.92,0.63) to (1.58,0.6);
\end{tikzpicture}
}
\caption{Forward deformation of $\xi_\kk^{(\rr+1)}$ for $(\kk_\rr,\kk_{\rr+1};s,s')=(\kk,\kk;\kk-1,\kk)$. The contours are schematic: only existing layers are intended, and contours carrying no remaining variables are to be omitted.}
\label{fig:leading-forward-left-deformation}
\end{figure}
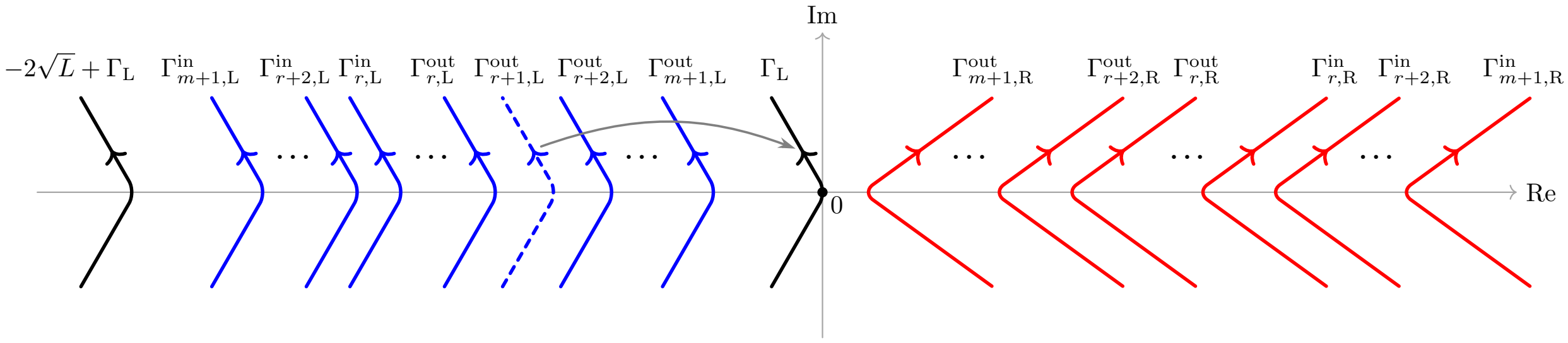

By Lemma~\ref{lm:backward-left-chain-phase-bounds}, specifically \eqref{eq:stopped-left-chain-phase-bound}, and the remaining estimates for fixed narrow-wedge series indices $\bk$, any left chain that stops before reaching the first layer contributes $o(\delta_L)$. We therefore keep only the complete left chain. For each chain, symmetry lets us choose the variable with index $\kk_\ell$ in every layer $\ell$ it visits. The surviving variable $\xi_{\kk_1}^{(1)}$ is integrated over $-2\sqrt L+\Gamma_{\rL}$, shown as the fixed leftmost black contour in Figures~\ref{fig:leading-backward-right-deformation} and~\ref{fig:leading-forward-left-deformation}.

\medskip
\noindent\emph{Proof of (i).}
Suppose $(\kk_\rr,\kk_{\rr+1};s,s')=(\kk+1,\kk;\kk,\kk)$. Since $\kk_{\rr+1}=s=\kk$, all $\xi^{(\rr+1)}$-variables are matched as well. Thus the contours $\Gamma^{\out}_{\rr+1,\rL}$ and $\Gamma^{\inn}_{\rr+1,\rR}$ are absent. The complete left chain leaves the central factor $\cC(\xi_{\kk_1}^{(1)};\eta_{\kk+1}^{(\rr)})$. We deform the $\eta_{\kk+1}^{(\rr)}$-contours to $\Gamma_{\rR}$, as shown in Figure~\ref{fig:leading-backward-right-deformation}.

The dashed contours are the original weighted contours for $\eta_{\kk+1}^{(\rr)}$. We move only these two contours. The matched variables $\xi_1^{(\rr)},\ldots,\xi_\kk^{(\rr)}$ and $\eta_1^{(\rr)},\ldots,\eta_\kk^{(\rr)}$ keep their original layer-$\rr$ contour integrals.

For $\rr\geq2$, the out-contour crosses no pole. The in-contour crosses the poles at $\eta_{\kk+1}^{(\rr)}=\eta_a^{(\rr-1)}$, $1\leq a\leq\kk_{\rr-1}$. There are no residues from later layers in the narrow-wedge formula, since $\eta_{\kk+1}^{(\rr)}$ does not occur in their Cauchy factors. As in \eqref{eq:backward-right-contour-residues}, we obtain
\begin{equation*}
\frac{1}{1-z_{\rr-1}}\int_{\Gamma_{\rr,\rR}^{\inn}}\cdots\frac{\rd\eta_{\kk+1}^{(\rr)}}{2\pi\rmi}
-\frac{z_{\rr-1}}{1-z_{\rr-1}}\int_{\Gamma_{\rr,\rR}^{\out}}\cdots\frac{\rd\eta_{\kk+1}^{(\rr)}}{2\pi\rmi}=\int_{\Gamma_{\rR}}\cdots\frac{\rd\eta_{\kk+1}^{(\rr)}}{2\pi\rmi}
+\frac{1}{1-z_{\rr-1}}\sum_{a=1}^{\kk_{\rr-1}}\Res_{\eta_{\kk+1}^{(\rr)}=\eta_a^{(\rr-1)}}(\cdots).
\end{equation*}
The two nonresidue integrals combine with coefficient $1$. The rotation from the initial ray directions $\pm\pi/5$ to those of $\Gamma_{\rR}$, namely $\pm\pi/3$, is justified by the cubic decay used in \eqref{eq:backward-right-phase-bound}. When $\rr=1$, the single contour $\Gamma_{1,\rR}$ is moved directly to $\Gamma_{\rR}$ without crossing a pole.

After each residue, we repeat the deformation with the surviving variable in the preceding layer. The resulting terms are indexed by $0\leq p\leq\rr-1$: the right chain takes residues in layers $\rr-1,\ldots,p+1$ and ends with $\eta_{\kk_{p+1}}^{(p+1)}$ on $\Gamma_{\rR}$. The case $p=\rr-1$ takes no residue and is labeled $\eta_{\kk+1}^{(\rr)}$ in the figure. The case $p=0$ reaches the first layer.

Fix $p$. Layer $\rr+1$ in the narrow-wedge formula has been removed. The matched variables from layer $\rr$ in the narrow-wedge formula form layer $\rr$ in the half-Brownian formula. Layers $\rr+2,\ldots,\mm+1$ in the narrow-wedge formula become layers $\rr+1,\ldots,\mm$ in the half-Brownian formula. Counting the remaining right variables gives the half-Brownian series indices
\begin{equation*}
\nn_\ell=
\begin{cases}
\kk_\ell,&1\leq \ell\leq p,\\
\kk_\ell-1,&p+1\leq \ell\leq\rr-1,\\
\kk,&\ell=\rr,\\
\kk_{\ell+1},&\rr+1\leq \ell\leq\mm.
\end{cases}
\end{equation*}
We relabel the remaining variables as the tuples $\widehat{\bxi}^{(\ell)}$ and $\vecteta^{(\ell)}$ from Definition~\ref{defn:q-coefficients}, preserving their order, and relabel the contours accordingly. We also relabel $z_\ell$ in the narrow-wedge formula as $z_{\ell-1}$ in the half-Brownian formula for $\ell\geq\rr+1$. This changes only the indices; the contours retain their positions, orientations, and weights.

Assume first that all $\nn_\ell$ are positive. Rescale the two chain endpoints by
\begin{equation*}
\xi_{\kk_1}^{(1)}=-2\sqrt L+L^{-1/4}\xi,
\qquad
\eta_{\kk_{p+1}}^{(p+1)}=L^{-1/4}\eta,
\qquad \xi\in\Gamma_{\rL},\quad\eta\in\Gamma_{\rR}.
\end{equation*}
The other variables remain on their fixed contours. Write the left tuples as $\bxi^{(\ell)}_{[1,\nn_\ell-1]}$ for $\ell\leq p$ and $\bxi^{(\ell)}$ for $\ell>p$, and the right tuples as $\vecteta^{(\ell)}$.

Fix $\xi,\eta$ and the remaining integration variables. We first compute the limits of the Cauchy factors containing the chain endpoints. For $p=0$, the first determinant and the central factor give
\begin{equation*}
\lim_{L\to\infty}
\frac{4L\,\cC\left(\vecteta^{(1)}\sqcup L^{-1/4}\eta;
\bxi^{(1)}\sqcup(-2\sqrt L+L^{-1/4}\xi)\right)}
{-2\sqrt L+L^{-1/4}(\xi-\eta)}=-\cC(\vecteta^{(1)};\bxi^{(1)})
\prod_{i=1}^{\nn_1}\frac{\eta_i^{(1)}}{\xi_i^{(1)}}.
\end{equation*}
For $p\geq1$, the corresponding limit is
\begin{equation}\label{eq:leading-boundary-limit-with-prefix}
\lim_{L\to\infty}
\frac{4L\,\cC\left(\vecteta^{(1)};
\bxi^{(1)}_{[1,\nn_1-1]}\sqcup(-2\sqrt L+L^{-1/4}\xi)\right)}
{-2\sqrt L+L^{-1/4}(\xi-\eta)}=-\cC(\vecteta^{(1)};\bxi^{(1)}_{[1,\nn_1-1]}\sqcup0)
\frac{\textstyle\prod_{i=1}^{\nn_1}\eta_i^{(1)}}
{\textstyle\prod_{i=1}^{\nn_1-1}\xi_i^{(1)}}.
\end{equation}

For $p\geq1$, the right endpoint also occurs in the Cauchy minor joining half-Brownian layers $p$ and $p+1$. Move its row before the $\nn_{p+1}$ remaining right rows; the sign of this permutation is included below. Then
\begin{equation}\label{eq:leading-interface-limit}
\cC\left(
\bxi^{(p)}_{[1,\nn_p-1]}\sqcup L^{-1/4}\eta\sqcup\vecteta^{(p+1)};
\vecteta^{(p)}\sqcup\bxi^{(p+1)}
\right)\longrightarrow
\cC\left(\bxi^{(p)}_{[1,\nn_p-1]}\sqcup0\sqcup\vecteta^{(p+1)};
\vecteta^{(p)}\sqcup\bxi^{(p+1)}\right).
\end{equation}
This gives the fixed zero in layer $p$. The remaining minors have the form in the Cauchy product \eqref{eq:q-cauchy-product}; for $p=0$, they are the factors in the definition \eqref{eq:q-coefficient-no-residue} of $Q^{(0)}$.

By \eqref{eq:translated-narrow-wedge-phase-ratios}, the exponential factors at the chain endpoints give
\begin{equation}\label{eq:leading-gaussian-phase}
\frac{f_\rr(-\sqrt L+L^{-1/4}\xi)}{f_\rr(\sqrt L+L^{-1/4}\eta)}
=\exp\left\{-\frac43L^{3/2}\right\}
\exp\left\{\xi^2+\eta^2+\frac{L^{-3/4}}{3}(\eta^3-\xi^3)\right\},
\end{equation}
multiplied by $g_p(L^{-1/4}\eta)\to1$ when $p\geq1$. All remaining exponential factors are $G_\ell$ or their reciprocals.

By symmetry, the label choices contribute
\begin{equation*}
\prod_{\ell=1}^{\rr-1}\kk_\ell\prod_{\ell=p+1}^{\rr-1}\kk_\ell
=\left(\prod_{\ell=1}^{p}\nn_\ell\right)
\left(\prod_{\ell=p+1}^{\rr-1}\kk_\ell^2\right).
\end{equation*}
As in the Cauchy residue calculation \eqref{eq:q-cauchy-residue}, the left and right residues at interface $\ell$ have cofactor signs $(-1)^{\kk_{\ell+1}+1}$ and $(-1)^{\kk_{\ell+1}-1}$. Including the row permutation at interface $p$, their combined sign is $1$ for $p=0$ and $(-1)^{p-1+\sum_{\ell=2}^{p}\nn_\ell}$ for $p\geq1$. The changes in the circle powers and the residue weights give
\begin{equation*}
\prod_{\ell=1}^{\rr-1}(1-z_\ell)^{\kk_\ell-\nn_\ell}(1-z_\ell^{-1})^{\kk_{\ell+1}-\nn_{\ell+1}}
\prod_{\ell=1}^{\rr-1}\frac{-z_\ell}{1-z_\ell}
\prod_{\ell=p+1}^{\rr-1}\frac{1}{1-z_\ell}=
\begin{cases}
1,&p=0,\\
\displaystyle(-1)^{p-1}\prod_{\ell=1}^{p-1}\frac{z_\ell}{1-z_\ell},&1\leq p\leq\rr-1.
\end{cases}
\end{equation*}
Thus, apart from the scalar prefactor and $\prod_{\ell=p+1}^{\rr-1}\kk_\ell^2$, the resulting coefficient is
\begin{equation*}
(-1)^{\sum_{\ell=2}^{p}\nn_\ell}
\left(\prod_{\ell=1}^{p}\nn_\ell\right)
\prod_{\ell=1}^{p-1}\frac{z_\ell}{1-z_\ell},
\end{equation*}
with the circle powers corresponding to $\bn$.

The central Cauchy factor and the far-left column each contribute $O(L^{-1/2})$, while all other entries are uniformly bounded. Thus contour separation bounds the Cauchy product by $O(L^{-1})$. After rescaling, \eqref{eq:complete-left-chain-phase-expansion} and \eqref{eq:backward-right-phase-bound}, including the first-layer case, give integrable Gaussian bounds, uniform in $L$, for the endpoint exponential factors as functions of $\xi,\eta$, after extracting the common factor $e^{-4L^{3/2}/3}$. Together with the Jacobian $L^{-1/2}$ and the bounds on the fixed exponential factors in Lemma~\ref{lm:half-brownian-absolute-convergence}, these estimates give an integrable bound after division by $\delta_L$. Dominated convergence therefore justifies taking the limits under the integrals.

With the prescribed orientations, rotation to $\rmi\realR$ gives
\begin{equation}\label{eq:leading-gaussian-integrals}
\int_{\Gamma_{\rL}}e^{\xi^2}\frac{\rd\xi}{2\pi\rmi}
=\int_{\Gamma_{\rR}}e^{\eta^2}\frac{\rd\eta}{2\pi\rmi}
=\frac{1}{2\sqrt\pi}.
\end{equation}
Combining the boundary factor $-1/(4L)$, the exponential factor and the Jacobian, we obtain
\begin{equation}\label{eq:leading-gaussian-normalization}
-\frac{L^{-1/2}}{4L}\exp\left\{-\frac43L^{3/2}\right\}
\left(\frac{1}{2\sqrt\pi}\right)^2
=-\delta_L(1+o(1)).
\end{equation}
The scalar sign is now $(-1)^\mm$. Comparing the remaining factors with \eqref{eq:q-coefficient-no-residue} and \eqref{eq:q-coefficient-with-residues}, we find that the contribution indexed by $p$ is
\begin{equation*}
\delta_L((\kk+1)!)^2
\left(\prod_{\ell=p+1}^{\rr-1}\kk_\ell^2\right)Q^{(p)}_{\bn}
+o(\delta_L).
\end{equation*}

If some half-Brownian series index $\nn_\ell$ is zero, choose the largest such $\ell$. Then $p<\ell\leq\rr-1$ and $\nn_{\ell+1}\geq1$. Both tuples in layer $\ell$ are empty, and neither endpoint occurs at interface $\ell$. We may therefore merge the layer-$(\ell+1)$ in- and out-contours. The resulting integral is independent of $z_\ell$, so the circle integral vanishes by the same cancellation as in Lemma~\ref{lm:half-brownian-empty-layer}. This agrees with $Q^{(p)}_{\bn}=0$. Summing over $0\leq p\leq\rr-1$ proves \eqref{eq:leading-r-first-limit}.

\medskip
\noindent\emph{Proof of (ii).}
Suppose $(\kk_\rr,\kk_{\rr+1};s,s')=(\kk,\kk;\kk-1,\kk)$. The only independent variable in layer $\rr+1$ in the narrow-wedge formula is $\xi_\kk^{(\rr+1)}$. The complete left chain leaves the central factor $\cC(\xi_{\kk_1}^{(1)};\xi_\kk^{(\rr+1)})$. We deform the $\xi_\kk^{(\rr+1)}$-contour $\Gamma_{\rr+1,\rL}^{\out}$ rightward to $\Gamma_{\rL}$, collecting residues at the intervening next-layer poles as described below. Figure~\ref{fig:leading-forward-left-deformation} shows this deformation.

The right-contour integrals remain unchanged. At layer $\ell$ in the narrow-wedge formula, $\rr+1\leq \ell\leq\mm$, we cross only the poles at $\xi_{\kk_\ell}^{(\ell)}=\xi_a^{(\ell+1)}$ with $\xi_a^{(\ell+1)}$ on $\Gamma_{\ell+1,\rL}^{\out}$. The integral term ends on $\Gamma_{\rL}$, while each residue continues the deformation in layer $\ell+1$. No further pole is crossed at layer $\mm+1$. The resulting endpoint is therefore $\xi_{\kk_{p+1}}^{(p+1)}\in\Gamma_{\rL}$ for some $\rr\leq p\leq\mm$, with $p=\rr$ corresponding to no forward residue.

At each crossed interface $\ell$, the label choices, the minus sign from the rightward deformation, the out-contour weight and the Cauchy cofactor give
\begin{equation*}
\kk_{\ell+1}(-1)\left(-\frac{z_\ell}{1-z_\ell}\right)(-1)^{\kk_{\ell+1}}
=(-1)^{\kk_{\ell+1}}\kk_{\ell+1}\frac{z_\ell}{1-z_\ell}.
\end{equation*}

Fix $p$ and relabel the remaining variables and circles as in part~(i). The half-Brownian series indices are
\begin{equation*}
\nn_\ell=
\begin{cases}
\kk_\ell,&1\leq \ell\leq\rr-1,\\
\kk,&\ell=\rr,\\
\kk_{\ell+1},&\rr+1\leq \ell\leq\mm.
\end{cases}
\end{equation*}
The remaining variables form the tuples $\widehat{\bxi}^{(\ell)}$ and $\vecteta^{(\ell)}$ from Definition~\ref{defn:q-coefficients}. The powers in the circle weights already correspond to the half-Brownian series indices $\bn$. The factor $\kk=\nn_\rr$ in the scalar prefactor, together with the backward and forward residue multipliers, gives
\begin{equation}\label{eq:leading-second-prefix-coefficient}
\kk\prod_{\ell=1}^{\rr-1}\left((-1)^{\nn_{\ell+1}}\nn_\ell\frac{z_\ell}{1-z_\ell}\right)
\prod_{\ell=\rr+1}^{p}\left((-1)^{\nn_\ell}\nn_\ell\frac{z_{\ell-1}}{1-z_{\ell-1}}\right)=(-1)^{\sum_{\ell=2}^{p}\nn_\ell}
\left(\prod_{\ell=1}^{p}\nn_\ell\right)
\prod_{\ell=1}^{p-1}\frac{z_\ell}{1-z_\ell}.
\end{equation}

Rescale the two endpoints by
\begin{equation*}
\xi_{\kk_1}^{(1)}=-2\sqrt L+L^{-1/4}\xi_1,
\qquad
\xi_{\kk_{p+1}}^{(p+1)}=L^{-1/4}\xi_2,
\qquad \xi_1,\xi_2\in\Gamma_{\rL}.
\end{equation*}
The boundary limit \eqref{eq:leading-boundary-limit-with-prefix} applies with $(\xi,\eta)=(\xi_1,\xi_2)$. In the Cauchy minor that becomes the factor at index $p$ in \eqref{eq:q-cauchy-product}, the forward endpoint already occupies the last left-row position, so the limit \eqref{eq:leading-interface-limit} supplies the fixed zero without a row permutation. The other minors have the form in the Cauchy product \eqref{eq:q-cauchy-product}.

By \eqref{eq:leading-gaussian-phase}, the exponential factor at the endpoints is the same as in part~(i), with $(\xi,\eta)=(\xi_1,\xi_2)$, multiplied by $g_p(L^{-1/4}\xi_2)\to1$. The remaining exponential factors are $G_\ell$ or their reciprocals. The same Cauchy bounds apply. For the second endpoint, \eqref{eq:forward-left-phase-bound} gives an integrable Gaussian bound, uniform in $L$, for its exponential factor after rescaling, also when $p=\rr$ and $\tau_\rr=1$. We may therefore apply dominated convergence as in part~(i).

Both rescaled variables lie on $\Gamma_{\rL}$. Hence \eqref{eq:leading-gaussian-integrals} gives the same normalization \eqref{eq:leading-gaussian-normalization}. Combining this with \eqref{eq:leading-second-prefix-coefficient} and comparing the remaining factors with \eqref{eq:q-coefficient-with-residues}, we obtain the contribution
\begin{equation*}
\delta_L(\kk!)^2Q^{(p)}_{\bn}+o(\delta_L).
\end{equation*}
Summing over $\rr\leq p\leq\mm$ proves \eqref{eq:leading-r-second-limit}.\qed

\subsubsection{Conditioning after the latest prescribed time}
We now consider $\rr=\mm+1$. This case requires separate treatment because there is no layer after the conditioning point, so the matching across layers $\rr$ and $\rr+1$ used in the $R$-decomposition is unavailable. Since $U^{(\mm+2)}=V^{(\mm+2)}=\emptyset$, no later layer contributes poles to the deformation of the layer-$(\mm+1)$ contours. We apply the backward-deformation arguments of Lemmas~\ref{lm:nonleading-r-coefficients} and~\ref{lm:leading-r-asymptotics} directly to $\hatD_{\bk}$.

\begin{lm}\label{lm:terminal-fixed-order-asymptotics}
Suppose $\rr=\mm+1$ and fix integers $\kk_1,\ldots,\kk_{\mm+1}\geq1$.
\begin{enumerate}[(i)]
\item If $\kk_{\mm+1}=1$, then
\begin{equation}\label{eq:terminal-single-pair-limit}
\lim_{L\to\infty}\frac{\hatD_{(\kk_1,\ldots,\kk_\mm,1)}}{\delta_L}
=\sum_{p=0}^{\mm}
\left(\prod_{\ell=p+1}^{\mm}\kk_\ell^2\right)
Q^{(p)}_{(\kk_1,\ldots,\kk_p,\kk_{p+1}-1,\ldots,\kk_\mm-1)}.
\end{equation}

\item If $\kk_{\mm+1}\geq2$, then
\begin{equation}\label{eq:terminal-multiple-pairs-negligible}
\lim_{L\to\infty}\frac{\hatD_{\bk}}{\delta_L}=0.
\end{equation}
\end{enumerate}
\end{lm}

\begin{proof}
We first translate all contours by $\sqrt L$ in \eqref{eq:narrow-wedge-integrated-coefficient} and make the change of variables \eqref{eq:narrow-wedge-shifted-variables}. Since $\tau_\mm<1$, the left-chain estimates in Lemma~\ref{lm:backward-left-chain-phase-bounds} and the right-chain estimate \eqref{eq:backward-right-phase-bound}, together with its first-layer case, apply.

\medskip
\noindent\emph{Proof of (i).}
Suppose $\kk_{\mm+1}=1$. As in Lemma~\ref{lm:leading-r-asymptotics}(i), all stopped left chains contribute $o(\delta_L)$. We keep the complete left chain, which reaches the first layer. The right chain ends in layer $p+1$ in the narrow-wedge formula for some $0\leq p\leq\mm$. The corresponding half-Brownian series indices are
\begin{equation*}
\nn_\ell=
\begin{cases}
\kk_\ell,&1\leq \ell\leq p,\\
\kk_\ell-1,&p+1\leq \ell\leq\mm.
\end{cases}
\end{equation*}
The label multiplicity is
\begin{equation*}
\prod_{\ell=1}^{\mm}\kk_\ell\prod_{\ell=p+1}^{\mm}\kk_\ell
=\left(\prod_{\ell=1}^{p}\nn_\ell\right)
\left(\prod_{\ell=p+1}^{\mm}\kk_\ell^2\right).
\end{equation*}

The narrow-wedge formula has one additional circle integral, in $z_\mm$. Combining its weight with the residue factors gives
\begin{equation*}
\oint_{>1}\frac{(1-z_\mm)^{\nn_\mm-1}}{z_\mm}\frac{\rd z_\mm}{2\pi\rmi}
=
\begin{cases}
1,&\nn_\mm\geq1,\\
0,&\nn_\mm=0.
\end{cases}
\end{equation*}
Here $\nn_\mm=\kk_\mm$ when $p=\mm$ and $\nn_\mm=\kk_\mm-1$ otherwise. If $\nn_\mm\geq1$ and some half-Brownian series index $\nn_\ell$ with $\ell<\mm$ is zero, choose the largest such $\ell$; the same circle-integral cancellation as in Lemma~\ref{lm:half-brownian-empty-layer} makes the corresponding circle integral vanish.

We compute the contribution from $\Gamma_{\mm+1,\rR}^{\out}$ explicitly. When $p=\mm$, the scalar factor $(1-z_\mm)^{\kk_\mm}(1-z_\mm^{-1})$ from the narrow-wedge formula, the left-residue weight $-z_\mm/(1-z_\mm)$, and the circle measure give $(1-z_\mm)^{\kk_\mm-1}\rd z_\mm/(2\pi\rmi z_\mm)$. The two right-contour weights then give
\begin{equation*}
\begin{aligned}
\Gamma_{\mm+1,\rR}^{\out}:\quad
&-\oint_{>1}(1-z_\mm)^{\kk_\mm-2}\frac{\rd z_\mm}{2\pi\rmi}
=\mathbf1_{\{\kk_\mm=1\}},\\
\Gamma_{\mm+1,\rR}^{\inn}:\quad
&\oint_{>1}(1-z_\mm)^{\kk_\mm-2}\frac{\rd z_\mm}{2\pi\rmi z_\mm}
=\mathbf1_{\{\kk_\mm\geq2\}}.
\end{aligned}
\end{equation*}
These contributions sum to $1$. When $\kk_\mm=1$, the entire term with no right residue comes from the outer contour. For $p<\mm$, the right chain instead starts with a residue from the inner contour.

For positive half-Brownian series indices, the identification with $Q^{(p)}_{\bn}$, Gaussian normalization, and dominated convergence follow as in the proof of Lemma~\ref{lm:leading-r-asymptotics}(i). At the terminal endpoint $p=\mm$, the fixed zero lies in the last Cauchy factor, $\cC(\widehat{\bxi}^{(\mm)}\sqcup0;\vecteta^{(\mm)})$. Including the sign $(-1)^{\mm+1}$ and the label multiplicity gives the contribution
\begin{equation*}
\delta_L\left(\prod_{\ell=p+1}^{\mm}\kk_\ell^2\right)Q^{(p)}_{\bn}
+o(\delta_L).
\end{equation*}
Summing over $0\leq p\leq\mm$ proves \eqref{eq:terminal-single-pair-limit}.

\medskip
\noindent\emph{Proof of (ii).}
Suppose $\kk_{\mm+1}\geq2$. We apply the backward deformations to all $\xi_i^{(\mm+1)}$ and $\eta_i^{(\mm+1)}$, $1\leq i\leq\kk_{\mm+1}$, excluding labels already assigned to another residue chain. Each of the $\kk_{\mm+1}$ left chains has the exponential decay in Lemma~\ref{lm:backward-left-chain-phase-bounds}, while each right-chain exponential factor has absolute integral $O(L^{-1/4})$. Combining the estimates for fixed narrow-wedge series indices $\bk$ in Lemma~\ref{lm:nonleading-r-coefficients} with \eqref{eq:gue-tracy-widom-upper-tail} gives
\begin{equation*}
\frac{|\hatD_{\bk}|}{\delta_L}
\leq\exp\left\{-\frac43(\kk_{\mm+1}-1)L^{3/2}+o(L^{3/2})\right\}
\longrightarrow0.
\end{equation*}
This proves \eqref{eq:terminal-multiple-pairs-negligible}.
\end{proof}

\subsection{Uniform bounds and convergence of the series}
\label{sec:uniform-bounds-and-series-limit}

To pass from limits with fixed narrow-wedge series indices $\bk$ to the limit of the full series, we need a summable bound for the normalized contour integrals that is uniform in $L$.

\begin{lm}\label{lm:narrow-wedge-uniform-domination}
There exist $C>0$ and $L_0\geq1$, depending only on the fixed
parameters $(\alpha_\ell,\tau_\ell,\beta_\ell)_{\ell=1}^{\mm}$, contours, and circle radii, such that, for all
$\kk_1,\ldots,\kk_{\mm+1}\geq1$,
\begin{equation*}
\sup_{L\geq L_0}
\frac{|\hatD_{\bk}|}
{\delta_L(\kk_1!\cdots\kk_{\mm+1}!)^2}
\leq
\frac{C^{\kk_1+\cdots+\kk_{\mm+1}}}
{\bigl((\kk_1+\cdots+\kk_{\mm+1})!\bigr)^{1/(\mm+1)}}.
\end{equation*}
The right-hand side is summable over all narrow-wedge series indices $\kk_1,\ldots,\kk_{\mm+1}\geq1$.
\end{lm}

\begin{proof}
We use the contour deformations in Section~\ref{sec:fixed-order-asymptotics} in a fixed order. For $\rr\leq\mm$, we first expand \eqref{eq:narrow-wedge-residue-decomposition} and interpret the factors $\binom{\kk_\rr}{s}\binom{\kk_{\rr+1}}{s}s!$ and $\binom{\kk_\rr}{s'}\binom{\kk_{\rr+1}}{s'}s'!$ in \eqref{eq:r-prefactor} as sums over labeled central matchings. For $\rr=\mm+1$, we apply the backward deformations directly to $\hatD_{\bk}$. Each subsequent residue joins two adjacent labels of the same type and deletes a Cauchy row and column. Since used labels give no further poles, the labels are partitioned into disjoint $U$- and $V$-chains, each occupying a consecutive interval of layers in the narrow-wedge formula. We count a variable not involved in a residue as a chain of length one. Each chain leaves one integration variable. If $d$ is the total number of chains, then
\begin{equation}\label{eq:residue-chain-count}
\frac{2}{\mm+1}\sum_{\ell=1}^{\mm+1}\kk_\ell
\leq d\leq 2\sum_{\ell=1}^{\mm+1}\kk_\ell.
\end{equation}
The remaining Cauchy determinant orders also sum to $d$. Indeed, their sum initially equals the number of integration variables, and each residue reduces both by one. The fixed-circle weights, in/out choices, and the last binomial factor in \eqref{eq:r-prefactor} contribute at most $C^{\kk_1+\cdots+\kk_{\mm+1}}$. Throughout the proof, $C$ may increase, but is independent of $L$, the narrow-wedge series indices $\bk$, and the residue choices.

Suppose first that $\rr\leq\mm$. The unmatched
$\xi^{(\rr)}$- and $\eta^{(\rr+1)}$-variables generate
$\kk_\rr-s$ and $\kk_{\rr+1}-s'$ chains, respectively.
The unmatched $\eta^{(\rr)}$- and $\xi^{(\rr+1)}$-variables
generate a total of $\kk_\rr+\kk_{\rr+1}-s-s'$ chains.
By Lemma~\ref{lm:backward-left-chain-phase-bounds} and
\eqref{eq:forward-left-phase-bound},
\eqref{eq:backward-right-phase-bound}, and
\eqref{eq:forward-right-phase-bound}, the absolute integrals of the exponential factors
per chain in these three groups are bounded, respectively, by
\begin{equation*}
CL^{-1/4}e^{-4L^{3/2}/3},\qquad
CL^{-1/4}e^{-L^{3/2}/(8\tau_\mm^2)},\qquad
CL^{-1/4}.
\end{equation*}
These bounds hold uniformly over the finitely many possible ending
layers and include $\tau_\rr=1$.
The remaining exponential factors are $G_\ell$ or their reciprocals
on fixed contours, and their absolute integrals are bounded by $C$.

Each chain variable remains in two Cauchy factors, once as a row variable and once as a column variable. The endpoints in the first two groups occupy rows of the central determinant and columns of the outer minors. In both occurrences, their contours are separated from the coupled contours by a fixed positive multiple of $\sqrt L$. All other coupled contours have a fixed positive separation. Only finitely many exponential factors and contour pairs occur, so the constants and $L_0$ can be chosen independently of the narrow-wedge series indices $\bk$ and residue choices. We apply Hadamard's inequality by rows to the central determinant and by columns to the outer minors, gaining $L^{-1/2}$ from each occurrence. This bounds the product of the Cauchy factors by
\begin{equation*}
C^{\kk_1+\cdots+\kk_{\mm+1}}d^{d/2}
L^{-(\kk_\rr+\kk_{\rr+1}-s-s')}.
\end{equation*}
Here we used that the determinant orders sum to $d$.
Combining this with the bounds on the exponential factors, we find that the absolute integral of each
labeled residue term, including its scalar weights, is at most
\begin{equation*}
C^{\kk_1+\cdots+\kk_{\mm+1}}d^{d/2}
L^{-\frac32(\kk_\rr+\kk_{\rr+1}-s-s')}
\exp\left\{
-\frac43(\kk_\rr-s)L^{3/2}
-\frac{\kk_{\rr+1}-s'}{8\tau_\mm^2}L^{3/2}
\right\}\leq
C^{\kk_1+\cdots+\kk_{\mm+1}}d^{d/2}\delta_L.
\end{equation*}
The last inequality follows from \eqref{eq:gue-tracy-widom-upper-tail},
since $\kk_\rr-s\geq1$ and $\kk_{\rr+1}-s'\geq0$.
The constants and $L_0$ can be chosen independently of $\bk$, $s$, and $s'$.

When $\rr=\mm+1$, there are $\kk_{\mm+1}$ backward left chains
and the same number of backward right chains. The Cauchy factor originating from $\cC(U^{(\mm+1)};V^{(\mm+1)})$
and the outer minors give the factor $L^{-\kk_{\mm+1}}$
by the same row and column estimates. Thus the absolute integral
of each labeled term is bounded by
\begin{equation*}
C^{\kk_1+\cdots+\kk_{\mm+1}}d^{d/2}
L^{-3\kk_{\mm+1}/2}e^{-\frac43\kk_{\mm+1}L^{3/2}}
\leq
C^{\kk_1+\cdots+\kk_{\mm+1}}d^{d/2}\delta_L,
\end{equation*}
using $\kk_{\mm+1}\geq1$ and the bounds on the backward-chain exponential factors as above.

We now sum over the labeled residue terms. We classify chains by their $U/V$ type and interval of layers in the narrow-wedge formula. There are at most $(\mm+1)(\mm+2)$ chain types.

Temporarily number the $d$ chains. Their types can be assigned in at most $((\mm+1)(\mm+2))^d$ ways. For each assignment compatible with the narrow-wedge series indices $\bk$, there are $(\kk_\ell!)^2$ ways to distribute the $\kk_\ell$ labels of each of $U$ and $V$ among the chains at layer $\ell$.

Each collection of labeled chains has $d!$ numberings. Since the deformation order is fixed, the chains determine every residue choice, including the central matchings. Thus division by $(\kk_1!\cdots\kk_{\mm+1}!)^2$ leaves the multiplicity bound
\begin{equation*}
\frac{((\mm+1)(\mm+2))^d}{d!}.
\end{equation*}

Consequently, for fixed $d$, the sum of the absolute
integrals, divided by $\delta_L(k_1!\cdots k_{m+1}!)^2$, is bounded by
\begin{equation*}
C^{\kk_1+\cdots+\kk_{\mm+1}}
d^{d/2}\frac{((\mm+1)(\mm+2))^d}{d!}
\leq C^{\kk_1+\cdots+\kk_{\mm+1}}d^{-d/2},
\end{equation*}
where Stirling's inequality and the upper bound in
\eqref{eq:residue-chain-count} allow the exponential factors in $d$
to be absorbed into $C^{\kk_1+\cdots+\kk_{\mm+1}}$.
Summing over $d$ and using the lower bound in
\eqref{eq:residue-chain-count}, we obtain
\begin{equation*}
\sup_{L\geq L_0}
\frac{|\hatD_{\bk}|}
{\delta_L(\kk_1!\cdots\kk_{\mm+1}!)^2}
\leq
C^{\kk_1+\cdots+\kk_{\mm+1}}
\left(\sum_{\ell=1}^{\mm+1}\kk_\ell\right)^{
-\frac{1}{\mm+1}\sum_{\ell=1}^{\mm+1}\kk_\ell}\leq
\frac{C^{\kk_1+\cdots+\kk_{\mm+1}}}
{\bigl((\kk_1+\cdots+\kk_{\mm+1})!\bigr)^{1/(\mm+1)}}.
\end{equation*}
The last bound is summable, since
\begin{equation*}
\frac{C^{\kk_1+\cdots+\kk_{\mm+1}}}
{\bigl((\kk_1+\cdots+\kk_{\mm+1})!\bigr)^{1/(\mm+1)}}
\leq
\prod_{\ell=1}^{\mm+1}
\frac{C^{\kk_\ell}}{(\kk_\ell!)^{1/(\mm+1)}},
\end{equation*}
and each factor on the right is summable over $\kk_\ell\geq1$.
\end{proof}

\bigskip

Proposition~\ref{prop:narrow-wedge-fixed-order-asymptotics}, Lemma~\ref{lm:terminal-fixed-order-asymptotics}, and Lemma~\ref{lm:narrow-wedge-uniform-domination} allow us to apply dominated convergence:
\begin{equation*}
\lim_{L\to\infty}\frac{1}{\delta_L}
\sum_{\kk_1,\ldots,\kk_{\mm+1}\geq1}
\frac{\hatD_{\bk}}{(\kk_1!\cdots\kk_{\mm+1}!)^2}
=
\sum_{\kk_1,\ldots,\kk_{\mm+1}\geq1}
\frac{1}{(\kk_1!\cdots\kk_{\mm+1}!)^2}
\lim_{L\to\infty}\frac{\hatD_{\bk}}{\delta_L}.
\end{equation*}

We next regroup the surviving terms according to the half-Brownian series indices $\bn$. Suppose first that $\rr\leq\mm$. By Proposition~\ref{prop:narrow-wedge-fixed-order-asymptotics},
only the cases $\kk_\rr=\kk_{\rr+1}+1$ and
$\kk_\rr=\kk_{\rr+1}$ contribute. Fix the half-Brownian series indices $\bn$, with $\nn_1,\ldots,\nn_\mm\geq1$. For $0\leq p\leq\rr-1$, the term $Q^{(p)}_{\bn}$ in the first case of Proposition~\ref{prop:narrow-wedge-fixed-order-asymptotics} comes from the narrow-wedge series indices
\begin{equation*}\begin{gathered}
\bk=(\nn_1,\ldots,\nn_p,
\nn_{p+1}+1,\ldots,\nn_{\rr-1}+1,
\nn_\rr+1,\nn_\rr,\nn_{\rr+1},\ldots,\nn_\mm),\\ 
\frac{((\nn_\rr+1)!)^2
\prod_{\ell=p+1}^{\rr-1}\kk_\ell^2}
{(\kk_1!\cdots\kk_{\mm+1}!)^2}
=\frac{1}{(\nn_1!\cdots\nn_\mm!)^2}.
\end{gathered}\end{equation*}
For $\rr\leq p\leq\mm$, the narrow-wedge series indices producing $Q^{(p)}_{\bn}$ in
the second case of Proposition~\ref{prop:narrow-wedge-fixed-order-asymptotics} and their normalized weight are
\begin{equation*}\begin{gathered}
\bk=(\nn_1,\ldots,\nn_{\rr-1},
\nn_\rr,\nn_\rr,\nn_{\rr+1},\ldots,\nn_\mm),\qquad
\frac{(\nn_\rr!)^2}{(\kk_1!\cdots\kk_{\mm+1}!)^2}
=\frac{1}{(\nn_1!\cdots\nn_\mm!)^2}.
\end{gathered}\end{equation*}
In these tuples, a block indexed from $a$ to $b$ is omitted when $a>b$; a product over such an index range equals $1$. The families with $0\leq p<\rr$ and $\rr\leq p\leq\mm$ give each pair $(\bn,p)$ with $\nn_\ell\geq1$ and $0\leq p\leq\mm$ exactly once.

When $\rr=\mm+1$, Lemma~\ref{lm:terminal-fixed-order-asymptotics} shows that only $\kk_{\mm+1}=1$ contributes. For each $0\leq p\leq\mm$, the term $Q^{(p)}_{\bn}$ in \eqref{eq:terminal-single-pair-limit} comes from
\begin{equation*}\begin{gathered}
\bk=(\nn_1,\ldots,\nn_p,
\nn_{p+1}+1,\ldots,\nn_\mm+1,1),\qquad 
\frac{\prod_{\ell=p+1}^{\mm}\kk_\ell^2}
{(\kk_1!\cdots\kk_{\mm+1}!)^2}
=\frac{1}{(\nn_1!\cdots\nn_\mm!)^2}.
\end{gathered}\end{equation*}
Terms for which some half-Brownian series index $\nn_\ell$ is zero vanish by Section~\ref{sec:fixed-order-asymptotics}, in agreement with Definition~\ref{defn:q-coefficients}. Lemma~\ref{lm:q-coefficient-summability} gives absolute convergence of the $Q$-series and justifies the regrouping. Applying Proposition~\ref{prop:half-brownian-residue-decomposition} and \eqref{eq:half-brownian-upper-tail}, we obtain, for every $0<\tau_1<\cdots<\tau_\mm$,
\begin{equation*}\begin{split}
\lim_{L\to\infty}\frac{1}{\delta_L}
\sum_{\kk_1,\ldots,\kk_{\mm+1}\geq1}
\frac{\hatD_{\bk}}{(\kk_1!\cdots\kk_{\mm+1}!)^2}
&=
\sum_{\nn_1,\ldots,\nn_\mm\geq1}
\frac{\textstyle\sum_{p=0}^{\mm}Q^{(p)}_{\bn}}
{(\nn_1!\cdots\nn_\mm!)^2}=
\sum_{\nn_1,\ldots,\nn_\mm\geq1}
\frac{\hatK_{\bn}}{(\nn_1!\cdots\nn_\mm!)^2}\\
&=
\rT_\mm(\beta_1,\ldots,\beta_\mm;
(\alpha_1,\tau_1),\ldots,(\alpha_\mm,\tau_\mm)).
\end{split}\end{equation*}

By Proposition~\ref{prop:half-brownian-upper-tail}, the last expression is the joint upper-tail probability of $\hKPZhB$. Together with \eqref{eq:conditional-upper-tail-series}, this proves Theorem~\ref{thm:main} at distinct positive times $\tau_1,\ldots,\tau_\mm$, which may lie before or after the conditioning time, or equal it.

It remains to allow coincident times. By Corollary~\ref{cor:half-brownian-joint-cdf-continuity}, the joint CDFs of $\hKPZhB$ are continuous in the space-time points and height parameters.

Fix $T>0$ and an arbitrary sequence $L_n\to\infty$. Apply \cite[Lemma~3.6]{Liu-Wang24} on $\realR\times(0,T)$ to the field $\{\mathcal H_{L_n}^{\mathrm{edge}}(\alpha,\tau)\}$ under the conditional law $\prob(\,\cdot\mid\hKPZ(0,1)\geq L_n)$, with limiting field $\{\hKPZhB\}$. The preceding convergence and continuity verify its hypotheses, so finite-dimensional convergence holds also at coincident times. Since $T$ and the sequence $(L_n)$ are arbitrary, this proves Theorem~\ref{thm:main}.

\def\cydot{\leavevmode\raise.4ex\hbox{.}}


\begin{thebibliography}{{Liu}22b}

\bibitem[ACH24a]{Aggarwal-Corwin-Hegde24a}
A. Aggarwal, I. Corwin, and M. Hegde.
\newblock {KPZ} fixed point convergence of the {ASEP} and stochastic six-vertex models, 2024.
\newblock arXiv:2412.18117.

\bibitem[ACH24b]{Aggarwal-Corwin-Hegde24}
A. Aggarwal, I. Corwin, and M. Hegde.
\newblock Scaling limit of the colored {ASEP} and stochastic six-vertex models, 2024.
\newblock arXiv:2403.01341.

\bibitem[ACS26]{Aggarwal-Corwin-Schmid26}
A. Aggarwal, I. Corwin, and D. Schmid.
\newblock Periodic directed landscape, 2026.
\newblock arXiv:2607.18104.

\bibitem[AFvM10]{Adler-Ferrari-van_Moerbeke10}
M. Adler, P.~L. Ferrari, and P. van Moerbeke.
\newblock Airy processes with wanderers and new universality classes.
\newblock {\em The Annals of Probability}, 38(2):714--769, 2010.

\bibitem[Bai23]{Baik23}
J. Baik.
\newblock {KPZ} limit theorems.
\newblock In D. Beliaev and S. Smirnov, editors, {\em International Congress of Mathematicians 2022, Vol. 6}, pages 4190--4211. EMS Press, 2023.

\bibitem[BBD08]{Baik-Buckingham-DiFranco08}
J. Baik, R. Buckingham, and J. DiFranco.
\newblock Asymptotics of {T}racy-{W}idom distributions and the total integral of a {P}ainlev\'e {II} function.
\newblock {\em Communications in Mathematical Physics}, 280(2):463--497, 2008.

\bibitem[BBP05]{Baik-Ben_Arous-Peche05}
J. Baik, G. Ben Arous, and S. P{\'e}ch{\'e}.
\newblock Phase transition of the largest eigenvalue for nonnull complex sample covariance matrices.
\newblock {\em The Annals of Probability}, 33(5):1643--1697, 2005.

\bibitem[BCF14]{Borodin-Corwin-Ferrari14}
A. Borodin, I. Corwin, and P.~L. Ferrari.
\newblock Free energy fluctuations for directed polymers in random media in {$1+1$} dimension.
\newblock {\em Communications on Pure and Applied Mathematics}, 67(7):1129--1214, 2014.

\bibitem[BCT25]{Baik-Cordaro-Tripathi25}
J. Baik, D. {Cordaro}, and T. {Tripathi}.
\newblock Conditional exponential directed last passage percolation under a one-point upper large deviation event, 2025.
\newblock arXiv:2508.04954.

\bibitem[BG23]{Basu-Ganguly23}
R. Basu and S. Ganguly.
\newblock Connecting eigenvalue rigidity with polymer geometry: Diffusive transversal fluctuations under large deviation.
\newblock {\em Annales de l'Institut Henri Poincar\'e, Probabilit\'es et Statistiques}, 59(2):1040--1073, 2023.

\bibitem[BL24]{Baik-Liu24}
J. Baik and Z. {Liu}.
\newblock Pinched-up periodic {KPZ} fixed point, 2024.
\newblock arXiv:2403.01624.

\bibitem[BLL26]{Baik-Liao-Liu26}
J. Baik, Y. Liao, and Z. {Liu}.
\newblock Periodic {KPZ} fixed point with general initial conditions, 2026.
\newblock arXiv:2603.01964.

\bibitem[BP08]{Borodin-Peche08}
A. Borodin and S. P{\'e}ch{\'e}.
\newblock Airy kernel with two sets of parameters in directed percolation and random matrix theory.
\newblock {\em Journal of Statistical Physics}, 132(2):275--290, 2008.

\bibitem[BR00]{Baik-Rains00}
J. Baik and E.~M. Rains.
\newblock Limiting distributions for a polynuclear growth model with external sources.
\newblock {\em Journal of Statistical Physics}, 100(3--4):523--541, 2000.

\bibitem[BT26]{Baik-Tripathi26}
J. Baik and T. {Tripathi}.
\newblock One-point fluctuations for exponential last passage percolation under upper-tail conditioning, 2026.
\newblock arXiv:2608.20552.

\bibitem[CFP10]{Corwin-Ferrari-Peche10}
I. Corwin, P.~L. Ferrari, and S. P{\'e}ch{\'e}.
\newblock Limit processes for {TASEP} with shocks and rarefaction fans.
\newblock {\em Journal of Statistical Physics}, 140(2):232--267, 2010.

\bibitem[CLW16]{Corwin-Liu-Wang16}
I. Corwin, Z. {Liu}, and D. Wang.
\newblock Fluctuations of {TASEP} and {LPP} with general initial data.
\newblock {\em The Annals of Applied Probability}, 26(4):2030--2082, 2016.

\bibitem[Cor12]{Corwin12}
I. Corwin.
\newblock The {K}ardar-{P}arisi-{Z}hang equation and universality class.
\newblock {\em Random Matrices. Theory and Applications}, 1(1):1130001, 76, 2012.

\bibitem[DDV26]{Das-Dauvergne-Virag26}
S. Das, D. Dauvergne, and B. Vir{\'a}g.
\newblock Upper tail large deviations of the directed landscape.
\newblock {\em Duke Mathematical Journal}, 175(13):2425--2493, 2026.

\bibitem[dNL17]{de_Nardis-Le_Doussal17}
J. de~Nardis and P. Le~Doussal.
\newblock Tail of the two-time height distribution for {KPZ} growth in one dimension.
\newblock {\em Journal of Statistical Mechanics: Theory and Experiment}, 2017(5):053212, 2017.

\bibitem[dNL18]{de_Nardis-Le_Doussal18}
J. de~Nardis and P. Le~Doussal.
\newblock Two-time height distribution for {1D KPZ} growth: the recent exact result and its tail via replica.
\newblock {\em Journal of Statistical Mechanics: Theory and Experiment}, 2018(9):093203, 2018.

\bibitem[DOV22]{Dauvergne-Ortmann-Virag22}
D. Dauvergne, J. Ortmann, and B. Vir{\'a}g.
\newblock {The directed landscape}.
\newblock {\em Acta Mathematica}, 229(2):201 -- 285, 2022.

\bibitem[DT24]{Das-Tsai24}
S. Das and L.-C. Tsai.
\newblock Solving marginals of the LDP for the directed landscape, 2024.
\newblock arXiv:2405.17041.

\bibitem[DV21]{Dauvergne-Virag21}
D. Dauvergne and B. Vir{\'a}g.
\newblock The scaling limit of the longest increasing subsequence, 2021.
\newblock arXiv:2104.08210.

\bibitem[DZ24]{Dauvergne-Zhang24}
D. Dauvergne and L. Zhang.
\newblock Characterization of the directed landscape from the {KPZ} fixed point, 2024.
\newblock arXiv:2412.13032.

\bibitem[DZ26]{Dauvergne-Zhang26}
D. Dauvergne and L. Zhang.
\newblock The directed landscape in half-space, 2026.
\newblock arXiv:2604.10020.

\bibitem[EJS20]{Emrah-Janjigian-Seppalainen20}
E. Emrah, C. Janjigian, and T. Sepp{\"a}l{\"a}inen.
\newblock Right-tail moderate deviations in the exponential last-passage percolation, 2020.
\newblock arXiv:2004.04285.

\bibitem[GH22]{Ganguly-Hegde22}
S. Ganguly and M. Hegde.
\newblock Sharp upper tail behavior of line ensembles via the tangent method, 2022.
\newblock arXiv:2208.08922.

\bibitem[GHZ23]{Ganguly-Hegde-Zhang23}
S. Ganguly, M. Hegde, and L. Zhang.
\newblock Brownian bridge limit of path measures in the upper tail of {KPZ} models, 2023.
\newblock arXiv:2311.12009.

\bibitem[IS04]{Imamura-Sasamoto04}
T. Imamura and T. Sasamoto.
\newblock Fluctuations of the one-dimensional polynuclear growth model with external sources.
\newblock {\em Nuclear Physics. B}, 699(3):503--544, 2004.

\bibitem[JR21]{Johansson-Rahman21}
K. Johansson and M. Rahman.
\newblock Multitime distribution in discrete polynuclear growth.
\newblock {\em Communications on Pure and Applied Mathematics}, 74(12):2561--2627, 2021.

\bibitem[JR22]{Johansson-Rahman22}
K. Johansson and M. Rahman.
\newblock On inhomogeneous polynuclear growth.
\newblock {\em The Annals of Probability}, 50(2):559--590, 2022.

\bibitem[Lia21]{Liao21}
Y. Liao.
\newblock {\em Topics in Interacting Particle Systems and Random Schr{\"o}dinger Operators}.
\newblock PhD thesis, University of Michigan, 2021.

\bibitem[{Liu}22a]{Liu22}
Z. {Liu}.
\newblock Multipoint distribution of {TASEP}.
\newblock {\em The Annals of Probability}, 50(4):1255--1321, 2022.

\bibitem[{Liu}22b]{Liu22c}
Z. {Liu}.
\newblock One-point distribution of the geodesic in directed last passage percolation.
\newblock {\em Probability Theory and Related Fields}, 184(1--2):425--491, 2022.

\bibitem[{Liu}22c]{Liu22b}
Z. {Liu}.
\newblock When the geodesic becomes rigid in the directed landscape.
\newblock {\em Electronic Communications in Probability}, 27:1 -- 13, 2022.

\bibitem[LL25]{Liao-Liu25}
Y. Liao and Z. {Liu}.
\newblock Multipoint distributions of the {KPZ} fixed point with compactly supported initial conditions, 2025.
\newblock arXiv:2509.03246.

\bibitem[LMT26]{Liu-Ma-Tripathi26}
Z. {Liu}, C. Ma, and T. {Tripathi}.
\newblock Limiting one-point fluctuations of the geodesic in the directed landscape near the endpoints when the geodesic length goes to infinity.
\newblock {\em Electronic Journal of Probability}, 31:1--27, 2026.

\bibitem[LO25]{Liu_Ortiz25}
Z. {Liu} and A. Ortiz.
\newblock On the multipoint distribution formulas of the parabolic {Airy} process, 2025.
\newblock arXiv:2507.20487.

\bibitem[LT25]{Lin-Tsai25}
Y. Lin and L.-C. Tsai.
\newblock Spacetime limit shapes of the {KPZ} equation in the upper tails.
\newblock {\em Communications in Mathematical Physics}, 406:Paper No. 113, 2025.

\bibitem[LW24]{Liu-Wang24}
Z. {Liu} and Y. Wang.
\newblock A conditional scaling limit of the {KPZ} fixed point with height tending to infinity at one location.
\newblock {\em Electronic Journal of Probability}, 29:Paper No. 30, 27, 2024.

\bibitem[LZ25]{Liu-Zhang25}
Z. {Liu} and R. Zhang.
\newblock An upper tail field of the {KPZ} fixed point.
\newblock {\em Communications in Mathematical Physics}, 406(8):Paper No. 198, 2025.

\bibitem[MQR21]{Matetski-Quastel-Remenik21}
K. Matetski, J. Quastel, and D. Remenik.
\newblock The {KPZ} fixed point.
\newblock {\em Acta Mathematica}, 227(1):115--203, 2021.

\bibitem[NQR20]{Nica-Quastel-Remenik20}
M. Nica, J. Quastel, and D. Remenik.
\newblock One-sided reflected {B}rownian motions and the {KPZ} fixed point.
\newblock {\em Forum of Mathematics. Sigma}, 8:Paper No. e63, 16, 2020.

\bibitem[NZ26]{Nissim-Zhang26}
R. Nissim and R. Zhang.
\newblock Edgeworth-type expansion for the one-point distribution of the {KPZ} fixed point with a large height at a prior location.
\newblock {\em Annales de l'Institut Henri Poincaré, Probabilités et Statistiques}, 62(2):1455--1477, 2026.

\bibitem[Qua26]{Quastel26}
J. Quastel.
\newblock The {KPZ} fixed point.
\newblock In {\em Proceedings of the International Congress of Mathematicians 2026, Volume 2: Plenary Lectures}, pages 272--301. SIAM, 2026.

\bibitem[Rah26]{Rahman26}
M. Rahman.
\newblock Revisiting the temporal law in {KPZ} random growth.
\newblock {\em Arkiv f{\"o}r Matematik}, 64(2):337--380, 2026.

\bibitem[Rem23]{Remenik23}
D. Remenik.
\newblock Integrable fluctuations in the {KPZ} universality class.
\newblock In D. Beliaev and S. Smirnov, editors, {\em International Congress of Mathematicians 2022, Vol. 6}, pages 4426--4450. EMS Press, 2023.

\bibitem[TW94]{Tracy-Widom94}
C.~A. Tracy and H. Widom.
\newblock Level-spacing distributions and the {A}iry kernel.
\newblock {\em Communications in Mathematical Physics}, 159(1):151--174, 1994.

\bibitem[Vir20]{Virag20}
B. Vir{\'a}g.
\newblock The heat and the landscape I, 2020.
\newblock arXiv:2008.07241.

\bibitem[Wu26]{Wu26}
X. Wu.
\newblock The {KPZ} equation and the directed landscape.
\newblock {\em The Annals of Probability}, 54(3):1212--1257, 2026.

\bibitem[Zha24]{Zhang24}
X. Zhang.
\newblock {TASEP} in half-space, 2024.
\newblock arXiv:2409.09974.

\bibitem[Zha25]{Zhang25}
X. Zhang.
\newblock Convergence from the log-gamma polymer to the directed landscape, 2025.
\newblock arXiv:2505.05685.

\bibitem[Zha26]{Zhang26}
X. Zhang.
\newblock Bulk directed landscape convergence and phase transitions
in half-space models, 2026.
\newblock arXiv:2602.23549.

\end{thebibliography}
\end{document}